\documentclass[10pt]{article}
\usepackage{amsgen,amsmath,amstext,amsbsy,amsopn, amssymb,amsfonts}
\usepackage{pstricks}
\usepackage{epsf}
\usepackage{float}
\usepackage{mathrsfs}
\usepackage{epsfig}

\usepackage{fullpage}
\usepackage{amsfonts}
\usepackage{amssymb}
\DeclareMathAlphabet{\mathpzc}{OT1}{pzc}{m}{it}

\newcommand{\al}{\alpha}
\newcommand{\R}{{\mathbf R}}
\newcommand{\ds}{\displaystyle}
\newcommand{\e}{\varepsilon}

\newcommand{\lra}{\longrightarrow}
\newcommand{\ra}{\rightarrow}
\newcommand{\p}{\partial}
\newcommand{\la}{\lambda}
\newcommand{\g}{\gamma}

\newcommand{\C}{{\mathbf C}}

\newcommand{\N}{{\mathbf N}}

\newcommand{\E}{{\mathbf E}}

\newcommand{\si}{\sigma}

\newcommand{\de}{\delta}

\newcommand{\Po}{{\mathcal P}}

\newcommand{\Lo}{{\mathcal L}}
\newcommand{\Do}{{\mathcal D}}
\newcommand{\Co}{{\mathcal C}}
\newcommand{\Eo}{{\mathcal E}}

\newcommand{\Io}{{\mathcal I}}

\newcommand{\wh}{\widehat}

\def \be{\begin{equation}}
\def \ee{\end{equation}}
\newcommand{\beq}{\begin{eqnarray}}
\newcommand{\eeq}{\end{eqnarray}}
\newcommand{\bea}{\begin{array}{c}}
\newcommand{\eea}{\end{array}}
\newcommand{\bi}{\begin{itemize}}
\newcommand{\ei}{\end{itemize}}

\newtheorem{Theorem}{Theorem}

\newtheorem{Lemma}[Theorem]{Lemma}

\newtheorem{rem}[Theorem]{Remark}
\newtheorem{lem}[Theorem]{Lemma}
\newtheorem{prop}[Theorem]{Proposition}
\newtheorem{cor}[Theorem]{Corollary}

\newtheorem{theo}[Theorem]{Theorem}

\def\N{\mathbb N}
\def\R{\mathbb R}

\def\C{\mathbb C}

\def\vp{\varphi}

\def\e{\varepsilon}

\def\Chi{\raise .3ex \hbox{\large $\chi$}} 
\def\vp{\varphi}
\def\s{\sigma}
\def\d{\delta}

\def\a{\alpha}

\def\g{\gamma}
\def\p{\partial}
\def\la{\lambda}

\def\G{\Gamma}

\def\carre{\ \hfill $\Box$}

\def\ii{\infty}
\def\BB{\mathcal{B}}

\def\BA{\mathcal{A}}

\def\BE{\mathcal{E}}
\def\BK{\mathcal{K}}

\def\BR{\mathcal{R}}

\def\BX{\mathcal{X}}
\def\BY{\mathcal{Y}}
\def\BD{\mathcal{D}}

\def\BC{\mathcal{C}}
\def\BO{\mathcal{O}}
\def\BW{\mathcal{W}}
\def\BU{\mathcal{U}}

\def\ra{\to}
\def\n{\|}

\def\ds{\displaystyle}
\def\bs{\boldsymbol}

\def\bi#1{\textbf{\em #1}}

\def\Fou{\mathscr{F}}
\def\bS{\mathscr{S}}
\def\cs{{\mathfrak c}_s}
\def\ex{{\sf e}}

\def\Co{\mathscr{C}}

\def\bY{\mathscr{Y}}

\title{Rarefaction pulses for the Nonlinear Schr\"odinger Equation\\ 
in the transonic limit}
\author{D. Chiron\footnote{Laboratoire J.A. Dieudonn\'e, Universit{\'e} 
de Nice-Sophia Antipolis, Parc Valrose, 06108 Nice Cedex 02, France. 
\quad \quad \quad \quad \quad \quad
{\sf e-mail}: chiron@unice.fr.} \quad \quad 
\& \quad \quad M. Mari\c{s}\footnote{Institut de Math\'ematiques de 
Toulouse, Universit\'e Paul Sabatier, 118 route de Narbonne 31062 Toulouse 
Cedex 9, France. {\sf e-mail}: mihai.maris@math.univ-toulouse.fr.}}

\date{}

\begin{document}

\maketitle

\begin{abstract}
We investigate the properties of finite energy  travelling waves to the nonlinear Schr\"odinger equation 
with nonzero conditions at infinity for a wide class of nonlinearities. In  space dimension two and three we 
prove that   travelling waves converge in the transonic limit (up to rescaling) 
 to ground states of the Kadomtsev-Petviashvili  equation. Our results  generalize an earlier result 
of F. B\'ethuel, P. Gravejat and J-C. Saut  for 
the two-dimensional Gross-Pitaevskii equation, and 
provide a rigorous proof to a conjecture by  C. Jones and P. H. Roberts about the existence of an upper branch of 
travelling waves  in dimension three.
\end{abstract}

\noindent {\bf MSC (2010)}  Main: 35C07, 35B40, 35Q55, 35Q53.  Secondary: 35B45,  35J20, 35J60, 35Q51, 35Q56, 35Q60.

\tableofcontents

\section{Introduction}

We consider the nonlinear Schr\"odinger equation in $\R^N$
\be
\tag{NLS}
i \frac{\p \Psi}{\p t} + \Delta \Psi +  F(|\Psi|^2)  \Psi  = 0 
\ee
with the condition  $|\Psi(t,x)| \to r_0$ as $ | x | \to \ii$,  
where $r_0>0$ and $F(r_0^2) =0$. This equation arises
as a relevant model in many physical situations, such as   the 
theory of Bose-Einstein condensates,  superfluidity ({see} 
\cite{Cos}, \cite{G}, \cite{IS}, \cite{JR}, \cite{JPR} and 
the surveys \cite{RB}, \cite{AHMNPTB}) or as an approximation of the Maxwell-Bloch system in Nonlinear 
Optics ({cf.} \cite{KL}, \cite{KivPeli}). When 
$F(\varrho) = 1- \varrho$, the corresponding (NLS) equation is called 
the Gross-Pitaevskii equation and is a common model for Bose-Einstein condensates. 
The so-called ``cubic-quintic'' (NLS), where
$$ F(\varrho) = - \a_1 + \a_3 \varrho - \a_5 \varrho^2 $$
 for some positive constants 
$\a_1 $, $\a_3 $ and $\a_5$
and $F $ has two positive roots, is also of high interest in Physics 
(see, {e.g.}, \cite{BP}). In Nonlinear Optics, the nonlinearity $F$ 
can take various forms (cf. \cite{KL}), for instance
\be
\label{nonlin}
F(\varrho) = - \alpha \varrho^\nu - \beta \varrho^{2\nu}, 
\quad \quad 
F(\varrho) = - \alpha 
\Big( 1 - \frac{1}{ (1+ \frac{\varrho}{\varrho_0} )^\nu} \Big), 
\quad \quad 
F(\varrho) = - \alpha \varrho \Big( 1 + \gamma \, {\rm tanh} 
( \frac{\varrho^2 - \varrho_0^2}{\sigma^2} ) \Big), \qquad \mbox{etc.,}
\ee
where $\alpha$, $\beta$, $\gamma $, $\nu$,  $\sigma > 0$ are 
given constants (the second formula, for instance, was proposed to take into account 
saturation effects). It is therefore important  to allow the nonlinearity 
to be as general as possible.

The travelling wave solutions propagating with speed $c$ in 
the $x_1$-direction are the solutions of the form $ \Psi (x,t) = U(x_1 - ct, x_2, \dots, x_N)$. The profile $U$ satisfies the  equation
\be
\tag{TW$_c$}
- i c \p_{x_1} U+ \Delta U +  F(|U|^2) U = 0.
\ee
They are supposed to play an important role in the dynamics of 
(NLS). 
Since $(U, c)$ is a solution of  (TW$_c$) 
 if and only if $(\overline{U}, -c)$ is also a solution, 
we may assume that  $c \geq 0$. The nonlinearities we consider are  general, and we will merely make use of  the following assumptions:\\

\noindent {\bf (A1)} The function $F$ is continuous on 
$[0,+\ii)$, of class $\BC^1$ near $r_0^2$,  $ F(r_0^2) = 0 $ and $ F'(r_0^2) <0$.

\medskip

\noindent
{\bf (A2)} There exist $C >0$ and $p_0 \in [1 , \frac{2}{N-2}) $ ($p_0 < \ii$ if $N=2$) such 
that $|F(\varrho)| \leq C ( 1 + \varrho^{p_0} )$ for all
$\varrho \geq 0$.

\medskip

\noindent
{\bf (A3)} There exist $C_0>0$, $\a_0>0$ and $ \varrho_0 > r_0 $ 
such that $ F(\varrho) \leq - C_0 \varrho^{\a_0}$ for  all 
$\varrho \geq \varrho_0$.\\

Assumptions (A1) and ((A2) or (A3)) are sufficient to guarantee the existence of travelling waves. 
However, in order to  get some sharp  results  we will need sometimes more information about the behavior 
of $F$ near $ r_0^2$, 
so we will  replace (A1) by 

\medskip

\noindent {\bf (A4)} The function $F$ is continuous on 
$[0,+\ii)$, of class $\BC^2$ near $r_0^2$, with 
$ F(r_0^2) = 0 , $ $ F'(r_0^2) <0 $ and
$$ F(\varrho) = F(r_0^2) + F'(r_0^2) (\varrho-r_0^2) 
+ \frac12 F''(r_0^2) ( \varrho-r_0^2)^2 + \BO((\varrho-r_0^2)^3) 
\quad \quad \mbox{ as } \quad \varrho \to r_0^2 . $$
If  $F$ is $\BC^2$ near $r_0^2$, we define, as in \cite{CM1}, 
\beq
\label{Gamma}
 \Gamma = 6 - \frac{4r_0^4}{\cs^2} F''(r_0^2) . 
\eeq

The coefficient $\G$ is positive for the Gross-Pitaevskii nonlinearity 
($F(\varrho) = 1 - \varrho$) as well as  for the cubic-quintic Schr\"odinger 
equation. However, for the nonlinearity 
$F(\varrho) = b \ex^{-\varrho/ \alpha} - a $, where $\a>0$ and $0 < a < b$ 
(which  arises  in nonlinear 
optics and takes into account saturation effects, see \cite{KL}),  
we have $ \G = 6 + 2 \ln(a/b) $,  so that $ \G$ 
can take any value in $(-\ii, 6)$, including zero. 
The coefficient $\G$ may also vanish for some polynomial nonlinearities 
(see \cite{C1d} for some examples and for the study of  travelling 
waves in dimension one in that case). In this paper  we shall be concerned only 
with the nondegenerate case $\G \not = 0$.

\medskip

{\bf Notation and function spaces.} 
For $ x = ( x_1, x_2, \dots, , x_N) \in \R^N$, we denote 
$x = ( x_1, x_{\perp})$, where  $ x_{\perp} = ( x_2, \dots, x_N) \in \R^{N-1}$.
Given a function $f$ defined on $ \R^N$, 
we denote $ \nabla _{x_{\perp}} f = ( \frac{ \p f}{\p x_2}, \dots,  \frac{ \p f}{\p x_N}).$
We will write $ \Delta_{x_{\perp} }= \frac{ \p^2}{\p x^2 } + \dots + \frac{ \p^2}{\p x^N }$. 
By "$f (t) \sim g(t)$ as $ t \to t_0$" we mean $ \lim_{t \ra t_0 } \frac{f(t)}{  g(t) } = 1$. 

\medskip

We denote by $\Fou $ the Fourier transform, 
defined by $ \Fou (f) (\xi )= \ds \int_{\R^N} \ex^{-i x.\xi } f(x) \, dx $ whenever 
$ f \in L^1 ( \R^N)$.

\medskip

Unless otherwise stated, the $L^p$ norms are computed on the  whole space $\R^N$. 

\medskip

We fix an odd function $ \chi : \R \to \R $ such that 
$\chi(s) = s $ for $0 \leq s \leq 2 r_0$, $ \chi(s) = 3 r_0$ for 
$s \geq 4 r_0 $ and $ 0 \leq \chi' \leq 1$ on $\R_+$.
As usually, we denote $ \dot{H}^1(\R^N) = \{ h \in L_{loc}^1(\R^N) \; | \; \nabla h 
\in L^2(\R^N) \}$.
We define the Ginzburg-Landau energy 
of a function $ \psi \in \dot{H}^1(\R^N) $ by 
$$
 E_{\rm GL} (\psi) =
\int_{\R^N} |\nabla \psi|^2 + (\chi^2(| \psi|)-r_0^2)^2 \ dx . 
$$
We will use the function space
$$ \BE = \left\{ \psi \in \dot{H}^{1} (\R^N) \;  \big| \; 
\chi^2 ( |\psi| ) - r_0^2 \in L^2(\R^N) \right\}
=  \left\{ \psi \in \dot{H}^{1} (\R^N) \;  \big| \;   E_{\rm GL}(\psi ) < \infty \right\} . $$
The basic properties of this space have been discussed in the Introduction of \cite{CM1}.
We will also consider the space 
$$ \BX = \left\{ u \in \BD^{1,2} (\R^N)  \;  \big| \; 
\chi^2 ( |r_0 -u| ) - r_0^2 \in L^2(\R^N) \right\},  $$
where $ \BD^{1,2} (\R^N)$ is the completion of $ \BC_c^{\infty} (\R^N)$ for the norm 
$ \| u  \|_{\BD^{1,2}} = \| \nabla u \|_{L^2(\R^N)}$. 
If $ N \geq 3$ it can be proved that 
$ \BE = \{ \al ( r_0 - u) \; \big| \; \, u \in \BX , \; \al \in \C, \;  |\al | = 1\}.$

\medskip

{\bf Hamiltonian structure. } 
The  flow associated to (NLS) formally preserves the energy
$$ E(\psi) = \int_{\R^N} |\nabla \psi|^2 + V(|\psi|^2) \ dx , $$
where $V$ is the antiderivative of $-F$ which vanishes at $r_0^2$, 
that is ${ V(s) = \int_s^{r_0^2} F(\varrho) \ d \varrho }$, 
as well as  the momentum.
The momentum (with respect to the direction of propagation $x_1$) 
is a functional $Q$ defined on $ \BE$ (or, alternatively, on $ \BX$) in the following way.
Denoting by $ \langle \cdot , \cdot \rangle$ the standard scalar product in $ \C$, 
it has been proven in  \cite{CM1} and \cite{Maris} that for any $ \psi \in \Eo $ we have 
 $ \langle i \frac{\p \psi}{\p {x_1} } , \psi \rangle \in \BY +  L^1(\R^N)  $, where 
$\BY = \{ \frac{ \p h}{ \p{x_1} } \; | \;  h \in \dot{H}^1(\R^N) \}$ and 
$ \BY$  is endowed  with the norm $ \| \p_{x_1} h \|_{\BY} =
\| \nabla h \|_{L^2(\R^N)}$. It is then possible to define 
the linear continuous functional $L$  on  $ \BY + L^1(\R^N) $ by 
$$
 L \left( \frac{ \p h}{\p{x_1}}  + \Theta \right) = \int_{\R^N} \Theta (x) \ dx 
\qquad \mbox{ for any }  \frac{ \p h}{\p{x_1}} \in \BY \mbox{ and }  \Theta \in L^1( \R^N).
$$
The momentum (with respect to the direction $x_1$) of a function $ \psi \in \BE $ is 
$
Q( \psi ) = L \left( \langle i \frac{\p \psi}{\p {x_1} } , \psi \rangle \right).
$
\\
If $ \psi \in \Eo $  does 
not vanish,  it can be lifted in the fom $\psi  = \rho \ex^{i\phi}$ and we have 
\be
\label{momentlift}
 Q(\psi) = \int_{\R^N} (r_0^2 - \rho^2) \frac{\p \phi}{\p x_1} \ dx . 
\ee

Any solution $U \in \BE$ of (TW$_c$) is a critical point of the functional $E _c = E + cQ$ 
and satisfies the standard 
Pohozaev identities (see Proposition 4.1 p. 1091 in \cite{M2}):
\be
\label{Pohozaev}
\left\{\begin{array}{ll}
P_c(U) = 0, \qquad \mbox{ where } 
\ds{ P_c(U) 
= E(U) + cQ(U)   - \frac{2}{N-1} \int_{\R^N} |\nabla_{x_\perp} U|^2 \ dx, } \qquad \mbox{ and} \\ 
 \\ 
\ds{ E(U) = 2 \int_{\R^N} |\p_{x_1} U|^2 \ dx }.
\end{array}\right.
\ee
We denote 
\be
\label{Cc}
\Co _c= \{ \psi \in \BE \; \big| \; \psi \mbox{ is not constant and } P_c( \psi ) = 0 \}.
\ee

Using the Madelung transform $ \Psi = \sqrt{ \varrho } \ex ^{ i \theta}$
(which makes sense in any domain where $ \Psi \neq 0$), equation (NLS) can be put into a hydrodynamical form.
In this context one may compute the associated speed of sound at infinity (see, for instance,  the introduction of \cite{M2}): 
$$ \cs = \sqrt{ - 2 r_0^2 F'(r_0^2) } > 0 . $$
    Under general assumptions it was  proved
that finite energy travelling waves to (NLS) with speed $c$ exist if and only if $|c| < \cs$
(see  \cite{M2, Maris}).

\medskip

Let us recall the  existence results of nontrivial 
travelling waves that  we  use.

\begin{theo} [\cite{CM1}] \label{th2dposit} Let  $N = 2$ and assume
that the nonlinearity $F$ satisfies (A2) and (A4)and that    $ \G \not = 0 $.

\medskip

(a) Suppose moreover that  $V$ is nonnegative on $[0, \infty)$. 
Then for any $q \in (-\ii , 0 )$ there exists 
$U\in \BE$ such that $Q(U)= q$ and
$$ E(U) = \inf \{ E(\psi) \; \big| \; \psi \in \BE, \; Q(\psi) = q \} . $$

(b) Without any assumption on the sign of $V$, there is $ q_{\infty} >0$ such that for any $ q \in (- q_{\infty}, 0 )$ there is 
$U\in \BE$ satisfying $Q(U)= q$ and
$$ E(U) = \inf \Big\{ E(\psi) \; \big| \; \psi \in \BE, \; Q(\psi) = q , \; \int_{\R^2} V(|\psi |^2 ) \, dx > 0 \Big\} . $$

For any  $U$ satisfying (a) or (b)  there exists $ c = c(U) \in (0,\cs)$ such 
that $U$ is a nonconstant solution to {\rm (TW$_{c(U)}$)}. 
Moreover, if $ Q( U_1 ) < Q( U_2 ) <0 $ we have 
$0 < c(U_1) < c(U_2 ) < \cs $ and
$c(U ) \to \cs $   as $q \to 0$. 
\end{theo}

\begin{theo} [\cite{CM1}] \label{th2d} Let $N = 2$. Assume that the nonlinearity $F$ 
satisfies (A2) and (A4) and  that $ \G \not = 0 $. Then  there exists 
$ 0 < k_{\infty} \leq  \infty $ such that for any $ k \in (0, k_{\infty})$, 
there is   $ \BU \in \BE $ such that 
$ \ds{\int_{\R^2} |\nabla \BU |^2 \ dx = k}$ and
$$ \int_{\R^2} V(| \BU |^2) \ d x + Q( \BU ) = 
\inf \left\{ \int_{\R^2} V(|\psi|^2) \ d x + Q(\psi) \; \Big| \;  \psi \in \BE, \;  
\int_{\R^2} |\nabla \psi|^2 \ dx = k \right\} . $$
For any such $ \BU $  there exists $c = c( \BU ) \in (0, \cs)$ 
such that the function $ U(x) = \BU ( x / c) $ is a solution to {\rm (TW$_c$)}.
Moreover, if $ \BU_1$, $\BU_2$ are as above and 
$\ds \int_{\R^2} |\nabla \BU_1|^2 \, dx  < \int_{\R^2} |\nabla \BU_2|^2 \, dx $, 
then $\cs > c(\BU_1) > c (\BU_2) > 0$ and  we have $ c( \BU ) \to \cs $ as $ k \to 0 $. 
\end{theo}

\begin{theo} [\cite{Maris}] \label{thM} Assume that $N \geq 3$ 
and  the nonlinearity $F$ satisfies (A1) and (A2). 
Then for any $0 < c < \cs $ there exists a
nonconstant $ \BU \in \BE $ such that $ P_c( \BU ) = 0 $ and 
$ E(\BU) + c Q (\BU) = \ds \inf_{ \psi \in \Co _c}  ( E( \psi )+ cQ( \psi )).$
If $N \geq 4$, any such $\BU  $ is a nontrivial solution 
to {\rm (TW$_c$)}. If $ N = 3$, for any  $\BU$ as above there exists 
$\s > 0 $ such that $U(x) = \BU ( x_1 , \s x_\perp ) \in \BE $ 
is a nontrivial solution to {\rm (TW$_c$)}.
\end{theo}

If (A3) holds it was proved that there is $C_0 >0$, depending only on $F$, such that for any $ c \in (0, \cs)$
 and for any solution $U\in \BE$ to (TW$_c$)   we have $|U| \leq C_0$ in $ \R^N$
 (see Proposition 2.2 p. 1079 in \cite{M2}).
If (A3) is satisfied but (A2) is not, one can modify $F$ in a neighborhood of infinity in such
a way  that the modified nonlinearity $ \tilde{F}$ satisfies (A2) and (A3) and $ F = \tilde{F}$ on $[0, 2 C_0]$. 
Then the solutions of  (TW$_c$) are the same as the solutions of  (TW$_c$) with  $F$ replaced by $\tilde{F}$. 
Therefore all the existence results above hold if (A2) is replaced by (A3); 
however, the minimizing properties hold only if we replace throughout $F$ and $V$ by $\tilde{F}$ and $\tilde{V}$, respectively, 
where $\tilde{V}(s)  = \ds \int_s^{ r_0 ^2} \tilde{F}( \tau ) \, d \tau$. 

\medskip

 The above results provide, under various assumptions,  travelling waves to (NLS) 
with speed close to the speed of sound $ \cs$.
We will study the behavior of travelling waves in the transonic limit $c \to \cs$ in each 
of the previous situations.

\subsection{Convergence to  ground states for (KP-I)}

In the transonic limit, the travelling waves are expected 
to be rarefaction pulses close, up to a rescaling, to  ground 
states  of the Kadomtsev-Petviashvili I (KP-I) equation. We refer to 
\cite{JR} in the case of the Gross-Pitaevskii equation ($F(\varrho)= 1 - \varrho$) 
in space dimension $N=2$ or $N=3$, and  to \cite{KL}, \cite{KAL}, 
\cite{KivPeli} in the context of Nonlinear Optics. 
In our setting, the (KP-I) equation associated to (NLS) is
\be
\tag{KP-I}
2 \p_\tau \zeta + \G \zeta \p_{z_1} \zeta 
- \frac{1}{\cs^2} \, \p_{z_1}^3 \zeta 
+ \Delta_{z_\perp} \p^{-1}_{z_1} \zeta = 0 ,
\ee
where $\Delta_{z_\perp} =  \ds{\sum_{j=2}^N \p_{z_j}^2}$ and   the coefficient $\G$ is related to the nonlinearity 
$F$ by (\ref{Gamma}).

The (KP-I) flow preserves (at least formally) the $L^2$ norm
$$ \int_{\R^N} \zeta^2 \ dz $$
and the energy
$$ \mathscr{E} (\zeta) =  \int_{\R^N} 
\frac{1}{\cs^2}\, (\p_{z_1} \zeta)^2 
+ |\nabla_{z_\perp} \p_{z_1}^{-1} \zeta|^2 
+ \frac{\G}{3} \, \zeta^3 \ dz . $$
A solitary wave of speed $ 1 / (2\cs^2) $, moving to the left in the 
$z_1$ direction,  is a particular solution of (KP-I) of the form 
$\zeta(\tau,z) = \BW(z_1 + \tau/ (2\cs^2), z_\perp)$. 
The profile  $\BW$  solves the equation 
\be
\tag{SW}
\frac{1}{\cs^2} \, \p_{z_1} \BW + \G \BW \p_{z_1} \BW 
- \frac{1}{\cs^2} \, \p_{z_1}^3 \BW + \Delta_{z_\perp} \p^{-1}_{z_1} \BW 
= 0 .
\ee
Equation (SW) has no nontrivial solution in the degenerate 
linear case $\G = 0$ or in space dimension $N \geq 4$ (see Theorem 1.1 p. 214 in \cite{dBSIHP} 
or the begining of section \ref{proofGS}).  If $\G \not = 0$,
since the nonlinearity is homogeneous,   one 
can construct solitary waves of any (positive) speed just by using the 
scaling properties of the equation.
The solutions of (SW) are critical points of the associated action
$$ \mathscr{S} (\BW) = \mathscr{E} (\BW) 
+ \frac{1}{\cs^2} \int_{\R^N} \BW^2\ dz . $$
The natural energy space for (KP-I) is  $ \mathscr{Y} (\R^N)$,  which is 
 the closure of $\p_{z_1} \BC^\ii_c(\R^N)$ 
for the (squared) norm
$$ 
\n \BW\n_{\mathscr{Y} (\R^N)}^2 =
\int_{\R^N} \frac{1}{\cs^2}\, \BW^2 + \frac{1}{\cs^2}\, (\p_{z_1} \BW)^2 
+ |\nabla_{z_\perp} \p_{z_1}^{-1} \BW|^2 \ dz . $$
From the anisotropic Sobolev embeddings (see \cite{BIN}, p. 323) it follows that 
$\mathscr{S}$ is well-defined and is a continuous functional on $\mathscr{Y}(\R^N)$ for $N=2$ and $N=3$. 
Here we are not interested in arbitrary solitary 
waves for (KP-I),  but only in {\it ground states.} A ground state 
of (KP-I) with speed $1/(2\cs^2)$ (or, equivalently, a ground state of (SW)) is 
a nontrivial solution of (SW) which minimizes the action 
$\mathscr{S}$ among all solutions of (SW). We shall 
denote $\bS_{\rm min}$ the corresponding action:
$$ \bS_{\rm min} =  \inf \Big\{ \mathscr{S} (\BW ) \; \Big| \;  \BW \in \mathscr{Y}(\R^N) \setminus \{ 0 \}, 
\ \BW \ {\rm solves \ (SW)} \Big\} . $$
The existence of ground states (with speed $1/(2\cs^2)$) for (KP-I) 
in dimensions $N=2$ and $N=3$ follows from Lemma 2.1 p. 1067 in \cite{dBSSIAM}. 
In dimension $N=2$, we may use the variational characterization 
provided by Lemma 2.2 p. 78 in \cite{dBS}:

\begin{theo} [\cite{dBS}] \label{gs2d} Assume that $N=2$ and 
$\G \not = 0$. There exists $\mu >0$ such that the set of solutions to  the 
minimization problem
\be
\label{minimiz}
{\mathscr{M}(\mu) } = \inf \left\{ \mathscr{E} (\BW)\; \Big| \;  \BW \in \mathscr{Y}(\R^2), \ 
\int_{\R^2} \frac{1}{\cs^2}\, \BW^2 \ dz = \mu \right\} ,
\ee
is precisely the set of ground states of {\rm (KP-I)} 
and it is not empty. 
Moreover, any  sequence $ (\BW _n)_{n \geq 1} \subset \mathscr{Y}(\R^2)$ 
such that $ \ds \int_{\R^2} \frac{1}{\cs^2}\, \BW_n^2 \ dz \to \mu $ and 
$\mathscr{E} (\BW_n) \to {\mathscr{M}(\mu) }$ 
contains a convergent subsequence
 in $\mathscr{Y}(\R^2)$ (up to 
translations). Finally, we have 
$$ \mu = \frac{3}{2} \, \bS_{\rm min} 
\quad \quad \quad {\it and} \quad \quad \quad 
{\mathscr{M}(\mu) } = - \frac{1}{2} \bS_{\rm min} . $$
\end{theo}

We emphasize  that  this characterization of ground states is specific to the two-dimensional case.
Indeed, since $   \mathscr{E} $ and the $L^2$ norm are conserved by (KP-I), 
it implies  the orbital stability of the set of ground states 
for (KP-I) if $N=2$ ({cf.} \cite{dBS}). On the other hand,  it is known that this set is orbitally unstable 
if $N=3$ (see \cite{Liu}). In the three-dimensional case we need the following  result, which
 shows that  ground states 
are minimizers of the action  under a Pohozaev type constraint. 
Notice that any solution of (SW) in $\mathscr{Y}(\R^N)$ 
satisfies the Pohozaev identity
$$ \int_{\R^N} \frac{1}{\cs^2}\, (\p_{z_1} \BW)^2 
+ |\nabla_{z_\perp} \p_{z_1}^{-1} \BW|^2 
+ \frac{\G}{3}\, \BW^3 + \frac{1}{\cs^2}\, \BW^2 \ dz = 
\frac{2}{N-1} \int_{\R^N} |\nabla_{z_\perp} \p_{z_1}^{-1} \BW|^2 \ dz , $$
which is (formally) obtained by multiplying (SW) by 
$z_\perp \cdot \nabla_{z_\perp} \p_{z_1}^{-1} \BW$ and integrating by 
parts (see Theorem 1.1 p. 214 in \cite{dBSIHP} for a rigorous justification). 
Taking into account how  travelling wave   solutions to (NLS)  
are constructed in Theorem \ref{thM} above, in the case $N =3$  we  consider  the minimization problem
\be
\label{miniGS}
\mathscr{S}_* = \inf \Big\{ \mathscr{S} (\BW) \; \Big| \; 
\BW \in \mathscr{Y}(\R^3) \setminus \{ 0 \}, \ \mathscr{S} (\BW) 
= \int_{\R^3} |\nabla_{z_\perp} \p_{z_1}^{-1} \BW|^2 \ dz \Big\} .
\ee

Our first result shows that in space dimension $N=3$ the ground states (with speed $1/(2\cs^2)$) 
of (KP-I) are  the solutions of  the minimization problem (\ref{miniGS}).

\begin{theo} \label{gs} 
 Assume that $ N = 3 $ and  $ \G \not = 0 $.
Then $ \mathscr{S}_* > 0 $ and the problem (\ref{miniGS}) has  minimizers. 
Moreover, $\BW_0$ is a minimizer for the problem (\ref{miniGS}) if and only if 
 there exist a ground state   $ \BW $  for {\rm (KP-I)} (with speed $1/(2\cs^2)$) and    $\s> 0$ 
such that $ \BW_0 (z) = \BW(z_1, \s  z_\perp) $.
In particular, we have $ \mathscr{S}_* = \bS_{\rm min} $.

\medskip
 
Furthermore, let  $(\BW_n)_{n \geq 1} \subset \mathscr{Y}(\R^3)$ be a sequence satisfying: 

(i) There are positive constants $ m_1, m_2 $ such that 
$ m_1 \leq \ds \int_{\R^3} \BW_n ^2 + (\p _{z_1} \BW_n )^2 \, dz \leq m_2$. 

(ii) $ \ds \int_{\R^3} \frac{ 1}{\cs ^2} \BW_n ^2 + \frac{ 1}{\cs ^2} (\p _{z_1} \BW_n )^2 + \frac{ \G}{3 } \BW_n ^3 \, dz \to 0$ as $ n \to  \infty$. 

\medskip

(iii) $ \ds \liminf_ {n \ra \infty} \mathscr{S}( \BW_n) \leq \mathscr{S}_*.$

\medskip

\noindent
Then there exist $\s>0$, $ \BW \in \mathscr{Y}(\R^3) \setminus \{ 0 \}$, 
a subsequence $ ( \BW_{n_j} )_{j \geq 0}$, and a sequence 
$(z^j)_{j\geq 0} \subset \R^3$ such that 
$z \mapsto \BW( z_1 , \s^{-1} z_\perp)$ is a ground state for 
{\rm (KP-I)} with speed $1/(2\cs^2)$ and
$$ \BW_{n_j} ( \cdot - z^j ) \to \BW \quad \quad \quad 
{\it in} \quad \mathscr{Y}(\R^3) . $$
\end{theo}

\bigskip

We will study the behavior of travelling waves to (TW$_c$) in the transonic limit $c \nearrow \cs$  
 in space dimension $N=2$ and $N=3$  under the assumption  $\G \not = 0$
(so that (KP-I) has nontrivial solitary waves). 
For $ 0 < \e < \cs$, we define $c(\e)>0$ by
$$ 
c(\e) = \sqrt{\cs^2 - \e^2} . 
$$
As already mentioned, in this asymptotic regime  the 
travelling waves are expected to be close to "the" ground state of (KP-I) 
(to the best of our knowledge, the uniqueness of this solution up to translations has not been proven yet). 
Let us give the formal derivation of this result, which follows the 
arguments given in  \cite{JR} for the 
Gross-Pitaevskii equation  
in dimensions $N=2$ and $N=3$. 
We insert the ansatz
\be
\label{ansatzKP}
U(x) = r_0 \Big( 1 + \e^2 A_\e (z) \Big) 
\exp \big( i\e \vp_\e (z) \big) \quad \quad \quad \mbox{ where } 
z_1 = \e x_1, \quad z_\perp = \e^2 x_\perp
\ee
in (TW$_{c(\e)}$), cancel the phase factor and separate the real 
and imaginary parts to obtain the system
\be
\label{MadTW}
\left\{\begin{array}{ll}
\ds{ - c(\e) \p_{z_1} A_\e 
+ 2 \e^2 \p_{z_1} \vp_\e \p_{z_1} A_\e 
+ 2 \e^4 \nabla_{z_\perp} \vp_\e \cdot \nabla_{z_\perp} A_\e 
+ (1 + \e^2 A_\e) \Big( \p_{z_1}^2 \vp_\e + 
\e^2 \Delta_{z_\perp} \vp_\e \Big)} = 0 \\ \ \\ 
\ds{ - c(\e) \p_{z_1} \vp_\e 
+ \e^2 (\p_{z_1} \vp_\e)^2 
+ \e^4 |\nabla_{z_\perp} \vp_\e |^2 
- \frac{1}{\e^2} F\Big( r_0^2 ( 1 +\e^2 A_\e )^2 \Big) 
- \e^2 \frac{ \p_{z_1}^2 A_\e + \e^2 \Delta_{z_\perp} A_\e}{1+\e^2 A_\e}} = 0.
\end{array}\right.
\ee
Formally,  if $A_\e \to A$ and $\vp_\e \to \vp$ as $\e \to 0$ 
in some reasonable sense, then to the leading order  we obtain 
$ - \cs \p_{z_1} A +\p^2_{z_1} \vp = 0 $ for the first equation 
in \eqref{MadTW}. Since $F$ is of  class $\BC^2$ near $r_0^2$,   using the Taylor expansion 
$$ F\Big( r_0^2 ( 1 +\e^2 A_\e )^2 \Big) = 
F(r_0^2) - \cs^2 \e^2 A_\e + \BO(\e^4) $$
with $F(r_0^2)=0$ and $ \cs^2 = -2 r_0^2 F'(r_0^2) $, the second equation in \eqref{MadTW} 
implies $ - \cs \p_{z_1} \vp + \cs^2 A = 0$. In both cases, we obtain 
the constraint
\be
\label{prepar}
\cs A = \p_{z_1} \vp .
\ee
We multiply  the first equation in \eqref{MadTW} by $c(\e) / \cs^2 $
and we apply the operator $ \frac{1}{ \cs ^2} \p_{z_1} $ to the second one, then we add the resulting equalities.
Using the Taylor  expansion
$$ F\Big( r_0^2 (1+ \a)^2 \Big) = - \cs^2 \a 
- \Big( \frac{\cs^2}{2} - 2 r_0^4 F''(r_0^2) \Big) \a^2 
+ F_3(\a),   \quad \quad \quad {\rm where} \ F_3(\a) = \BO(\a^3) 
\ {\rm as} \ \a \to 0 , $$
we get 
\begin{align}
\label{desing}
\frac{\cs^2 - c^2(\e) }{\e^2 \cs^2} \, \p_{z_1} A_\e & \ 
- \frac{1}{\cs^2} \, \p_{z_1} \Big( 
\frac{\p_{z_1}^2 A_\e + \e^2 \Delta_{z_\perp} A_\e}{1 + \e^2 A_\e} \Big) 
+ \frac{c(\e) }{ \cs^2} ( 1 + \e^2 A_\e ) \Delta_{z_\perp} \vp_\e 
\nonumber \\ & + 
\Big\{ 2 \frac{c(\e)}{\cs^2} \p_{z_1} \vp_\e \p_{z_1} A_\e 
+ \frac{c(\e)}{\cs^2}\, A_\e \p^2_{z_1} \vp_\e 
+ \frac{1}{\cs^2} \, \p_{z_1} \left( ( \p_{z_1} \vp_\e )^2 \right) 
+ \Big[ \frac{1}{2} - 2 r_0^4 \frac{F''(r_0^2)}{\cs^2} \Big] 
\p_{z_1} ( A_\e^2 ) \Big\} \nonumber \\ & = 
- 2 \e^2 \frac{c(\e)}{\cs^2} \nabla_{z_\perp} \vp_\e \cdot \nabla_{z_\perp} A_\e 
- \frac{\e^2}{\cs^2} \p_{z_1}  \left( \, |\nabla_{z_\perp} \vp_\e |^2 \right) 
- \frac{1}{\cs^2 \e^4}\, \p_{z_1} \left( F_3(\e^2 A_\e) \right)  .
\end{align}
If $A_\e \to A$ and $\vp_\e \to \vp$ as $\e \to 0$ 
in a suitable sense, we have   $\cs^2 - c^2(\e) = \e^2$ and  
$\p_{z_1}^{-1} A = \vp/\cs$ by \eqref{prepar}, and then 
\eqref{desing} gives
$$ \frac{1}{\cs^2}\, \p_{z_1} A - \frac{1}{\cs^2} \, \p^3_{z_1} A 
+ \G A \p_{z_1} A + \Delta_{z_\perp} \p_{z_1}^{-1} A = 0 , $$
which is (SW).\\

\medskip

The main result of this paper is as follows. 

\begin{theo}
\label{res1}
Let $N \in \{ 2, 3 \}$ and assume that the nonlinearity $F$ satisfies (A2) and (A4) with $ \G \neq 0$. 
Let $(U_n, c_n)_{n \geq 1}$ be any sequence such that $U_n \in \BE$  is a nonconstant solution of (TW$_{c_n}$), 
$ c_n \in (0, \cs )$ and $ c_n \to \cs $ as $ n \to \infty $  and one of the following situations occur:

\medskip

(a) 
$N=2$   and $U_n$ minimizes $E$ under the constraint $Q = Q(U_n)$, as in Theorem \ref{th2dposit} (a) or (b). 

\medskip

(b) $N=2$ and  $U_n (c_n \cdot) $ minimizes the functional 
$ I (\psi ) : = Q( \psi ) + \ds \int_{\R^N} V(|\psi |^2) \, dx $ 
under the constraint $ \ds \int_{ \R^N} |\nabla \psi |^2 \, dx = \int_{ \R^N} |\nabla U_n |^2 \, dx$, 
as in Theorem \ref{th2d}.

\medskip

(c) $N=3$ and $U_n$ minimizes $E_{c_n} = E + c_n Q$ under the constraint $P_{c_n} =0$, 
as in Theorem \ref{thM}. 

\medskip

Then there exists $ n _0 \in \N$ such that $ |U_n| \geq r_0 /2$ in $ \R^N$ for all $ n \geq n_0$ and, 
denoting $ \e_n = \sqrt{ \cs ^2 - c_n ^2}$ (so that $ c_n = c( \e_n)$), we have 
\beq
\label{energy}
 E(U_n) \sim - \cs Q (U_n) \sim 
r_0^2 \cs^4 (7-2N) \bS_{\rm min} \Big( \cs^2 - c_n^2 \Big)^{\frac{5-2N}{2}} 
= r_0^2 \cs^4 (7-2N) \bS_{\rm min} \e_n^{5-2N} 
\eeq
and
\beq
\label{Ec}
 E(U_n) + c_n Q (U_n) \sim \cs^2 r_0^2 \bS_{\rm min} \e_n^{7-2N} 
\qquad \mbox{ as } n \to \infty.
\eeq
Moreover, $U_n$ can be written in  the form
$$ U_n(x) = r_0 \Big( 1 + \e_n^2 A_n (z) \Big) 
\exp \big( i \e_n \vp_n (z) \big) ,  \quad \quad \quad \mbox{ where } \quad 
z_1 = \e_n x_1, \quad z_\perp = \e_n^2 x_\perp , $$
and there exist  a subsequence $(U_{n_k}, c_{n_k})_{k \geq 1}$, a ground state $\BW$ of {\rm (KP-I)} 
and a sequence $ (z^k )_{k \geq 1}  \subset \R^N$ 
such that, denoting $ \tilde{A}_k = A_{n_k}( \cdot - z^k)$, $ \tilde{\vp}_k = \vp_{n_k } (\cdot - z^k)$, 
 for any  $1 < p < \ii$ we have
 $$
 \tilde{A}_k \to \BW, \qquad 
 \p_{z_1}  \tilde{A}_k \to \p_{z_1} \BW, \qquad 
\p_{z_1} \tilde{\vp}_ k \to \cs \BW  \quad \mbox{\it and} \quad 
\p_{z_1} ^2 \tilde{\vp}_ k \to \cs \p_{z_1} \BW  \quad {\it in} \quad W^{1,p}(\R^N)   \mbox{ as } k \to \ii.
$$
\end{theo}

As already mentioned, 
if $F$ satisfies (A3) and (A4) it is possible to modify $F$ 
in a neighborhood of infinity such that the modified nonlinearity $\tilde{F}$ also satisfies (A2)
and (TW$_c$) has the same solutions as the same equation with $\tilde{F}$ instead of $F$. 
Then one may use Theorems \ref{th2dposit}, \ref{th2d} and \ref{thM} to construct travelling waves for (NLS). 
It is obvious that Theorem \ref{res1} above also applies to the solutions constructed in this way. 

\medskip

Let us mention that in the case of the Gross-Pitaevskii 
nonlinearity $F(\varrho) = 1-\varrho$ and  in dimension $N=2$,
 F. B\'ethuel, P. Gravejat and 
J-C. Saut proved  in \cite{BGS1}  the same type of convergence for 
the solutions constructed in \cite{BGS2}. Those solutions are  global 
minimizers of the energy with prescribed momentum, which allows to 
derive {\it a priori} bounds: for instance, their  energy is small. In fact, 
if $V$ is nonnegative and $N=2$, Theorem \ref{th2dposit} provides travelling 
wave solutions with speed $\simeq \cs$ for  $|q|$ small  and the proof 
of Theorem \ref{res1}  is quite similar to \cite{BGS1}, and therefore we  will focus on the 
other cases. However, if the potential $V$ achieves negative values, 
the minimization of the energy under the constraint of fixed 
momentum on the whole space $ \BE$ is no longer possible, 
hence the approach in Theorem \ref{th2d} or the local minimization approach in 
Theorem \ref{th2dposit} $(b)$. In dimension  $N=3$ (even for the 
Gross-Pitaevskii nonlinearity $F(\varrho) = 1-\varrho$), the travelling waves we deal with 
have high energy and momentum and are {\it not} minimizers of the energy at fixed momentum 
(which are the vortex rings, see \cite{BOS}). 
In particular, we have to show that the $U_n$'s are vortexless 
($|U_n|\geq r_0 /2$). For the Gross-Pitaevskii nonlinearity, 
Theorem \ref{res1} provides a rigorous proof to the existence of 
the upper branch in the so-called Jones-Roberts curve in 
dimension three (\cite{JR}). This upper branch was conjectured by 
formal expansions and numerical simulations (however limited to not 
so large momentum). In dimension $N=3$, the solutions on this 
upper branch are expected to be unstable (see \cite{BR}), and these 
rarefaction pulses should evolve by creating vortices (cf. \cite{B}).\\

It is also  natural to investigate  the one dimensional case. 
Firstly, the (KP-I) equation has to be replaced by the (KdV) 
equation
\be
\tag{KdV}
2 \p_\tau \zeta + \G \zeta \p_{z} \zeta 
- \frac{1}{\cs^2} \, \p_{z}^3 \zeta = 0 ,
\ee
and (SW) becomes
$$ \frac{1}{\cs^2}\, \p_z \BW + \G \BW \p_{z} \BW 
- \frac{1}{\cs^2} \, \p_{z}^3 \BW = 0 . $$
If  $ \G \not = 0 $,  the only 
nontrivial travelling wave for (KdV) (up to space translations) is given by
$$ {\rm w} (z) = - \frac{3}{\cs^2 \G \cosh^2(z/2)} , $$
and there holds
$$ \bS({\rm w}) = \int_{\R} \frac{1}{\cs^2}\, (\p_{z} {\rm w} )^2 
+ \frac{\G}{3} \, {\rm w}^3 \ dz + \frac{1}{\cs^2} \int_{\R} {\rm w}^2\ dz 
= \int_{\R} \frac{2}{\cs^2}\, (\p_{z} {\rm w} )^2 \ dz 
= \frac{48}{5 \cs^6 \G^2} . $$
The following 
result, which corresponds to Theorem \ref{res1} in dimension $N=1$,
was proved in \cite{C1d} by using  ODE techniques.

\begin{theo} [\cite{C1d}] \label{res2} 
Let $N = 1 $ and  assume that $F$ satisfies (A4)  with  $ \G \not = 0 $. 
Then, there are  $\d>0 $ and $ 0 <  \mathfrak{c}_0 < \cs $ with 
the following properties. For any $ \mathfrak{c}_0 \leq c < \cs $, 
there exists a solution $U_c$ to {\rm (TW$_{c}$)} satisfying 
$ \| \, |U_c| - r_0 \|_{L^\ii(\R)} \leq \d$. 
Moreover, for
$ \mathfrak{c}_0 \leq c < \cs $ any  nonconstant solution  $u$
of {\rm (TW$_{c}$)} verifying $ \| \, |u| - r_0 \|_{L^\ii(\R)} \leq \d$ is of the form  
$ u(x) = \ex^{i\theta} U_c (x-\xi)$
for some  $ \theta \in \R $ and $ \xi \in \R$.  The map $U_c$ can be written in
the form  
$$ U_c(x) = r_0 \left( 1 + \e^2 A_\e (z) \right)
\exp ( i \e \vp_\e (z) )  , \qquad \mbox{ where }    z = \e x  \quad \mbox{ and } 
\quad \e = \sqrt{\cs^2 - c^2} $$
and for any  $1 \leq p \leq \ii$,
$$ \p_{z} \vp_\e \to \cs {\rm w} 
\quad \quad \quad {\it and} \quad \quad \quad 
A_\e \to {\rm w} \quad {\it in} \quad W^{1,p}(\R) 
\quad \quad \quad {\it as} \quad \e \to 0 . $$
Finally, as $ \e \to  0 $,
$$ E(U_{c(\e)} ) \sim - \cs Q (U_{c(\e)} ) \sim 
5 r_0^2 \cs^4 \bS( {\rm w}) \Big( \cs^2 - c^2(\e) \Big)^{\frac{3}{2}} 
= \e^3 \frac{48 r_0^2 }{\cs^2 \G^2} $$
and
$$ E(U_{c(\e)} ) + c(\e) Q (U_{c(\e)} )  \sim \cs^2 r_0^2 \bS( {\rm w}) \e^{5} 
= \frac{48 r_0^2 }{5 \cs^4 \G^2} \e^5. $$
\end{theo}

\begin{rem} \rm
In the one-dimensional case it can be easily shown that 
the mapping $(\mathfrak{c}_0 , \cs) \ni c \mapsto ( A_c - r_0 , \p_z \phi ) \in W^{1,p}(\R) $, 
where $ U_c = A_c \exp( i \phi ) $, is continuous for every $1 \leq p \leq \ii$.
\end{rem}

A natural question is to investigate the dynamical counterparts of Theorems 
\ref{res1} and \ref{res2}. If $\Psi_\e^0$ is an initial datum for (NLS) 
of the type
$$ \Psi_\e^0 (x) = r_0 \Big( 1 + \e^2 A_\e^0 (z) \Big) 
\exp \Big( i \e \vp_\e^0(z) \Big) , $$
with $z=(z_1,z_\perp) = (\e x_1, \e^2 x_\perp)$ and 
$ \cs A_\e^0 \simeq \p_{z_1} \vp_\e^0$, we use for $\Psi_\e$ 
the ansatz at time $t>0$, for some functions 
$A_\e$, $\vp_\e$ depending on $(\tau,z)$,
$$ \Psi_\e (t,x) = r_0 \Big( 1 + \e^2 A_\e (\tau,z)\Big) 
\ex^{i \e \vp_\e(\tau,z) } , \quad \quad \quad 
\tau = \cs \e^3 t , \quad z_1 = \e ( x_1 - \cs t) , 
\quad z_\perp = \e^2 x_\perp  . $$
Similar computations  imply  that, for times $\tau$ of order one 
(that is $t$ of order $\e^{-3}$), we have  $ \cs A_\e \simeq \p_{z_1} \vp_\e$ 
and $A_\e$ converges to a solution of the (KP-I) equation.
This (KP-I) asymptotic dynamics for 
the Gross-Pitaevskii equation 
in dimension $N=3$ is formally derived in \cite{BR} and is used to investigate 
the linear instability 
of the solitary waves of speed close to $ \cs = \sqrt{2} $. The one-dimensional 
analogue, where the (KP-I) equation has to be replaced by the corresponding 
Korteweg-de Vries equation, can be found in \cite{ZK} and \cite{KAL}. 
The rigorous mathematical proofs of these regimes have been provided 
in \cite{CR2} in arbitrary space  dimension and for a general nonlinearity $F$ 
(the coefficient $\G$ might even vanish), respectively in \cite{BGSS} for the one 
dimensional Gross-Pitaevskii equation by using 
 the complete integrability of the equation (more precisely,  the existence of 
sufficiently many conservation laws).

\subsection{Scheme of the proof of Theorem \ref{res1}}

In case $(a)$  there is a direct proof of Theorem \ref{res1} which  is quite  similar 
to the one in \cite{BGS1}.  Moreover, it follows  from 
 Proposition 5.12 in \cite{CM1}  that if $(U_n, c_n)$ satisfies $(a)$
then it also satisfies $(b)$, so it suffices to prove Theorem \ref{res1} 
in cases $(b)$ and $(c)$. 

\medskip

The first step is to give sharp asymptotics for the quantities 
minimized in \cite{CM1} and \cite{Maris}  in order to prove the existence 
of travelling waves, namely to estimate 
$$ I_{\rm min} (k) =
\inf \Big\{ \int_{\R^2} V(|\psi|^2) \ d x + Q(\psi )  \; \big\vert \;    \psi \in \BE, \ 
\int_{\R^2} |\nabla \psi|^2 \ dx = k \Big\}  \qquad \mbox{ as } k \to 0$$
and
$$ T_c = \inf \Big\{ E(\psi ) + c Q (\psi ) \; \big\vert \;  \psi \in \BE, \; \psi \mbox{ is not constant, }
  E(\psi) + c Q (\psi) =
\int_{\R^3} | \nabla_{x_\perp} \psi |^2 \ dx\} \qquad \mbox{ as } c \to \cs.$$
These bounds are obtained by plugging test functions  with the 
ansatz \eqref{ansatzKP} into the corresponding  minimization problems, 
where  $(A_\e,\vp_\e ) \simeq (A, \cs^{-1} \p_{z_1}^{-1} A )$ 
and $A$ is a ground state for (KP-I).
A similar upper bound for $I_{\rm min}(k)$  was already a 
crucial point in \cite{CM1} to rule out the dichotomy of
minimizing sequences.

\begin{prop}
\label{asympto}
Assume that $F$ satisfies (A2)  and (A4) with  $ \G \not = 0 $. Then: 

\medskip

(i) If $N=2$, we have  as $ k \to 0 $
$$ I_{\rm min} (k) \leq - \frac{k}{\cs^2} 
- \frac{4 k^3}{27 r_0^4 \cs^{12} \bS^2_{\rm min}} + \BO(k^5) . $$

(ii) If $ N = 3 $, the following upper bound holds  as 
$\e \to 0$ (that is, as $c(\e) \to \cs$):
$$ T_{c(\e)} \leq 
\cs^2 r_0^2 \bS_{\rm min} ( \cs^2 - c^2(\e) )^{\frac{1}{2}} 
+ \BO\Big( ( \cs^2 - c^2(\e) )^{\frac{3}{2}} \Big) 
= \cs^2 r_0^2 \bS_{\rm min} \e + \BO(\e^3) . $$
\end{prop}

The second step is to derive upper bounds for the energy and the momentum. 
In space dimension three (case $(c)$) this is   tricky. 
Indeed, if $U_c$ is a minimizer of $E_c$ under the constraint $ P_c = 0$, the only 
information we have is about $ T_c = \ds \int_{\R^N} |\nabla _{x_{\perp}} U_c |^2 \, dx $ 
(see the first identity in (\ref{Pohozaev})). 
In particular, we have no {\it a priori } bounds on 
$\ds \int_{\R^N} \Big| \frac{ \p U_c}{\p x_1 } \Big|^2 \, dx$, $Q(U_c)$ and the potential energy
$\ds \int_{\R^N} V(|U_c|^2) \, dx$. 
Using an averaging argument we infer that there is a sequence $(U_n, c_n)$ for which we have "good" bounds 
on the energy and the momentum. 
Then we prove a rigidity property of "good sequences": 
 any  sequence $(U_n, c_n)$  that satisfies the "good bounds" 
has a subsequence that satisfies the conclusion 
of Theorem \ref{res1}. 
This rigid behavior combined with the existence of a sequence with "good bounds"
 and a continuation argument allow us to conclude that 
  Theorem \ref{res1} holds for {\it any } sequence $(U_n, c_n)$  with $ c_n \to \cs$ (as in (c)).
More precisely, we will prove:

\begin{prop} \label{monoto} 
Let  $N \geq 3$ and assume that $F$ satisfies (A1) and (A2). Then: 

\medskip

(i) For any $ c \in (0, \cs )$ and any  minimizer $ U  $ of $E_c$ in $ \Co _c $  we have $Q(U) <0$. 

\medskip

(ii) The function $ (0, \cs) \ni c \longmapsto T_c \in \R_+$ is  decreasing, 
 thus  has a derivative almost everywhere.

\medskip
 
(iii) The function $  c \longmapsto T_c $ is left continuous on $(0, \cs)$.
 If it 
has a derivative at $c_0$, then for any minimizer $U_0$ of  $ E_{c_0}$ under the constraint $P_{c_0} = 0$,  
scaled so that $U_0$ solves {\rm (TW}$_{c_0}${\rm)}, there holds
$$ \frac{d T_c}{dc}_{|c=c_0} = Q(U_0) . $$

(iv) Let $ c_0 \in ( 0 , \cs )$. 
Assume that there is a sequence $ (c_n) _{n \geq 1}$ such that $ c_n > c_0$, $ c_n \to c_0$ 
and for any $n$ there is a minimizer $ U_n \in \BE$ 
of $E_{c_n}$ on $ \Co _{c_n}$ which solves (TW$_{c_n}$) and the sequence $(Q(U_n))_{n \geq 1 }$ is bounded. 
Then $ c \longmapsto T_c$ is continuous at $ c_0$.

(v) Let $0 < c_1 <c_2 < \cs$. Let $U_i $ be minimizers of $E_{c_i} $ 
on $ \Co_{c_i}$, $ i =1,2$, 
such that $ U_i$ solves {\rm (TW}$_{c_i}${\rm)}.
Denote $ q_1 = Q(U_1)$ and $q_2 = Q(U_2).$ Then we have
$$
\frac{ T_{c_1}^2}{q_1^2} - c_1 ^2 \geq \frac{ T_{c_2}^2}{q_2^2} - c_2 ^2 .
$$

(vi) If $N=3$, $F$ verifies (A4)  and $\G \not = 0$, there exist a constant $C>0$ and 
a sequence  $\e_n \to 0$ such that for any minimizer $U_n \in \BE$ of $ E_{c(\e_n)} $ on $ \Co_{c(\e_n)}$ 
which solves {\rm (TW}$_{c(\e_n)}${\rm)} we have
$$  E(U_n) \leq \frac{C}{\e_n} \qquad \mbox{ and } \qquad |Q(U_n) | \leq \frac{C}{\e_n} . $$
\end{prop}

\begin{prop} \label{convergence} 
Assume that $N=3$, (A2) and (A4) hold and $ \G \neq 0$. 
Let $( U_n, \e _n)_{n \geq 1}$  be a sequence such that $ \e_n \to 0$, 
$U_n$ minimizes $E_{c(\e_n)}$ on $ \Co _{c(\e_n)}$,  satisfies {\rm (TW}$_{c(\e_n)}${\rm)}  
and there exists a constant $C>0$ such that 
$$
E(U_n)  \leq \frac{C}{\e_n} \qquad \mbox{ and } \qquad |Q(U_n) | \leq \frac{C}{\e_n} 
\qquad \mbox{ for all } n.
$$
Then there is a subsequence of $(U_n, c( \e _n))_{n \geq 1}$ which satisfies the conclusion of Theorem \ref{res1}.
\end{prop}

\begin{prop} \label{global3}
Let $ N =3$ and suppose that (A2) and (A4) hold with $ \G \neq 0$. 
There are $ K > 0$ and $ \e _* > 0 $ such that for any $ \e \in ( 0, \e _* )$ and for any 
minimizer $U$ of $E_{c(\e)} $ on $ \Co_{c (\e)}$ scaled so that $U$ satisfies {\rm (TW}$_{c(\e )}${\rm)}
we have
$$
E(U) \leq \frac{K}{\e} \qquad \mbox{ and } \qquad |Q(U) | \leq \frac{K}{\e}. 
$$
\end{prop}

It is now obvious that the proof of Theorem \ref{res1} in the three-dimensional 
case follows directly from Propositions \ref{convergence} and \ref{global3} above. 

\medskip

The most difficult and technical point in the above program is to prove Proposition \ref{convergence}. 
Let us describe our strategy to carry out  that proof, as well as the proof of Theorem \ref{res1} in the two-dimensional case. 

\medskip

Once we have a sequence of travelling waves to (NLS) with "good bounds" on the energy and the momentum 
and speeds that tend to $ \cs$, we need to show that those solutions do not vanish and can be lifted. 
We recall the following result, which is a consequence of Lemma 7.1 in \cite{CM1}:

\begin{lem} [\cite{CM1}] \label{liftingfacile} 
Let $ N \geq 2 $ and suppose that the nonlinearity 
$F$ satisfies (A1) and ((A2) or (A3)). 
Then for any $ \d > 0 $ there is $ M( \d)> 0$ such that for all $ c \in [0, \cs]$ and for all solutions 
$ U \in \BE$ of {\rm (TW}$_c${\rm )} such that $ \| \nabla U  \| _{L^2( \R^N)} < M( \d )$ we have 
$$ \n \, |U| - r_0 \n_{L^\ii(\R^N)} \leq \d . $$
\end{lem}

In the two-dimensional case the lifting properties follow immediately from Lemma \ref{liftingfacile}. 
However, in dimension $N=3$,  for  travelling waves $U_{c(\e)} $
 which minimize $E_{c(\e)}$ on $ \Co _{c(\e)} $ the quantity
  $ \Big\| \frac{ \p  U_{c(\e)}}{\p x_1 } \Big\|_{L^2} ^2  $ 
 is large, of order $\simeq \e^{-1}$  as $ \e \to 0$.
 We give a lifting
result for those solutions, based on the fact that 
$ \| \nabla _{x_{\perp}} U_{c(\e)} \|_{L^2} ^2= \frac{N-1}{2} T_{c(\e)} $ is sufficiently small.

\begin{prop} \label{lifting} 
We consider a nonlinearity $F$ satisfying (A1) and ((A2) or (A3)). 
Let $ U \in \BE$ be a travelling wave to {\rm (NLS)} of speed $ c \in [0, \cs ]$.

\medskip

(i) If $N \geq 3$, for any $0 < \d < r_0$  there exists 
$ \mu = \mu ( \d ) > 0 $ such that 
$$
 \Big\| \frac { \p U}{\p x_1} \Big\|_{L^2( \R^N)}  \cdot \|\nabla _{x_{\perp}} U \| _{L^2( \R^N)}^{N-1} \leq \mu( \d) 
\qquad \mbox{ implies } \qquad
\n \, |U| - r_0 \n_{L^\ii(\R^N)} \leq \d . 
 $$
 
(ii) If $ N \geq 4$ and, moreover, (A3) holds or 
$\ds  \Big\| \frac { \p U}{\p x_1} \Big\|_{L^2( \R^N)}  \cdot \|\nabla _{x_{\perp}} U \|_{L^2( \R^N)} ^{N-1} \leq 1$, 
then  for any $ \d > 0 $ there is $ m(\d) > 0$ such that 
$$
\int_{\R^N} |\nabla _{x_{\perp}} U|^2 \, dx \leq m( \d) 
\qquad \mbox{ implies } \qquad
\n \, |U| - r_0 \n_{L^\ii(\R^N)} \leq \d . 
$$
\end{prop}

As an immediate consequence, the three-dimensional travelling 
wave solutions provided by Theorem \ref{thM}  have modulus close to $r_0$
(hence do not vanish) as $ c \to \cs$:

\begin{cor}
\label{sanszero} 
Let $N=3$ and consider a nonlinearity $F$ satisfying (A2) 
and (A4)   with $\G \not= 0$.
Then, the travelling wave solutions $U_{c(\e)} $ to (NLS)
provided by Theorem \ref{thM} 
which satisfy an additional bound $E(U_{c(\e)} ) \leq \frac{C}{\e}$ 
(with $C$ independent on $ \e$)
verify
$$ \n \, |U_{c(\e)}| - r_0 \n_{L^\ii(\R^3)} \to 0 
\quad \quad \mbox{ as } \quad \e \to 0. $$
In particular, for $\e$ sufficiently close to $ 0$  we have $|U_{c(\e)}| \geq r_0 /2 $ 
in $\R^3$.
\end{cor}

\noindent {\it Proof.} 
By the the second identity in (\ref{Pohozaev})   we have 
$$ 
\int_{\R^3} \Big|\frac{\p U_{c(\e)}}{\p x_1} \Big|^2 \, dx = \frac 12  E ( U_{c(\e)}  ) \leq \frac{C}{\e}.
$$
 Moreover, the first identity in \eqref{Pohozaev} and Proposition  \ref{asympto}  $(ii)$  imply
$$ \int_{\R^3} |\nabla _{x_{\perp}} U_{c(\e)} |^2 \, dx  = E_{c(\e)}(U_{c(\e)}) = T_{c(\e)} \leq C \e . $$
Hence
$ \Big\| \frac { \p U_{c(\e)}}{\p x_1} \Big\|_{L^2( \R^3)}   \|\nabla _{x_{\perp}} U_{c(\e)} \| _{L^2( \R^3)}^{2} \leq {C}{\sqrt{\e} } $
and the result follows from Proposition \ref{lifting} $(ii)$. \carre \\

We  give now  some properties of the two-dimensional  travelling wave solutions provided by Theorem \ref{th2d}. 

\begin{prop} \label{prop2d} 
Let  $N = 2 $ and assume that $F$ verifies (A2)   and (A4) 
with $\G \not = 0$. Then there exist   constants $C_1, \, C_2, \, C_3 , \, C_4>0$ and 
$ 0 < k_* < k_\ii$ such that all travelling wave solutions 
$U_k$ provided by  Theorem \ref{th2d} with 
$ 0 < k = \ds \int_{\R^2} |\nabla U_k|^2 \ dx < k_* $ 
satisfy $|U_k| \geq r_0 / 2 $ in $\R^2$,
\be
\label{estim2d}
C_1 k \leq  - Q(U_k) \leq C_2 k, \qquad
C_1 k \leq \int_{\R^2} V(|U_k|^2) \, dx \leq C_2 k, \qquad
C_1 k \leq \int_{\R^2} (\chi ^2(|U_k|) - r_0 ^2) ^2\, dx   \leq C_2 k
\ee
and have a speed $c(U_k) = \sqrt{\cs^2 - \e_k^2}$ satisfying
\be
\label{kifkif}
 C_3 k  \leq \e_k \leq C_4 k . 
\ee
\end{prop}

At this stage, we know that the travelling waves provided by Theorems \ref{th2d} and \ref{thM} 
do not vanish  if  their speed is sufficiently close to $ \cs$. 
Using  the above  lifting results, we  may write such a solution $U_c$  in the form
\be
\label{ansatz}
 U_c (x) = \rho (x) \ex^{i\phi(x)} 
= r_0 \sqrt{1+\e^2 \BA_{\e}(z) }\ \ex^{i\e \vp_{\e} (z)}, 
\quad \quad \mbox{ where } \quad \e = \sqrt{ \cs ^2 - c^2}, \quad z_1 = \e x_1 , \ z_\perp = \e^2 x_\perp , 
\ee
and we use the same scaling as in  \eqref{ansatzKP}. The interest of
writing  the modulus in this way (and not as in  \eqref{ansatzKP}) is just to 
simplify a little   bit the algebra and to have expressions similar to 
those in \cite{BGS1}. Since  $\BA _{\e } = 2 A_{\e} + \e^2 A_{\e}^2$, 
bounds in  Sobolev spaces for $\BA _{\e} $ imply similar Sobolev bounds for $A_{\e }$ and conversely. 
We shall now  find  Sobolev bounds for $\BA _{\e} $ and $\vp _{\e } $. 
It is easy to see that  (TW$_{c}$) is equivalent to the following system for the phase $ \vp$ 
and the modulus $ \rho$ (in the original variable $x$):
\be
\label{phasemod}
\left\{\begin{array}{l}
\ds{ c \frac{\p}{\p {x_1}} (\rho^2 - r_0^2 ) = 
2 \mbox{div} ( \rho^2 \nabla \phi ) },  \\ \\ 
\ds{ \Delta \rho  - \rho |\nabla \phi|^2 +  \rho F(\rho^2) 
= - c  \rho \frac{\p  \phi}{\p {x_1}}  } . 
\end{array}\right.
\ee
Multiplying the second equation by $ 2 \rho$, we write (\ref{phasemod}) in the form
\be
\label{phasemod2}
\left\{\begin{array}{l}
 2 \mbox{div} ( (\rho^2 - r_0 ^2)  \nabla \phi )
- \ds{ c \frac{\p}{\p {x_1}} (\rho^2 - r_0^2 ) }  = - 2 r_0 ^2 \Delta \phi , \\   \\
\ds{ \Delta (\rho^2 - r_0 ^2) - 2 |\nabla U_c|^2 + 2 \rho^2 F(\rho^2) 
+ 2c( \rho ^2 - r_0 ^2) \frac{\p \phi}{\p {x_1}} 
= - 2 c  r_0^2 \frac{\p \phi}{\p {x_1}} } . 
\end{array}\right.
\ee
Let $\eta  =  \rho^2 - r_0 ^2$.  
We apply the operator $ \ds - 2c \frac{ \p }{\p x_1}$ to the first equation in (\ref{phasemod2}) 
and we take the Laplacian of the second one, then we add the resulting equalities to get 
\be
\label{fond}
\left[ \Delta^2 - \cs^2 \Delta + c^2 \frac{ \p ^2}{\p x_1^2} \right] \eta
=  \Delta \left( 2 |\nabla U_c|^2 - 2 c  \eta \frac{\p \phi}{\p x_1}  
- 2 \rho ^2 F( \rho ^2) - \cs^2 \eta \right)
+ 2 c \frac{\p}{ \p x_1} (\mbox{div} (\eta \nabla \phi ) ).
\ee
Since $\cs^2 = - 2 r_0^2 F'(r_0^2)$,   using the Taylor expansion  
$$ 
2 (s + r_0^2) F( s + r_0^2) + \cs^2 s
= - \frac{ \cs ^2}{r_0^2} \left( 1 - \frac{ r_0^4 F'' (r_0^2)}{\cs ^2} \right) s^2 + r_0 ^2 \tilde{F}_3(s), 
$$
where $ \tilde{F}_3(s) = \BO(s^3)$ as $s \to 0$,  we 
see that the right-hand side in \eqref{fond} is at least quadratic in 
$(\eta , \phi) $. Then we perform a scaling and  pass to the variable $ z = ( \e x_1, \e^2 x_{\perp})$ 
(where $ \e = \sqrt{ \cs ^2 - c^2}$), so that \eqref{fond} becomes
\begin{align}
\label{Fonda}
\Big\{ \p_{z_1}^4 - \p_{z_1}^2 - \cs^2 \Delta_{z_\perp} 
+ 2 \e ^2 \p_{z_1}^2 \Delta_{z_\perp} + \e^4 \Delta^2_{z_\perp} \Big\} 
\BA _{\e}  = & \ \BR _{\e},
\end{align}
where $\BR _{\e}$ contains terms at least quadratic in $(\BA _{\e} ,\vp _{\e})$:
\begin{align*}
\BR _{\e} = & \ 
\{ \p_{z_1}^2 + \e ^2 \Delta_{z_\perp} \} \Big[ 
2(1 + \e^2 \BA _{\e}) \Big( (\p_{z_1} \vp _{\e} )^2 
+ \e ^2 |\nabla_{z_\perp} \vp _{\e} |^2 \Big) + \e^2 \frac{(\p_{z_1} \BA _{\e} )^2 
+ \e ^2 |\nabla_{z_\perp} \BA _{\e} |^2}{2(1+ \e ^2 \BA _{\e} )} \Big]
\\ & \ 
- 2 c   \e ^2 \Delta_{z_\perp} ( \BA _{\e}  \p_{z_1} \vp _{\e}) 
+ 2 c   \e ^2 \ds \sum_{j=2}^N \p_{z_1} \p_{z_j} ( \BA _{\e} \p_{z_j} \vp _{\e} )
\\ & \ 
+ \{ \p_{z_1}^2 + \e ^2 \Delta_{z_\perp} \} \Big[ 
\cs^2 \Big( 1 - \frac{r_0^4F''(r_0^2)}{\cs^2} \Big) \BA _{\e}^2 
- \frac{1}{\e ^4} \tilde{F}_3(r_0^2 \e ^2 \BA _{\e}) \Big] .
\end{align*}
In the two-dimensional case, 
uniform bounds (with respect to $\e$) in Sobolev spaces have been derived   in \cite{BGS1}
by using  \eqref{Fonda} and a bootstrap  argument. This technique is  based upon the fact that some kernels related 
to the linear part in \eqref{Fonda}, such as
$$ \Fou^{-1} \Big( \frac{\xi_1^2}{ \xi_1^4 + \xi_1^2 + \cs^2 |\xi_{\perp}|^2 
+ 2 \e ^2 \xi_1^2 |\xi_{\perp}|^2 + \e ^4 |\xi_{\perp}|^4 } \Big) 
\quad \quad {\rm and} \quad \quad 
\Fou^{-1} \Big( \frac{\e ^2 |\xi_\perp|^2}{ \xi_1^4 + \xi_1^2 
+ \cs^2 |\xi_{\perp}|^2 + 2 \e ^2 \xi_1^2 |\xi_{\perp}|^2 
+ \e ^4 |\xi_{\perp}|^4 } \Big) $$
are bounded in $L^p(\R^2)$ for $p$ in some interval $[2, \bar{p})$, 
 uniformly with respect to  $\e$. However, this is no longer true 
in dimension $N=3$: the above mentioned kernels are not in $L^2(\R^3)$ 
(but their Fourier transforms are uniformly bounded), and from the analysis 
in \cite{G}, the kernel
$$ \Fou^{-1} \Big( \frac{\xi_1^2}{ \xi_1^4 + \xi_1^2 + \cs^2 |\xi_{\perp}|^2} \Big) $$
is presumably too 
singular near the origin  to be in $L^p( \R^3) $ if $p \geq 5/3$.
This lack of integrability of the kernels makes the analysis in the three dimensional 
case much more diffcult than in the case $ N=2$. 

\medskip

  One of the main difficulties in the three dimensional case is to 
 prove that for $ \e $ sufficiently small, $ \BA_{\e}$ is uniformly bounded  in $L^p$ for some $ p >2$. 
To do this we use  a suitable decomposition of $\BA _{\e }$ in the Fourier space 
(see the proof of Lemma \ref{Grenouille} below).
 Then  we improve the exponent $p$ 
by using a bootstrap argument,  combining the iterative argument in \cite{BGS1}  
(which uses the quadratic nature of $\BR _{\e}$ in \eqref{Fonda}) and the
appropriate decomposition of $\BA _{\e}$ in the Fourier space.  This  leads to 
some $L^p$ bound with $p> 3 = N$. Once this bound is proved, the proof 
of the $W^{1,p}$ bounds follows the scheme in  \cite{BGS1}. We get:

\begin{prop}
\label{Born} Under the assumptions of Theorem \ref{res1}, there is $\e_0 >0 $ such that 
$ \BA_{\e}  \in W^{4, p}(\R^N)$ and $ \nabla \vp _{e } \in W^{3, p}(\R^N)$ 
for all $ \e \in ( 0, \e_0)$ and all $ p \in ( 1, \infty)$. 
Moreover,   for any $ p \in (1, \ii ) $ 
there exists $C_{p} >0$ satisfying for all $ \e \in (0, \e _0) $
\beq
\label{goodestimate}
\n \BA_{\e } \n_{L^p} + \n \nabla \BA_{\e } \n_{L^p}  + \n \p^2_{z_1} \BA_{\e } \n_{L^p} 
+ \e  \n \p_{z_1} \nabla_{z_\perp} \BA_{\e } \n_{L^p} 
+ \e ^2 \n \nabla_{z_\perp}^2 \BA_{\e } \n_{L^p} \leq C_p  \qquad \mbox{ and } 
\eeq
\beq
\label{goodestimate2}
\begin{array}{l}
\| \p _{z_1} \vp _{\e } \|_{L^p} + \e \| \nabla_{z_{\perp }} \vp _{\e } \| _{L^p} 
+ \| \p _{z_1}^2  \vp _{\e } \|_{L^p} + \e \| \nabla_{z_{\perp }} \p _{z_1} \vp _{\e } \| _{L^p}
 + \e ^2  \| \nabla_{z_{\perp }} ^2 \vp _{\e } \| _{L^p} 
\\
\\
+ \| \p _{z_1}^3  \vp _{\e } \|_{L^p} + \e \| \nabla_{z_{\perp }} \p _{z_1} ^2 \vp _{\e }\| _{L^p} 
+ \e ^2  \| \nabla_{z_{\perp }} ^2 \p _{z_1} \vp _{\e }\| _{L^p}    \leq C_p.
\end{array}
\eeq
The estimate \eqref{goodestimate} is also valid with $A_{\e }$ instead of $ \BA_{\e }$. 
\end{prop}

Once these bounds are established, the estimates in Proposition \ref{asympto}  show  that 
$(\cs^{-1} \p_{z_1} \vp_n )_{n\geq 0}$ is a minimizing sequence for the problem \eqref{minimiz} if $N=2$, 
respectively for  the  problem (\ref{miniGS}) if $N=3$.  Since 
Theorems \ref{gs2d} and \ref{gs} provide compactness properties  for 
minimizing sequences, we get  (pre)compactness of 
$ ( \cs^{-1} \p_{z_1} \vp_n )_{n\geq 0}$ in 
$\mathscr{Y}(\R^N) \hookrightarrow L^2(\R^N)$, and then we complete 
the proof of Theorem \ref{res1} by standard interpolation 
in Sobolev spaces.

\subsection{On the higher dimensional case}

It is  natural to ask   what happens in the transsonic limit  in dimension $N\geq 4$.
Firstly, it should be  noticed  that even for the 
Gross-Pitaevskii nonlinearity  the problem is critical if $N=4$ and 
supercritical in higher dimensions, hence Theorem \ref{thM}  does not apply directly. 

The first crucial step is to investigate the behaviour 
of $T_c$ as $c \to \cs$. In particular, in order to be able to 
use Proposition \ref{lifting}  to show that the solutions 
are vortexless in this limit, we would need to prove that $T_c \to 0$ 
as $c \to \cs$. We have not been able to prove (or disprove)  
this in dimension $N=4$ and $N=5$, except for the case $ \G = 0$.
 Quite surprisingly, for nonlinearities satisfying (A3) and (A4) (this is the case for both the Gross-Pitaevskii and the cubic-quintic nonlinearity), 
 this is not true  in dimension higher than $5$, as shown by the following

\begin{prop}
\label{dim6} 
Suppose that $F$ satisfies  (A3) and (A4) (and $\G$ is arbitrary). 
If $ N \geq 6$,  there exists $\d >0 $ 
such that for any $ 0 \leq c \leq \cs$ and for any  nonconstant solution   $U \in \BE$ 
of {\rm (TW$_{c}$)}, we have
$$ 
E(U) + c Q(U) \geq \d . 
$$
In particular,
$$
 \inf_{0 < c < \cs} T_c > 0 . 
$$
The same conclusion holds if $ N \in \{ 4, 5 \}$ provided that $ \G = 0$. 
\end{prop}

Therefore  we do not know if the solutions constructed in 
Theorem \ref{thM} (for a subcritical nonlinearity) may vanish or not as 
$c\to \cs$ if $N\geq 6$. On the other hand  we can show, in 
any space dimension $N\geq 4$, that we cannot scale the 
solutions in order to have compactness and convergence to a localized 
and nontrivial object in the transonic limit as soon as the quantity 
$E+cQ$ tends to zero.

\begin{prop}
\label{vanishing} 
Let $N\geq 4$ and
suppose that $F$ satisfies (A2), (A3)  and (A4) (and $\G$ is arbitrary). 
 Assume that there exists  a 
sequence $(U_n, c_n)$ such that $ c_n \in (0, \cs]$, $U_n \in \BE$ is 
a nonconstant solution of {\rm (TW$_{c_n}$)} and $E_{c_n}(U_n) \to 0$ 
as $n\to \ii$. Then, for $n$ large enough, there exist 
$\a_n,\beta_n, \la_n, \s_n \in \R$, $ A_n \in H^1( \R^N)$ and $ \vp_n \in \dot{H}^1( \R^N)$
 uniquely determined such that
$$ U_n(x) = r_0 \Big( 1 + \a_n A_n (z) \Big) 
\exp \Big( i \beta_n \vp_n (z) \Big) , \quad \quad \quad \mbox{ where } \quad 
z_1 = \la_n x_1, \quad z_\perp = \s_n x_\perp , $$
$$
\a_n \to 0 \qquad \mbox{  and } \qquad
 \n A_n \n_{L^\ii(\R^N)} = \n A_n \n_{L^2(\R^N)} 
= \n \p_{z_1} \vp_n \n_{ L^2(\R^N)} 
= \n \nabla_{z_\perp} \vp_n \n_{L^2(\R^N)} = 1 . 
$$
Then we have   $c_n \to \cs$ and
$$ \n \p_{z_1} A_n \n_{L^2(\R^N)} \to 0 \qquad \mbox{ as } n \to +\ii. $$
\end{prop}

Consequently, even if one could show 
that $T_c \to 0$ as $ c \to \cs $ 
 in space dimension $4$ or $5$, 
 we would not have  a nontrivial limit (after rescaling) of the  corresponding  rarefaction pulses.

\section{Three-dimensional ground states for (KP-I) 
\label{proofGS}}

We recall  the 
anisotropic Sobolev inequality (see \cite{BIN}, p. 323): for $N \geq 2$ and 
for any $2 \leq p < \frac{2(2N-1)}{2N-3}$, there exists $C=C(p,N)$ 
such that for all $\Theta \in \BC_c^\ii(\R^N)$ we have
\be
\label{Pastis}
\n \p_{z_1} \Theta \n_{L^p(\R^N)} \leq C 
\n \p_{z_1} \Theta \n_{L^2(\R^N)}^{1 - \frac{(2N-1)(p-2)}{2p}} 
\n \p^2_{z_1} \Theta \n_{L^2(\R^N)}^{\frac{N(p-2)}{2p}} 
\n \nabla_{z_\perp} \Theta \n_{L^2(\R^N)}^{\frac{(N-1)(p-2)}{2p}} .
\ee
This shows that the energy $\mathscr{E}$ is well-defined on $\mathscr{Y}(\R^N)$ 
if $N=2$ or $N=3$. By (\ref{Pastis}) and the density of  $\p_{z_1}  C_c^{\infty}(\R^3)$  in $\mathscr{Y}(\R^3)$ 
we get for any $ w \in \mathscr{Y}(\R^3)$:
\be
\label{SobAnis}
\n w \n_{L^3(\R^3)} \leq C 
\n w \n_{L^2(\R^3)}^{\frac16} 
\n \p_{z_1} w \n_{L^2(\R^3)}^{\frac12} 
\n \nabla_{z_\perp} \p_{z_1}^{-1} w \n_{L^2(\R^3)}^{\frac13} .
\ee

On the other hand, the following identities hold for any solution 
$\BW \in \mathscr{Y} (\R^N)$ of (SW):
\be
\label{identites}
\left\{\begin{array}{ll}
\ds{ \int_{\R^N} \frac{1}{\cs^2}\, (\p_{z_1} \BW)^2 
+ |\nabla_{z_\perp} \p_{z_1}^{-1} \BW|^2 
+ \frac{\G}{2}\, \BW^3 + \frac{1}{\cs^2}\, \BW^2 \ dz = 0} \\ \ \\ 
\ds{ \int_{\R^N} \frac{-1}{\cs^2}\, (\p_{z_1} \BW)^2 
+3  |\nabla_{z_\perp} \p_{z_1}^{-1} \BW|^2 
+ \frac{\G}{3}\, \BW^3 + \frac{1}{\cs^2}\, \BW^2 \ dz } = 0 \\ \ \\ 
\ds{ \int_{\R^N} \frac{1}{\cs^2}\, (\p_{z_1} \BW)^2 
+ |\nabla_{z_\perp} \p_{z_1}^{-1} \BW|^2 
+ \frac{\G}{3}\, \BW^3 + \frac{1}{\cs^2}\, \BW^2 \ dz = 
\frac{2}{N-1} \int_{\R^N} |\nabla_{z_\perp} \p_{z_1}^{-1} \BW|^2 \ dz }.
\end{array}\right.
\ee
The first identity is obtained by multiplying (SW) by $\p_{z_1}^{-1} \BW$ 
and integrating, whereas the two other equalities  are the Pohozaev 
identities associated to the scalings in the $z_1$ and 
$z_\perp$ variables respectively. Formally, they are obtained   by multiplying (SW) 
by $z_1 \BW$ and $z_\perp \cdot \nabla_{z_\perp} \p_{z_1}^{-1} \BW$ 
respectively and integrating by parts (see \cite{dBSIHP} for   a
complete justification). Combining the equalities in \eqref{identites} we get 
\be
\label{Ident}
\left\{\begin{array}{c}
\ds{ \int_{\R^N} \frac{1}{\cs^2} \, (\p_{z_1} \BW )^2 \ dz = 
\frac{N}{N-1} \int_{\R^N} | \nabla_{z_\perp} \p_{z_1}^{-1} \BW|^2 \ dz } \\ \\ 
\ds{ \frac{\G}{6} \int_{\R^N} \BW^3 \ dz = 
- \frac{2}{N-1} \int_{\R^N} | \nabla_{z_\perp} \p_{z_1}^{-1} \BW|^2 \ dz } \\ \\
\ds{ \int_{\R^N} \frac{1}{\cs^2} \, \BW^2 \ dz 
= \frac{7-2N}{N-1} \int_{\R^N} | \nabla_{z_\perp} \p_{z_1}^{-1} \BW|^2 \ dz }.
\end{array}\right.
\ee
Notice that for  $N\geq 4$  we have $ 7-2N <0$ and the last equality  implies $\BW=0$.

\medskip

We recall the following results about the ground states of (SW) 
and the compactness of minimizing sequences in $ \bY(\R^3)$.

\begin{Lemma} [\cite{dBSIHP}, \cite{dBSSIAM}]
\label{gs3}
Let $ N=3$ and $ \G \neq 0$. 

\medskip

(i) For $ \la \in \R^*$, denote 
$I_{\la } = \inf \Big\{ \| w \|_{\bY (\R^3)} ^2 \; | \; \ds \int_{\R^3} w^3 (z)\, dz = \la \Big\}.$
Then for any $ \la \in \R^*$ we have $ I_{\la } > 0$ and there is $ w_{\la } \in \bY(\R^3)$ such that 
$\ds \int_{\R^3} w_{\la} ^3 (z)\, dz = \la$ and $\|  w _{\la} \|_{\bY (\R^3)} ^2 = I_{\la }$. 
Moreover, any sequence $(w_n)_{n \geq 1} \subset \bY(\R^3)$ such that 
$\ds \int_{\R^3} w_n^3 (z)\, dz \to \la $ and $\|  w _{n} \|_{\bY (\R^3)} ^2 \to I_{\la }$
has a subsequence that converges in $ \bY(\R^3)$ (up to translations) to a minimizer of $ I_{\la }$.

\medskip

(ii) There is $ \la^* \in \R^*$ such that $ w^* \in \bY(\R^3)$ is a ground state for {\rm (SW)} 
(that is, minimizes the action $ \bS$ among all solutions of {\rm (SW)}) 
if and only if $w^*$ is a minimizer of  $I_{\la ^*}$. 
\end{Lemma}

The first part of Lemma \ref{gs3} is a consequence of 
the proof of Theorem 3.2 p. 217 in \cite{dBSIHP} 
and the second part follows from Lemma 2.1 p. 1067 in \cite{dBSSIAM}.

\medskip

{\it Proof of Theorem \ref{gs}.} 
Given $ w \in \bY(\R^3)$ and $ \si > 0$, we denote 
$P(w) = \ds \int_{\R^3} \frac{1}{\cs ^2} w^2 +  \frac{1}{\cs ^2} |\p_{z_1} w|^2 + \frac{\G}{3} w^3 \, dz$ 
and $ w_{\si} (z) = w( z_1, \frac{ z_{\perp}}{\si} )$. It is obvious that 
$$
\begin{array}{c}
\ds \int_{\R^3} w_{\si }^p \, dz = \si ^2 \int_{\R^3} w^p \, dz , \qquad
\int_{\R^3} |\p _{z_1} (w_{\si })|^2 \, dz = \si ^2 \int_{\R^3} |\p_{z_1}w|^2 \, dz  \qquad \mbox{ and } 
\\ \\
\ds \int_{\R^3} |\nabla_{z_{\perp}} \p_{z_1}^{-1} (w_{\si })|^2 \, dz 
=  \int_{\R^3} |\nabla_{z_{\perp}} \p_{z_1}^{-1} (w)|^2 w|^2 \, dz .
\end{array}
$$

Let $w^*$ be a ground state for (SW) (the existence of $ w^*$ is guaranteed by Lemma \ref{gs3} above). 
Since $ w^*$ satisfies (\ref{identites}), we have $ P(w^*) = 0 $ and 
$ \bS (w^*) = \ds \int_{\R^3} |\nabla_{z_{\perp}} \p_{z_1}^{-1} (w^*)|^2 w|^2 \, dz.$
Consider $ w \in \bY(\R^3)$ such that $ w \neq 0$ and $ P(w) = 0$. 
Then 
$ \ds \frac{ \G}{3} \int_{\R^3} w^3 \, dz = - \frac{1}{\cs ^2} \int_{\R^3} w^2 + | \p _{z_1} w|^2 \, dz < 0 $
and it is easy to see that there is $ \si > 0 $ such that 
$\ds \int_{\R^3} w_{\si }^3 \, dz = \int_{\R^3} (w ^*)^3 \, dz = \la ^*$. 
From Lemma \ref{gs3} it follows that $\| w_{\si} \|_{\bY(\R^3)} ^2\geq \| w^* \|_{\bY(\R^3)} ^2, $ that is 
$$
\frac{ \si ^2} { \cs ^2} \int_{\R^3} w^2 + |\p _{z_1} w|^2 \, dz  + 
\int_{\R^3} |\nabla_{z_{\perp}} \p_{z_1}^{-1} w |^2  \, dz 
\geq 
\frac{ 1} { \cs ^2} \int_{\R^3} (w^*)^2 + |\p _{z_1} w^*|^2 \, dz  + 
\int_{\R^3} |\nabla_{z_{\perp}} \p_{z_1}^{-1} w^*|^2  \, dz .
$$
Since $ P(w) = 0 $ and $ P(w^*) = 0$ we have 
$$
\frac{ \si ^2} { \cs ^2} \int_{\R^3} w^2 + |\p _{z_1} w|^2 \, dz 
= - \si ^2 \frac { \G}{3}   \int_{\R^3} w^3 \, dz 
= -  \frac { \G}{3}   \int_{\R^3} (w^*)^3 \, dz 
= \frac{1}{ \cs ^2} \int_{\R^3} (w^*)^2 + |\p _{z_1} w^*|^2 \, dz 
$$
and the previous inequality gives 
$
\ds \int_{\R^3} |\nabla_{z_{\perp}} \p_{z_1}^{-1} w|^2  \, dz 
\geq 
\int_{\R^3} |\nabla_{z_{\perp}} \p_{z_1}^{-1} w^*|^2  \, dz , 
$
that is $ \bS(w) \geq \bS( w^*)$. 
So far we have proved that the set $ \Po = \{ w \in \bY(\R^3) \; |  \; w \neq 0, \; P(w) = 0 \}$ 
is not empty and any ground state $w^*$ of (SW) minimizes the action $ \bS$ in this set. 
It is then clear that for any $ \si > 0$, $ w_{\si }^* $ also belongs to $ \Po$ and mnimizes $ \bS$ on $ \Po$.

\medskip

Conversely, let $ w \in \Po $ be such that $ \bS(w) = \bS_*$. 
Let $w^*$ be a ground state for (SW). 
It is clear that  $ \ds \int_{\R^3} |\nabla_{z_{\perp}} \p_{z_1}^{-1} w|^2  \, dz = \bS_* 
= \int_{\R^3} |\nabla_{z_{\perp}} \p_{z_1}^{-1} w^*|^2  \, dz $. 
As above, there is a unique $ \si  >0 $ such that 
$\ds \int_{\R^3} w_{\si}^3 \, dz =  \int_{\R^3} (w^*)^3 \, dz = \la ^*$ and then we have 
$\ds \int_{\R^3} w_{\si } ^2 + |\p _{z_1} w_{\si }|^2 \, dz 
= \int_{\R^3} (w^*)^2 + |\p _{z_1} w^*|^2 \, dz $. We find 
$ \| w_{\si }\|_{\bY(\R^3)}^2 = \| w^*\|_{\bY(\R^3)}^2 = I_{\la ^*} $, 
thus $ w_{\si }$ is a minimizer for $ I_{\la ^*} $ and  Lemma \ref{gs3} (ii) implies that 
$ w_{\si }$ is a ground state for (SW). 

\medskip

Let $(\BW_n)_{n \geq 1} $ be a sequence satisfying $(i)$, $(ii)$ and  $(iii)$. 
We have $ P(\BW_n ) \to 0 $ and 
$$
\frac{\G}{3} \int_{\R^3} \BW_n^3 \, dz 
= P(\BW_n) - \frac{1}{\cs ^2} \int_{\R^3} \BW_n ^2 + | \p _{z_1} \BW_n |^2 \, dz 
\in \left[ \frac{ - 2 m_2}{\cs ^2}, - \frac{ m_1}{2 \cs ^2} \right]
\qquad \mbox{ for all $n$ sufficiently large.}
$$
We infer that there are $ n_0 \in \N$, $ \underline{ \si}, \bar{\si} > 0 $ 
and a sequence $(\si _n )_{ n \geq n_0}\subset [ \underline{ \si}, \bar{\si} ] $ such that 
$ \ds \int_{\R^3} ( (\BW_n)_{\si _n}) ^3 \, dz = \la ^*$ for all $ n \geq n_0$. 
Moreover, 
$$
\begin{array}{rcl}
\| (\BW_n)_{\si _n} \| _{\bY(\R^3)} ^2 
& = & \ds  \frac{ \si _n ^2}{\cs ^2}  \int_{\R^3}  \BW_n ^2 + | \p _{z_1} \BW_n |^2 \, dz   + 
    \int_{\R^3} |\nabla_{z_{\perp}} \p_{z_1}^{-1} \BW_n|^2  \, dz 
\\
\\
& = & \ds \si _n ^2 \left( P(\BW_n) - \frac{ \G}{3} \int_{\R^3}  \BW_n ^3 \right) + ( \bS( \BW _n) - P(\BW_n) )
\\
& = & \ds ( \si _n ^2 - 1) P(\BW_n) + \bS( \BW _n) -  \frac{ \G}{3} \int_{\R^3}  (\BW_n)_{\si _n} ^3 \, dz.
\end{array}
$$
Passing to the limit in the above equality we get 
$$
\ds \liminf_{n \ra \infty } \| (\BW_n)_{\si _n} \| _{\bY(\R^3)} ^2 
= \liminf_{n \ra \infty } \bS( \BW _n) - \frac{ \G}{3} \la ^* \leq \bS_* - \frac{ \G}{3} \la ^* 
= \bS( w^*) -  \frac{ \G}{3} \int_{\R^3} (w^*)^3\, dz 
= \| w^* \|_{\bY(\R^3)} ^2 = I_{\la ^*}. 
$$
Hence there is a subsequence of $((\BW_n)_{\si _n})_{n \geq 1}$ which is a minimizing sequence for $I_{\la ^*}$. 
Using Lemma \ref{gs3} we infer that there exist a subsequence $(n_j)_{ j \geq 1}$ such that 
$ \si_{n _j} \to \si \in [\underline{\si}, \bar{\si}]$, 
 a sequence $ (z_j)_{j \geq 1} \subset \R^3$ and a minimizer $ \BW$ of $ I_{\la ^*}$ 
 (hence a ground state for (SW)) such that 
 $(\BW_{n_j})_{\si_{n_j}} ( \cdot - z_j) \to \BW$ in $ \bY(\R^3)$. 
 It is then straightforward that $ \BW_{n_j}( \cdot - z_j) \to \BW_{\frac{1}{\si}}$ in $ \bY(\R^3)$. 
\carre \\

We may give an alternate proof of Theorem \ref{gs} which does not rely directly 
on the analysis in \cite{dBSIHP}, \cite{dBSSIAM} by following the strategy 
of \cite{Maris}, which can be adapted to our problem up to some details.

\section{Proof of Theorem \ref{res1}}

\subsection{Proof of Proposition \ref{asympto}}
\label{preuveasympto}

For some given real valued functions $A_\e$ and $\vp_\e$, 
we consider the mapping
$$ U_{\e}(x) = |U_{\e}| (x) \ex^{i\phi(x)} = 
r_0 \Big( 1+ \e^2 A_\e(z) \Big) \ex^{i\e \vp_\e(z)} , 
\quad \quad \mbox{ where } \quad z=(z_1,z_\perp) = ( \e  x_1 , \e^2 x_\perp) . $$
It is obvious that $ U_{\e } \in \BE$ provided that
 $ A_{\e} \in H^1 ( \R^N)$ and $ \nabla  \vp _{\e} \in L^2( \R^N)$. 
If $\e$ is small and $A_\e$ is uniformly bounded in $\R^N$, $U_{\e}$ does not 
vanish and the momentum $Q(U_{\e})$ is   given by
$$ Q(U_{\e}) = - \int_{\R^N} ( |U_{\e}|^2 - r_0^2 ) \frac{\p \phi}{\p x_1} \ dx 
= - \e^{5-2N} r_0^2 \int_{\R^N} \Big( 2 A_\e + \e^2 A_\e^2 \Big) 
\frac{\p \vp_\e}{\p z_1} \ dz , $$
while the energy of $U_{\e }$ is 
\begin{align*}
E(U_{\e} ) = & \ \int_{\R^N} |\nabla U_{\e}|^2 + V(|U_{\e}|^2) \ dx \\
= & \ \e^{5-2N} r_0^2 \int_{\R^N} 
(\p_{z_1} \vp_\e)^2 \Big( 1 + \e^2 A_\e \Big)^2 
+ \e^2 |\nabla_{z_\perp} \vp_\e|^2 \Big( 1 + \e^2 A_\e \Big)^2 
+ \e^2 (\p_{z_1} A_\e)^2 
+ \e^4 |\nabla_{z_\perp} A_\e|^2 \\
& \hspace{2cm} + \cs^2 A_\e^2 
+ \e^2 \cs^2 \Big( 1 - \frac{4r_0^4}{3\cs^2} F''(r_0^2) \Big) A_\e^3 
+ \frac{\cs^2}{\e^4} V_4 \Big( \e^2 A_\e \Big) \ dz ,
\end{align*}
where we have used  the Taylor expansion 
\be
\label{V}
 V \Big( r_0^2 (1 + \a)^2 \Big) = r_0^2 \Big\{ 
\cs^2 \a^2 + \cs^2 \Big( 1 - \frac{4r_0^4}{3\cs^2} 
F''(r_0^2) \Big) \a^3 + \cs^2 V_4(\a) \Big\}
= r_0 ^2  \cs ^2 \Big\{ \al ^2 + \Big(\frac{\G}{3} - 1 \Big) \al ^3 + V_4 ( \al ) \Big\}
\ee
 with 
$ V_4(\a) = \BO(\a^4) $ as   $\a \to 0 . $
Consequently, with   $\cs^2 = c^2(\e) + \e^2$ we get
\begin{align}
\label{develo}
 E_{c(\e)} (U_{\e }) & =  \ E(U_{\e} ) + c(\e) Q (U_{\e}) \nonumber \\
& =  \e^{5-2N} r_0^2 \int_{\R^N} 
(\p_{z_1} \vp_\e)^2 \Big( 1 + \e^2 A_\e \Big)^2 
+ \e^2 |\nabla_{z_\perp} \vp_\e|^2 \Big( 1 + \e^2 A_\e \Big)^2 
+ \e^2 (\p_{z_1} A_\e)^2 
+ \e^4 |\nabla_{z_\perp} A_\e|^2 \nonumber \\
& \hspace{1cm} + \cs^2 A_\e^2 
+ \e^2 \cs^2 \Big( 1 - \frac{4r_0^4}{3\cs^2} F''(r_0^2) \Big) A_\e^3 
+ \frac{\cs^2}{\e^4} V_4 \Big( \e^2 A_\e \Big)
- c(\e) \Big( 2 A_\e + \e^2 A_\e^2 \Big) \p_{z_1} \vp_\e \ dz \nonumber \\
& =   \e^{7-2N} r_0^2 \int_{\R^N} \frac{1}{\e^2} 
\Big( \p_{z_1} \vp_\e - c(\e) A_\e \Big)^2 
+ (\p_{z_1} \vp_\e)^2 ( 2 A_\e + \e^2 A_\e^2 ) 
+ |\nabla_{z_\perp} \vp_\e|^2 ( 1 + \e^2 A_\e )^2 + (\p_{z_1} A_\e)^2 
\nonumber \\
& \hspace{1cm} + \e^2 |\nabla_{z_\perp} A_\e|^2 + A_\e^2 
+ \cs^2 \Big( 1 - \frac{4r_0^4}{3\cs^2} F''(r_0^2) \Big) 
A_\e^3 + \frac{\cs^2}{\e^6} V_4( \e^2 A_\e) 
- c(\e) A_\e^2 \p_{z_1} \vp_\e \ dz .
\end{align}
Since the first term  in the last integral is penalised by $\e^{-2}$, 
in order to get sharp estimates on $E_{c(\e)} $ one needs
$\p_{z_1} \vp_\e \simeq c(\e) A_\e$.

\medskip

Let  $N=3$. By Theorem \ref{gs}, 
 there exists a  ground state $ A \in \bY(\R^3)$ for (SW). 
It follows from Theorem 4.1 p. 227  in \cite{dBSSIAM} that  $ A  \in H^s(\R^3)$  for any $s \in \N$.  
Let $ \vp = \cs \p_{z_1 }^{-1} A$. 
We use  \eqref{develo} with 
$A_\e(z) = \frac{\la  \cs}{  c(\e) } A( \la z_1 , z_\perp)$ 
and $\vp_\e(z) = \vp( \la z_1 , z_\perp)$. 
 For $\e>0$ small 
and $\la \simeq 1$ (to be chosen later) we define
$$
 U_\e (x) = |U_\e| (x) \ex^{i\phi_\e(x)} =
r_0 \Big( 1+ \e^2 \frac{\cs}{c(\e)} \la A(z) \Big) \ex^{i\e \vp(z)} , 
\quad \quad \quad \mbox{ where } \qquad z=(z_1,z_\perp) = ( \e \la x_1 , \e^2 x_\perp) .
$$
Notice that $U_{\e}$  does not vanish if  $\e$ is sufficiently small. 
 Since $\p_{z_1} \vp = \cs A$, we have 
$ \p_{z_1} \vp_\e (z)= \la \p_{z_1} \vp ( \la z_1, z_{\perp} )= \la \cs A (\la z_1, z_{\perp} )=  c(\e) A_\e (z) $
and therefore 
\begin{align*}
\la E_{c(\e)} (U_\e) = & \ \cs^2 r_0^2 \e \int_{\R^3} 
\la^3 \frac{\cs}{c(\e)} A^2 \Big( 2 A + \e^2 \frac{\cs}{c(\e)} \la A^2 \Big) 
+ \la^2 |\nabla_{z_\perp} \p_{z_1}^{-1} A|^2 
\Big( 1 + \e^2 \frac{\cs}{c(\e)} \la A \Big)^2 
+ \frac{\la^4}{c^2(\e)} (\p_{z_1} A)^2 \nonumber \\
& \hspace{2cm} + \e^2 \frac{\la^2}{c^2(\e)} |\nabla_{z_\perp} A|^2 
+ \frac{\la^2}{c^2(\e)} A^2 
+ \frac{\cs^3}{c^3(\e)} \la^3 \Big( 1-\frac{4r_0^4}{3\cs^2} F''(r_0^2) \Big) A^3 
+ \frac{1}{\e^6} V_4\Big( \e^2 \frac{\cs}{c(\e)} \la A \Big) 
\nonumber \\ 
& \hspace{2cm} - \la^3 \frac{\cs}{c(\e)} A^3 \ dz \nonumber \\ 
= & \ \cs^2 r_0^2 \e \int_{\R^3} 
\la^3 \frac{\cs}{c(\e)} \Big( 1 + \frac{\cs^2}{c^2(\e)} 
\Big[ \frac{\G}{3} - 1 \Big] \Big) A^3 
+ \la^2 |\nabla_{z_\perp} \p_{z_1}^{-1} A|^2 
\Big( 1 + \e^2 \frac{\cs}{c(\e)} \la A \Big)^2 
+ \frac{\la^4}{c^2(\e)} (\p_{z_1} A)^2 
\nonumber \\
& \hspace{2cm} + \frac{\la^2}{c^2(\e)} A^2 
+ \e^2 \frac{\la^2}{c^2(\e)} |\nabla_{z_\perp} A|^2 
+ \e^2 \la^4 \frac{\cs^2}{c^2(\e)} A^4 
+ \frac{1}{ \e^6} V_4\Big( \e^2 \frac{\cs}{c(\e)} \la A \Big) \ dz .
\end{align*}
On the other hand,
\begin{align*}
 \la \int_{\R^3} |\nabla_{\perp} U_\e|^2 \ dx = & \, 
r_0^2 \e \int_{\R^3} |\nabla_{z_\perp} \vp|^2 
\Big( 1+ \e^2 \la \frac{\cs}{c(\e)} A \Big)^2 
+ \e^2 \la^2 \frac{\cs^2}{c^2(\e)} |\nabla_{z_\perp} A|^2 \ dz \\ 
= & \, \cs^2 r_0^2 \e \int_{\R^3} 
|\nabla_{z_\perp} \p_{z_1}^{-1} A|^2 
\Big( 1+ \e^2 \la \frac{\cs}{c(\e)} A \Big)^2 
+ \e^2 \frac{\la^2}{c^2(\e)} |\nabla_{z_\perp} A|^2\ dz .
\end{align*}
Hence  $U_\e $ satisfies the constraint $P_{c(\e )} (U_{\e}) = 0 $ (or equivalently 
$ \ds E_{c(\e)} (U_\e ) = \int_{\R^3} |\nabla_{\perp} U_\e |^2 \ dx $) 
if and only if $G(\la, \e ^2) = 0$, where 
\begin{align*}
 G (\la ,  \e^2 ) = & \, \int_{\R^3} 
\la^3 \frac{\cs}{c(\e)} \Big( 1 + \frac{\cs^2}{c^2(\e)} 
\Big[ \frac{\G}{3} - 1 \Big] \Big) A^3 
+ \la^2 |\nabla_{z_\perp} \p_{z_1}^{-1} A|^2 
\Big( 1 + \e^2 \frac{\cs}{c(\e)} \la A \Big)^2 
+ \frac{\la^4}{c^2(\e)} (\p_{z_1} A)^2 \nonumber \\
& \hspace{2cm} + \frac{\la^2}{c^2(\e)} A^2 
+ \e^2 \frac{\la^2}{c^2(\e)} |\nabla_{z_\perp} A|^2 
+ \e^2 \la^4 \frac{\cs^2}{c^2(\e)} A^4 
+ \frac{1}{\e^6} V_4\Big( \e^2 \frac{\cs}{c(\e)} \la A \Big) \ dz \\ 
& \ - \int_{\R^3} |\nabla_{z_\perp} \p_{z_1}^{-1} A|^2 
\Big( 1+ \e^2 \la \frac{\cs}{c(\e)} A \Big)^2 
+ \e^2 \frac{\la^2}{c^2(\e)} |\nabla_{z_\perp} A|^2\ dz . 
\end{align*}
Denote $\epsilon = \e^2$.
Since $A$ is a ground state for (SW),   it satisfies 
 the Pohozaev identities \eqref{identites}. The last of these identities  is  $ \mathscr{S}(A) = 
\ds \int_{\R^3} |\nabla_{z_\perp} \p_{z_1}^{-1} A |^2 \ dz $, or equivalently 
$$ 
G( \la = 1, \epsilon = 0 ) = 0 . 
$$
A straightforward computation  using \eqref{Ident} gives
$$ \frac{\p G}{\p \la}_{|(\la = 1, \epsilon = 0)} = 
\int_{\R^3} \G A^3 + 2 |\nabla_{z_\perp} \p_{z_1}^{-1} A|^2 
+ \frac{4}{\cs^2} (\p_{z_1} A)^2 + \frac{2}{\cs^2} A^2 \ dz 
= 3 \int_{\R^3} |\nabla_{z_\perp} \p_{z_1}^{-1} A|^2 \not = 0 . $$
Then the implicit function theorem implies that there 
exists a function $\epsilon \longmapsto  \la(\epsilon) = 1 + \BO ( \epsilon) = 1 + \BO(\e^2)$ such that 
for all $\epsilon$ sufficiently small we have  $G(\la(\epsilon),\epsilon ) = 0$, 
that is $U_{c(\e)} $ satisfies the Pohozaev  identity   $P_{c(\e)}(U_{\e}) = 0$. 
Choosing $ \la = \la ( \e ^2)$ and taking into account the last indetity in \eqref{identites}, we find
$$ T_{c(\e)} \leq E_{c(\e)}(U_\e) = 
\int_{\R^3} |\nabla_{\perp} U_\e |^2 \ dx = 
\cs^2 r_0^2 \e \int_{\R^3} |\nabla_{z_\perp} \p_{z_1}^{-1} A|^2 + \BO(\e^3) 
= \cs^2 r_0^2 \e \bS_{\rm min} + \BO(\e^3)
$$
and the proof of $(ii)$ is complete. 

\medskip

Next we turn our attention to the case $N=2$. 
 Let $A = \cs^{-1} \p_{z_1} \vp \in \bY(\R^2)  $ 
 be   a ground state of (SW). The existence of $A$ is given   by Theorem  \ref{gs2d}.
By Theorem 4.1 p. 227 in \cite{dBSIHP}  we have $A \in H^s(\R^2)$ for all  $s \in \N$. 
For $\e$ small,   we define the map
$$
 U_\e(x) = |U_\e| (x) \ex^{i\phi_\e(x)} = 
r_0 \Big( 1+ \e^2 \frac{\cs}{c(\e)} A(z) \Big) \ex^{i\e \vp(z)} , 
\quad \quad \quad \mbox{ where } \qquad z=(z_1,z_2) = ( \e x_1 , \e^2 x_2) . 
$$
From the above computations and \eqref{Ident} we have
\begin{align*}
k_\e = & \ \int_{\R^2} |\nabla U_\e|^2 \ dx = 
r_0^2 \e \int_{\R^2} (\p_{z_1} \vp_\e)^2 \Big( 1 + \e^2 A_\e \Big)^2 
+ \e^2 (\p_{z_1} A_\e)^2 
+ \e^2 (\p_{z_2} \vp_\e)^2 \Big( 1 + \e^2 A_\e \Big)^2 
+ \e^4 (\p_{z_2} A_\e)^2 \ dz \\
= & \ r_0^2 \cs^2 \e \int_{\R^2} 
A^2 \Big( 1 + \frac{\e^2 \cs}{c(\e)} A \Big)^2 
+ \frac{\e^2}{c^2(\e)} (\p_{z_1} A)^2 
+ \e^2 (\p_{z_2} \p_{z_1}^{-1} A)^2 \Big( 1 + \frac{\e^2 \cs}{c(\e)} A \Big)^2 
+ \frac{\e^4}{c^2(\e)} (\p_{z_2} A)^2 \ dz \\ 
= & \ r_0^2 \cs^2 \Big\{ \e \int_{\R^2} A^2 \ dz 
+ \e^3 \int_{\R^2} \Big( 2 A^3 + \frac{(\p_{z_1} A)^2}{\cs^2} 
+ (\p_{z_2} \p_{z_1}^{-1} A)^2 \Big) \ dz 
+ \BO(\e^5) \Big\} \\ 
= & \ r_0^2 \cs^2 \Big\{ \e  \frac32 \cs ^2 \bS(A) 
+ \e^3 \Big( 2 - \frac{12}{\G} - \frac12 \Big) \bS(A) 
+ \BO(\e^5) \Big\}
\end{align*}
It is easy to see  that $\e \mapsto k_\e$ is a 
smooth increasing diffeomorphism from an interval $[0,\bar{\e}]$ onto an interval 
$[0,\bar{k}= \bar{k}_{\bar{\e}}]$, and that 
$ \e = \ds{\frac{k_\e}{r_0^2 \cs^2 \| A \|_{L^2}^2}} + \BO(k_\e^3)
= \frac{k_{\e }}{\frac 32 r_0 ^2 \cs ^4 \bS (A) } + \BO(k_{\e}^3) $ 
as $\e \to 0$.  Moreover, denoting  $U_{\e}^\s (x) =  U_{\e}(x/ \s)$  we have 
$$ \int_{\R^2} |\nabla U_\e^\s|^2 \ dx = \int_{\R^2} |\nabla U_\e|^2 \ dx $$
because $N=2$.  
Using the test function  $U_\e^\s$, it follows that 
$$  I_{\rm min} (k_\e) \leq I(U_\e^\s )    \qquad \mbox{ for any } \s >0 . $$
Since $Q(U_{\e}) < 0$, the mapping 
$$
\si \longmapsto  I(U_\e^\s ) = Q( U _\e^\s) + \int_{\R^2} V(|U_\e^\s|^2) \ dx 
= \s Q(U_\e) + \s^2 \int_{\R^2} V(|U_\e|^2) \ dx 
$$
achieves its  minimum at  $ \s_0 = 
\ds{\frac{- Q(U_\e)}{\ds{2 \int_{\R^2} V(|U_\e|^2)}}} > 0 $, and the minimum value is
$
 I(U_\e^{\s_0}) = \ds \frac{- \ds Q^2(U_\e)}{ \ds{4\int_{\R^2} V(|U_\e|^2) \ dx}}. 
$
Hence
$$ 
 I_{\rm min} (k_\e) \leq I(U_\e^{\s_0} ) 
=  \frac{- Q^2(U_\e)}{ \ds{4\int_{\R^2} V(|U_\e|^2) \ dx}} . 
$$
Using (\ref{V}) and (\ref{Ident}) we find 
$$
\begin{array}{rcl}
\ds \int_{\R^2} V(|U_\e|^2) \ dx & = & \ds \cs^2 r_0^2 \e \int_{\R^2} A^2 
+ \e^2 \Big( \frac{\G}{3} - 1 \Big) A^3 
+ \frac{1}{\e^4} V_4 ( \e^2 A ) \ dz
\\
\\
& = & \ds\frac 32 \cs ^4 r_0 ^2 \bS(A) \e - \cs ^2 r_0 ^2 \Big( \frac{\G}{3} - 1 \Big) \frac{ 6}{\G} \bS(A) \e ^3 + \BO(\e^5)
\end{array}
$$
and
\begin{align*}
Q(U_\e) = - \e r_0^2 \cs \int_{\R^2} \Big( 2 A^2 + \e^2 A^3 \Big) \ dz 
= -3 r_0 ^2 \cs ^3 \bS(A) \e + r_0 ^2 \cs \frac{ 6}{\G} \bS(A)  \e^3.
\end{align*}
Finally we obtain 
$$
\begin{array}{l}
\ds I_{\rm min} (k_\e)  +\ds  \frac{ k_{\e}}{\cs ^2}   
 \leq  \frac{- \ds Q^2(U_\e)}{ \ds{4 \int_{\R^2} V(|U_\e|^2) \ dx}} + \frac{1}{\cs^2} \ds \int_{\R^2} |\nabla U_{\e} |^2 \, dx 
\\
\\
\ds =  - \frac{\left(-3 \cs ^2 + \frac{6}{\G} \e^2 \right) ^2 r_0 ^4 \cs ^2 \bS^2(A) \e^2}
{4 \Big[ \frac 32 \cs ^2 - \Big( 2 - \frac{6}{\G} \Big) \e^2 + \BO(\e^4) \Big]r_0 ^2 \cs ^2 \bS(A) \e }
+ 
\left[ \frac 32 r_0 ^2 \cs ^2 \e + r_0 ^2 \Big( \frac 32 - \frac{12}{\G} \Big) \e^3 + \BO( \e ^5 ) \right] \bS(A)
\\ 
\\
= \ds - \frac{ \left(3  r_0 ^2 \cs ^2 \e ^3 + \BO( \e ^5) \right) \bS(A) }
{2 \left[3 \cs ^2 - \left( 4 - \frac{12}{\G} \right) \e^2 + \BO( \e^4) \right] }
= - \frac 12 r_0 ^2 \bS(A) \e ^3 + \BO(\e^5) 
\\
\\
\ds  =  - \frac 12 r_0 ^2 \bS(A) \left[ \frac{ k_{\e}}{\frac 32 r_0 ^2 \cs ^4 \bS(A) } + \BO(k_{\e}^3) \right]^3 + 
\BO \left( \left(\frac{ k_{\e}}{\frac 32 r_0 ^2 \cs ^4 \bS(A) } + \BO(k_{\e}^3) \right)^5 \right)
= \frac {-4k_{\e} ^3  }{27 r_0 ^4 {\cs } ^{12} \bS_{\rm min}^2 } + \BO(k_{\e}^5 ).
\end{array}
$$
Since $ \e \longmapsto k_{\e}$ is a diffeomorphism from  $[0,\bar{\e}]$ onto 
$[0,\bar{k}]$,  Proposition \ref{asympto} (i) is proven.
\carre

\subsection{Proof of Proposition \ref{monoto}}

Given a function $ f $ defined on $ \R^N$ and $ a, \,  b > 0$, 
we denote $ f_{a, b}(x) = f( \frac{ x_1}{a}, \frac{ x_{\perp}}{b}).$

\medskip

By Proposition 2.2 p. 1078 in \cite{M2}, any solution 
of (TW$_{c}$) belongs to $W_{loc}^{2, p }(\R^N)$ for all $ p \in [2, \infty)$, hence to $C^{1, \al }(\R^N)$ for all $ \al \in (0, 1)$. 
 
\medskip

($i$) Let $ U $ be a minimizer of $ E_c= E + cQ$ on $ \Co _c$ (where $\Co _c$ is as in \eqref{Cc})
such that $ \psi $ solves (TW$_c$). 
Then $U $  satisfies the Pohozaev identities \eqref{Pohozaev}. 

\medskip

If $ Q( U ) > 0$, let $ \tilde{ U }(x) = U ( - x_1, x_{\perp})$, 
so that $ Q( \tilde{ U }) = - Q( U ) < 0$ and 
$ P_c( \tilde{ U } ) = P_c( U ) - 2 c Q( U ) = - 2cQ( U ) < 0$. 
Since for any function $ \phi \in \Eo $  we have 
\be
\label{Pca}
P_c( \phi _{a, 1}) = \frac 1a \int_{\R^N} \Big| \frac{ \p \phi }{\p x_1} \Big| ^2 \, dx 
+ a \frac{N-3}{N-1} \int_{\R^N}  |\nabla_{x_{\perp}} \phi |^2 \, dx + c Q( \phi ) + a \int_{\R^N} V( |\phi |^2) \, dx, 
\ee
we see that there is $ a_0 \in (0, 1) $ such that $ P_c( \tilde{U}_{a_0 ,1 }) = 0$ . 
We infer that 
$$
T_c \leq E_c( \tilde{U}_{a_0, 1} ) 
= \frac{2}{N-1} \int_{\R^N}  |\nabla_{x_{\perp}} \tilde{U}_{a_0, 1} |^2 \, dx
= a_0 \frac{2}{N-1} \int_{\R^N}  |\nabla_{x_{\perp}} U |^2 \, dx
= a_0 E_c( U) 
= a_0 T_c, 
$$
contradicting the fact that $ T_c > 0$. 
Thus  $ Q( U ) \leq 0$. 

\medskip

Assume that $ Q( U ) = 0. $ 
From the identities (\ref{Pohozaev}) with $ Q( U ) = 0$ we get 
\be
\label{pr}
\ds \int_{\R^N} \Big| \frac{ \p U}{\p x_1} \Big| ^2 \, dx = - \frac{1}{N-2}  \int_{\R^N} V(|U |^2) \, dx  
 \qquad \mbox{  and  } \qquad
\ds \int_{\R^N} | \nabla _{x_{\perp}} U  | ^2 \, dx = - \frac{N-1}{N-2}  \int_{\R^N} V(|U |^2) \, dx .
\ee
Since $ U \in \Eo $ and $ U $ is not constant, 
 necessarily 
$ \ds \int_{\R^N} V(|U |^2) \, dx = - (N-2)  \int_{\R^N} \Big| \frac{ \p U}{\p x_1} \Big| ^2 \, dx <0$
and this implies that the potential $V$ must achieve negative values. 
Then it follows from Theorem 2.1 p. 100 in \cite{brezis-lieb} that there is $ \tilde{\psi }_0 \in \Eo $ 
such that 
$\ds \int_{\R^N} |\nabla \tilde{\psi}_0 |^2 \, dx = 
\inf \Big\{ \int_{\R^N} |\nabla \phi |^2 \, dx \; \Big| \; \phi \in \BE, \; 
\int_{\R^N} V(|\phi |^2) \, dx = -1 \Big\}.$
Using Theorem 2.2 p. 102 in \cite{brezis-lieb} we see that there is $ \si > 0$ such that, denoting 
$ \psi _0 = (\tilde{\psi }_0)_{\si, \si}$ and 
$ - v_0 = \ds \int_{\R^N} V(|\psi _0|^2) \, dx = - \si ^N$, 
we have $ \Delta \psi _0 + F( |\psi _0 |^2 ) \psi _0 = 0 $ in $ \R^N$. 
Hence  $\psi _0 $ solves (TW$_0$) and 
$$
\ds \int_{\R^N} |\nabla {\psi}_0 |^2 \, dx = \inf \Big\{ \int_{\R^N} |\nabla \phi |^2 \, dx \; \Big| \; 
\phi \in \BE, \; \int_{\R^N} V(|\phi |^2) \, dx = - v_0 \Big\}.
$$
Since all minimizers of this problem solve (TW$_{0}$) (after possibly rescaling), 
we know that they are $C^1$ in $\R^N$ and then Theorem 2 p. 314 in \cite{MarisARMA} 
imply that they are all radially symmetric (after translation). In particular, we have $Q (\psi _0 ) = 0 $ 
and  $  \ds \int_{\R^N} \Big| \frac{ \p \psi _0}{\p x_j} \Big| ^2 \, dx 
= \frac 1N  \ds \int_{\R^N} |\nabla \psi _0 | ^2 \, dx $ for $ j = 1, \dots, N$.
By Lemma 2.4 p. 104 in \cite{brezis-lieb} we know that $ \psi _0 $ satisfies the Pohozaev identity  
$  
  \ds \int_{\R^N} |\nabla \psi _0 | ^2 \, dx  = - \frac{ N}{N-2} v_0.
$
It follows that $P_c( \psi _0 ) = 0 $, hence $ \psi _0 \in \Co_c $ and we infer that 
$ E_c ( \psi _0 ) \geq T_c$, 
that is 
$ \ds \frac{2}{N-1}  \int_{\R^N} |\nabla_{x_{\perp}} \psi _0 |^2 \, dx 
\geq \frac{2}{N-1} \ds \int_{\R^N} |\nabla_{x_{\perp}} U  |^2 \, dx$.
Taking into account \eqref{pr} and the radial symmetry of $\psi _0$, 
this gives   $ v_0 \geq - \ds \int_{\R^N}  V( |U  |^2 ) \, dx$. 

On the other hand, by scaling it is easy to see that $ \psi _0 $ is a minimizer of the functional 
$ \phi \longmapsto \| \nabla \phi \|_{L^2( \R^N)} ^2 $ in the set 
$ \Po = \Big\{ \phi \in \BE \; \Big| \;  \ds \int_{\R^N} |\nabla \phi | ^2 \, dx = - \frac{N}{N-2} \int_{\R^N} V(|\phi |^2) \, dx \Big\}$. 
By \eqref{pr} we have $ U \in \Po$, hence 
$\| \nabla U \|_{L^2( \R^N)} ^2 \geq \| \nabla \psi _0\|_{L^2( \R^N)} ^2$ and consequently 
$  - \ds \int_{\R^N}  V( |U  |^2 ) \, dx \geq v_0$. 
Thus 
$\| \nabla U \|_{L^2( \R^N)} ^2 = \| \nabla \psi _0\|_{L^2( \R^N)} ^2$, 
$   \ds \int_{\R^N}  V( |U  |^2 ) \, dx = \int_{\R^N}  V( |\psi _0 |^2 )$ and $ U $ minimizes  
$ \| \nabla \cdot \|_{L^2(\R^N)}^2$ 
in the set $\Big\{ \phi \in \BE \; \Big| \; \ds \int_{\R^N}  V( |\phi  |^2 ) \, dx = - v_0 \Big\}$. 
By Theorem 2.2 p. 103 in \cite{brezis-lieb}, 
$ U $ solves the equation $ \Delta U + \la  F(|U |^2) U = 0 $ in $ \Do '(\R^N)$ for some $ \la>0$ 
and using the Pohozaev identity associated to this equation we see that $ \la = 1$, hence $ U $ solves (TW$_0$). 
Since $U $ also solves (TW$_c$) for some $ c >0$ and $ \frac{\p U}{\p x_1}$ is continuous, 
 we must have $   \frac{\p U}{\p x_1} = 0 $ in $ \R^N$. 
Together with the fact that $ U \in \BE$, this implies that $ U $ is constant, a contradiction. 
Therefore we cannot have $ Q( U ) = 0$ and we conclude that  $Q(U ) < 0$. 

\medskip

($ii$)  Fix  $ c_0 \in (0, \cs)$ and let $ U_0 \in \BE$ be a minimizer of $ E_{c_0}$ on $ \Co _{c_0}$, 
as given by Theorem \ref{thM}. 
It follows from \eqref{Pca} that $P_c ((U_0)_{a, 1} ) = \frac 1a R_{c, U_0}(a), $ where 
\be
\label{Rca}
R_{c, U_0}(a) = \int_{\R^N} \Big| \frac{ \p U_0}{\p x_1} \Big|^2 \, dx + a c Q( U_0) + 
a^2 \left[ \frac{N-3}{N-1} \int_{\R^N} |\nabla _{x_{\perp} } U_0 |^2 \, dx + \int_{\R^N} V(|u_0|^2) \, dx \right]
\ee
is a polynomial in $a$ of degree at most 2. 
It is clear that $R_{c, U_0} (0 ) > 0$, $R_{c_0, U_0} (1) = P_{c_0} (U_0) = 0 $ and for any $ c > c_0 $ we have
$R_{c, U_0}(1) = P_{c_0}(U_0) + (c - c_0) Q(U_0) < 0 $ because $ Q(U_0) <0$. 
Hence there is a unique $ a(c) \in (0,1)$ such that $R_{c, U_0} (a(c)) = 0$, which means $P_{c}((U_0)_{a(c), 1}) =0$. 
We infer that
\be 
\label{ac}
T_c \leq E_c( (U_0)_{a(c), 1}) = \frac{2}{N-1} \int_{\R^N} |\nabla _{x_{\perp} } (U_0)_{a(c), 1} |^2 \, dx 
= a(c)  \frac{2}{N-1} \int_{\R^N} |\nabla _{x_{\perp} } U_0 |^2 \, dx   
= a(c) T_{c_0} .
\ee
Since $ a(c) \in (0,1)$, we have  proved that $T_c < T_{c_0}$ whenever $ c_0 \in (0, \cs ) $ and $ c \in (c_0, \cs)$, 
thus $ c \longmapsto T_c$ is decreasing. By a well-known result of Lebesgue, the function 
$ c \longmapsto T_c $  has a derivative a.e.

\medskip

($iii$) 
Notice that \eqref{ac} holds whenever $c_0$, $U_{c_0}$ are as above and $a(c)$ is a positive 
root of $R_{c, U_0}$. Using the Pohozaev identities \eqref{Pohozaev} we find
$$
 2\int_{\R^N} \Big| \frac{ \p U_0}{\p x_1} \Big|^2 \, dx 
 = \frac{2}{N-1} \int_{\R^N} |\nabla _{x_{\perp} } U_0 |^2 \, dx - c_0 Q (U_0) 
 = T_{c_0} - c_0 Q( U_0) \qquad \mbox{ and then}
$$
\be
\label{ps}
\frac{N-3}{N-1} \int_{\R^N} |\nabla _{x_{\perp} } U_0 |^2 \, dx + \int_{\R^N} V(|u_0|^2) \, dx 
= - c_0 Q( U_0) - \int_{\R^N} \Big| \frac{ \p U_0}{\p x_1} \Big|^2 \, dx 
= - \frac 12 c_0 Q(U_0) - \frac 12 T_{c_0} .
\ee
We now distinguish two cases: $R_{c, U_0}$ has degree one or two.

\medskip

Case $(a)$:
If $ \ds \frac{N-3}{N-1} \int_{\R^N} |\nabla _{x_{\perp} } U_0 |^2 \, dx + \int_{\R^N} V(|u_0|^2) \, dx = 0$, 
then $R_{c, U_0}$ has degree one  and we have 
$\ds \int_{\R^N} \Big| \frac{ \p U_0}{\p x_1} \Big|^2 \, dx +  c_0 Q( U_0) =0$ because 
$ P_{c_0}(U_0) = 0$.  Since $R_{c, U_0}$ is an affine function, we find $a(c) = \frac{ c_0}{c}$
for all $ c > 0$, hence $ a( c_0) = 1$.  Moreover,  the left-hand side in \eqref{ps} is zero, thus we have
$ c_0 Q(U_0) + T_{c_0} = 0$  and  consequently
$ a'( c_0) = - \frac{ 1}{c_0} = \frac{Q(U_0)}{T_{c_0}}$.

\medskip

Case $(b)$: If $ \ds \frac{N-3}{N-1} \int_{\R^N} |\nabla _{x_{\perp} } U_0 |^2 \, dx + \int_{\R^N} V(|u_0|^2) \, dx 
\not = 0$, $ R_{c, U_0} $ has degree two, and the discriminant of this second-order polynomial 
is equal to
$$
\Delta _{c, U_0} = ( c^2 - c_0 ^2) Q^2 ( U_0) + T_{c_0}^2. 
$$
Consequently $R_{c, U_0}$ has real roots as long as $ ( c^2 - c_0 ^2) Q^2 ( U_0) + T_{c_0}^2 \geq 0$. 
It is easy to see that if there are real roots, at least one of them is positive. 
Indeed, $R_{c, U_0}(0) > 0 > R_{c, U_0}'(0) $. If $ \Delta _{c, U_0} \geq 0 $, 
no matter of  the sign of the leading order coefficient 
$ \frac{N-3}{N-1} \int_{\R^N} |\nabla _{x_{\perp} } U_0 |^2 \, dx + \int_{\R^N} V(|u_0|^2) \, dx \not = 0 $,
the smallest positive root $a(c)$ of $ R_{c, U_0} $ is given by the formula
\be
\label{root}
a(c) = \frac{- c Q( U_0) - \sqrt{ ( c^2 - c_0 ^2) Q^2 ( U_0) + T_{c_0}^2}}{- c_0 Q(U_0) - T_{c_0}}
= \frac{ - c_0 Q( U_0) + T_{c_0}}{- cQ(U_0) + \sqrt{ ( c^2 - c_0 ^2) Q^2 ( U_0) + T_{c_0}^2}} .
\ee
Therefore, the function $ c \longmapsto a(c)$ is defined on the interval 
$ [\tilde{c}_0, \infty )$ where $ \tilde{c}_0 = \sqrt{ c_0 ^2 - \frac{ T_{c_0^2}}{Q^2(U_0)}}< c_0$, 
it is differentiable  on $ (\tilde{c}_0, \infty )$ and  $ a( c_0) = 1$. 
Moreover, a straightforward computation gives $ a'( c_0) =  \frac{Q(U_0)}{T_{c_0}}$. 
Note that in Case $(a)$, the last expression in \eqref{root} is equal to 
$ \frac{c_0}{c} $, which is then indeed $ a(c)$.

By \eqref{ac} we have $ T_c \leq a(c) T_{c_0}$ and passing to the limit we get 
$ \ds \lim_{c \ra c_0,\,  c < c_0} T_c \leq \lim_{c \ra c_0, \, c < c_0} a(c) T_{c_0} = T_{c_0}$. 
Since $ c \longmapsto T_c$ is decreasing, $T_c > T_{c_0} $ for $ c < c_0 $ and  we see that it is left contiuous at $ c_0$. 
Moreover, we have 
$$
\frac{ T_c - T_{c_0}}{c - c_0 } \leq \frac{ a(c) - a(c_0)}{c - c_0} T_{c_0 } \quad \mbox{  for } c > c_0, 
\quad \mbox{ respectively } \quad
\frac{ T_c - T_{c_0}}{c - c_0 } \geq \frac{ a(c) - a(c_0)}{c - c_0} T_{c_0 } \quad \mbox{  for } c \in [\tilde{c}_0, c_0). 
$$
Passing to the limit in the above inequalities we obtain, since $ a'( c_0) =  \frac{Q(U_0)}{T_{c_0}}$ 
in Cases $(a)$ and $(b)$, 
$$
\limsup_{ c \ra c_0, \, c > c_0} \frac{ T_c - T_{c_0}}{c - c_0 } \leq  a'( c_0) T_{c_0} = Q(U_0), 
\qquad \mbox{ respectively } \qquad
\liminf_{ c \ra c_0, \, c < c_0} \frac{ T_c - T_{c_0}}{c - c_0 } \geq  a'( c_0) T_{c_0} = Q(U_0).
$$
It is then clear that if $ c \longmapsto T_c$ is differentiable at $ c_0$,  necessarily
$ \ds \frac{d T_c}{dc}_{|c=c_0} = Q(U_0) .$
 
\medskip

($iv$) Fix $ c_* \in ( c_0, \cs)$. Passing to a subsequence we may assume that 
$ c_0 < c_n < c_*$ for all $n$ and $ Q( U_n ) \to - q _0 \leq 0$. 
Then $ T_{c_0 } > T_{c_n} > T_{c_*} > 0 $ 
and $ (c_0 ^2 - c_n ^2) Q^2( U_n) + T_{c_n}^2 >  ( c_0 ^2 - c_n ^2) Q^2( U_n) + T_{c_*}^2 > 0$ 
for all sufficiently large $n$. 
Hence for large $n$ we may use \eqref{ac} and \eqref{root} with $( c_n, c_0)$ instead of $ ( c_0, c)$ and we get 
$$
T_{c_0} \leq 
\frac{ - c_n Q( U_n) + T_{c_n}}{- c_0 Q(U_n) + \sqrt{ ( c_0^2 - c_n ^2) Q^2 ( U_n) + T_{c_n}^2}}
T_{c_n}.
$$
Since $T_{c_n}$ has a positive limit, passing to the limit as $ n \to \infty $ in the 
above inequality and using the monotonicity of $ c \longmapsto T_c$ we get 
$ \ds T_{c_0} \leq \liminf_{ n \ra \infty} T_{c_n} = \liminf_{c  \ra c_0, \, c > c_0} T_c$. 
This and the fact that $ T_c$ is decreasing and left continuous imply that $T_c$ is continuous at $ c_0$.

\medskip

($v$) Let $ 0 < c_1 < c_2 < \cs $ and $U_1, \; U_2,$  $ q_1 = Q(U_1)< 0 $,  $ q_2 = Q(U_2) < 0$ 
be as in Proposition \ref{monoto} ($v$).
If $ c_1 ^2 \leq c_2 ^2 - \frac {T_{c_2}^2}{q_2^2}$, the inequality  in Proposition \ref{monoto} ($v$) 
obviously holds. From now on we assume that $ c_1 ^2 >  c_2 ^2 - \frac {T_{c_2}^2}{q_2^2} $. 
The two discriminants $ \Delta_{c_2 , U_1} = ( c_2^2 - c_1^2 ) q_1^2 + T_{c_1}^2 $ 
and $ \Delta_{c_1 , U_2} = ( c_1^2 - c_2^2 ) q_2^2 + T_{c_2}^2 $ are positive: since 
$ 0 < c_1 < c_2 $ for the first one, and by the assumption $ c_1 ^2 >  c_2 ^2 - \frac {T_{c_2}^2}{q_2^2}$ 
for the second one. Therefore, we may use \eqref{ac} and \eqref{root} with the couples $(c_1, c_2)$, 
respectively $(c_2, c_1)$ instead of $(c_0, c)$ to get
$$
T_{c_2} \leq \frac{ - c_1 q_1 + T_{c_1}}{- c_2 q_1 + \sqrt{ (c_2 ^2 - c_1 ^2 ) q_1 ^2 + T_{c_1}^2}} T_{c_1}, 
\qquad \mbox{ respectively } \qquad
T_{c_1} \leq \frac{ - c_2 q_2 + T_{c_2}}{ - c_1 q_2 + \sqrt{ (c_1 ^2 - c_2 ^2 ) q_2 ^2 + T_{c_2}^2}} T_{c_2}.
$$
Since $ T_{c_i} > 0 $, we must have
$$
\frac{ - c_1 q_1 + T_{c_1}}{ - c_2 q_1 + \sqrt{ (c_2 ^2 - c_1 ^2 ) q_1 ^2 + T_{c_1}^2}}
\cdot
\frac{ - c_2 q_2 + T_{c_2}}{ - c_1 q_2 + \sqrt{ (c_1 ^2 - c_2 ^2 ) q_2 ^2 + T_{c_2}^2}} 
\geq 1.
$$
We set $ y_1 = - \frac{ T_{c_1}}{c_1 q_1} > 0 $, and recast this inequality as
\be
\label{ineqmagique}
\frac{ 1 + y_1}{\frac{c_2}{c_1} + \sqrt{ \frac{ c_2^2}{c_1^2} - 1 + y_1^2}} 
\geq 
\frac{ - c_1 q_2 + \sqrt{ (c_1 ^2 - c_2 ^2 ) q_2 ^2 + T_{c_2}^2}}{ - c_2 q_2 + T_{c_2}} 
= \frac{ 1 + \sqrt{ 1 - \frac{c_2^2}{c_1^2} + \frac{T_{c_2}^2}{c_1^2 q_2^2} }}{ 
\frac{c_2}{c_1} - \frac{T_{c_2}}{c_1 q_2} } .
\ee
Denoting, for $y \in \R$, 
$ g (y) = \ds \frac{ 1 + y}{ \frac{c_2}{c_1} + \sqrt{ \frac{c_2^2}{c_1^2} - 1 + y^2}} $, 
\eqref{ineqmagique} is exactly
$$
g \Big( - \frac{ T_{c_1}}{c_1 q_1} \Big) = g (y_1) \geq 
g \Big( \sqrt{ 1 - \frac{c_2^2}{c_1^2} + \frac{T_{c_2}^2}{c_1^2 q_2^2} } \Big) .
$$
If we show that $g$ is increasing, then we obtain
$$
- \frac{ T_{c_1}}{c_1 q_1} \geq \sqrt{ 1 - \frac{c_2^2}{c_1^2} + \frac{T_{c_2}^2}{c_1^2 q_2^2} } , 
\quad \quad \quad {\rm or} \quad \quad \quad 
\frac{ T_{c_1}^2}{q_1^2} - c_1 ^2 \geq \frac{ T_{c_2}^2}{q_2^2} - c_2 ^2 ,
$$
which is the desired inequality. To check that $g$ is increasing, we simply compute
$$
 g' (y) = \ds \frac{ \ds \frac{c_2^2}{c_1^2} - 1 + \ds \frac{c_2}{c_1} \sqrt{ \ds \frac{c_2^2}{c_1^2} - 1 + y^2} - y}{ \Big( \ds \frac{c_2}{c_1} + \sqrt{ \frac{c_2^2}{c_1^2} - 1 + y^2} \Big)^2 
 \sqrt{ \frac{c_2^2}{c_1^2} - 1 + y^2}} ,
$$
which is positive since $ \frac{c_2}{c_1} > 1 $ and 
$ \sqrt{ \frac{c_2^2}{c_1^2} - 1 + y^2} > |y| $.

\medskip

($vi$) 
Since $ c \longmapsto - T_c$ is increasing, by a well-known result of Lebesgue this map  is differentiable a.e., 
the function $  c \longmapsto  \frac{ d T_c}{d c}$ belongs to $L_{loc}^1( 0, \cs)$ and for any $ 0 < c_1 < c_2 < \cs $ we have 
$ \ds \int_ {c_1  }^{c_2} - \frac{ d T_c}{d c} \, dc \leq -T_{c_2} + T_{c_1}. $

We recall that $ c( \e ) = \sqrt{ \cs ^2 - \e^2}$ for all $ \e \in ( 0, \cs )$. 
If $N=3$, (A2) and ( A4) hold and $ \G \neq 0$, by Proposition \ref{asympto} ($ii$) 
there is $ K > 0 $ such that $T_{c(\e)} \leq K \e $ for all sufficiently small $ \e$. 
Thus for $ n \in \N$ large we have 
$$
 \int_{c(2/n)}^{c(1/n)} - \frac{d T_c}{dc} \ dc 
\leq T_{c(2/n)} - T_{c(1/n)} \leq T_{c(2/n)} 
\leq \frac{2K}{n} . 
$$
Hence there exists  $ c_n \in (c(2/n) , c(1/n))$ 
such that $c \mapsto T_c $ is differentiable  at $ c_n $ and
$$ 
- \frac{d T_c}{dc}_{|c=c_n} 
\leq \frac{1}{c(\frac 1n)- c(\frac 2n) } \cdot \frac{2K}{n}
\leq K' n . $$
Let $\e_n = \sqrt{\cs^2 - c_n^2}  $, so that 
$c(\e_n) = c_n$. Since $ c(2/n) \leq c_n \leq c(1/n) $, we have
$
\frac1n \leq \e_n \leq \frac2n , 
$
so that $ \e _n \to 0 $ as $ n \to \infty$. 
Let $U_n$ be a minimizer 
of $ E_{c_n} $ on $ \Co _{c_n}$, scaled so
that $U_n$ solves (TW$_{c_n}$). From ($i$)  and ($iii$) we get 
$$ |Q(U_n) | = - Q( U_n) = - \frac{d T_c}{dc}_{|c=c_n }  
\leq K'n \leq \frac{ 2K'}{\e _n}  . $$
Since
$ E(U_n) + c_n Q(U_n) = T_{c_n } = \BO(\e_n) $, it follows that
$$
 E(U_n)\leq - c_n Q( U_n) + T_{c_n} \leq \frac{K''}{ \e_n} 
$$
and the proof is complete. \carre

\subsection{Proof of Proposition \ref{global3}}

We postpone the proof of Proposition \ref{convergence} and we prove Proposition \ref{global3}.

\medskip

Let $ (\e_n)_{n \geq 1}$ be the sequence given by Proposition \ref{monoto} ($vi$). 
For each $n$ let $ U_n \in \BE$ be a minimizer of $ E_{c_n}$ on $ \Co _{c_n}$ which solves (TW$_{c_n}$).
Passing to a subsequence if necessary and using Proposition \ref{convergence}, we may assume that 
$ (\e_n)_{n \geq 1}$ is strictly decreasing, that  $(\e _n , U_n)_{n \geq 1}$ 
satisfies the conclusion of Theorem \ref{res1} and
\be
\label{en}
\frac 12 r_0 ^2 \cs ^4 \bS_{\rm min} \frac{ 1}{\e _n } < E( U_n) < 2 r_0 ^2 \cs ^4 \bS_{\rm min} \frac{ 1}{\e _n }, 
\ee
\be
\label{mom}
\frac 12 r_0 ^2 \cs ^3 \bS_{\rm min} \frac{ 1}{\e _n } < - Q( U_n) < 2 r_0 ^2 \cs ^3 \bS_{\rm min} \frac{ 1}{\e _n } 
\qquad \mbox{ for all } n.
\ee

We shall argue by contradiction. More precisely, we shall prove by contradiction 
that there exists $ \e _* > 0 $ such that for any $ \e \in ( 0, \e _* )$ and for any 
minimizer $U$ of $E_{c(\e)} $ on $ \Co_{c (\e)}$ scaled so that $U$ satisfies 
(TW$_{c(\e )}$), we have
$$ |Q(U) | \leq \frac{5r_0^2 \cs^3 \bS_{\rm min}}{\e} . $$
In view of Proposition \ref{asympto} ($ii$), we then infer that
$$
 E(U) = T_{c(\e)} - c(\e) Q(U) \leq \frac{K}{\e} 
$$
for some constant $ K $ depending only on $ r_0$, $ \cs $ and $ \bS_{\rm min} $, which 
is the desired result. We thus assume that there exist infinitely many $n$'s such that 
there is $ \tilde{\e} _n \in ( \e _n, \e _{n-1})$ and there is a minimizer 
$ \tilde{U}_n $ of $ E_{c( \tilde{\e}_n)}$ on $ \Co _{c( \tilde{\e}_n)}$ which 
satisfies (TW$_{c( \tilde{\e}_n)}$) and 
\be
\label{mauvais}
|Q(\tilde{U}_n )| = - Q(\tilde{U}_n ) > 5 r_0 ^2 \cs ^3 \bS_{\rm min} \frac{ 1}{\tilde{\e} _n }.
\ee
Passing again to a subsequence of $(\e_n)_{n \geq 1}$, we may assume that \eqref{mauvais} holds for 
all $ n \geq 1$. Then for each $ n \in \N^*$ we define
$$
\begin{array}{rcl}
I_n & = & \Big\{ \e \in ( \e _n, \e_{n-1}) \; \Big| \; \mbox{ for all } \e ' \in [ \e_n, \e ] 
\mbox{ and for any minimizer } U_{\e'} \mbox{ of } 
E_{c(\e')} \mbox{ on } \Co _{c(\e')} 
\\
& & 
\mbox{ which solves (TW$_{c(\e')}$) there holds  }
 | Q( U_{\e'})| \leq 4 r_0 ^2 \cs ^3 \bS_{\rm min} \cdot \frac{ 1}{ \e' } \Big\}
\end{array}
$$
and 
$$
\e_n^{\#} = \sup I_n.
$$
By Proposition \ref{monoto} ($v$), for   $ \e'\in ( \e _n , \cs)$ and  for any minimizer $ U_{\e'}$ 
of $E_{c(\e')}$  on $ \Co _{c(\e')} $ which solves (TW$_{c(\e')}$) we have
$$
\frac{ T_{c( \e')}^2}{Q^2(U_{\e'})} + (\e')^2 \geq \frac{ T_{c( \e _n)}^2}{Q^2(U_{n})} + \e_n ^2, 
$$
which can be written as 
$\ds \frac{Q^2(U_{\e'})}{ T_{c( \e')}^2} \leq \frac{Q^2(U_{n})}{T_{c( \e _n)}^2 + ( \e _n ^2 - (\e ')^2) Q^2( U_n) } \; $
  and this gives
\be
\label{estimate1}
(\e')^2  Q^2(U_{\e'}) \leq \frac{(\e')^2 Q^2(U_{n}) T_{c( \e')}^2}{T_{c( \e _n)}^2 + ( \e _n ^2 - (\e ')^2) Q^2( U_n) }.
\ee
The mapping $ \e \longmapsto T_{c(\e)}$ is right continuous (because $ c \longmapsto T_c$ is left continuous) and using \eqref{mom} we find
$$
\lim_{\e'\ra \e_n, \, \e'> \e _n} 
\frac{(\e')^2 Q^2(U_{n}) T_{c( \e')}^2}{T_{c( \e _n)}^2 + ( \e _n ^2 - (\e ')^2) Q^2( U_n) } 
= \e_n ^2 Q^2 ( U_n) < ( 2 r_0 ^2 \cs ^3 \bS_{\rm min} )^2.
$$
Thus all $ \e '\in ( \e_n, \e_{n-1}) $ sufficiently close to $ \e_n$ belong to $I_n$. 
In particular, $ I_n$ is not empty. 
On the other hand, \eqref{mauvais} implies that any $ \e'\in ( \tilde{\e }_n , \e_{n-1}) $ does not belong to $ I_n$, hence 
$\e_n ^{\#} = \sup I_n \in ( \e _n, \tilde{\e}_n ]\subset ( \e_n, \e_{n-1}).$

Let $ U_n^{\#}$ be a minimizer of $E_{c( \e_n^{\#})}$ on $\Co _{c( \e_n^{\#})} $ which solves (TW$_{c( \e_n^{\#})}$).
We claim that 
\be
\label{claim1}
| Q(U_n^{\#}) | = 4 r_0 ^2 \cs ^3  \bS_{\rm min} \frac{1}{\e_n^{\#} } .
\ee
Indeed, proceeding as in \eqref{estimate1} we have for any $ \e'\in ( \e_n, \e_n^{\#})$ and any minimizer $  U_{\e'} $ of 
$E_{c(\e')}$  on $ \Co _{c(\e')} $ which satisfies (TW$_{c(\e')}$) 
\be 
\label{estimate2}
(\e_n^{\#})^2 Q^2(U_n^{\#}) \leq 
\frac{\left( \frac{\e_n^{\#}}{\e' } \right)^2 (\e')^2 Q^2(U_{\e '}) T_{c( \e_n^{\#})}^2}
{T_{c( \e ')}^2 + \left( 1 - \left( \frac{\e_n^{\#}}{\e'} \right)^2\right)(\e ')^2 Q^2( U_{\e'}) }.
\ee
Notice that $ (\e')^2 Q^2(U_{\e '}) \leq ( 4 r_0 ^2 \cs ^3  \bS_{\rm min} )^2 $ because $ \e' \in I_n$. 
In particular, $ Q(U_{\e'})$ is bounded as $ \e' \in ( \e_n, \,  \e_n^{\#}).$
Since $ c( \e') \searrow c( \e_n^{\#}) $ as $ \e' \nearrow \e_n^{\#}$,  Proposition \ref{monoto} ($iv$) 
implies that $ c \longmapsto T_c$ is continuous at $c( \e_n^{\#})$. 
Then passing to $\ds \liminf $ as $ \e' \nearrow \e_n^{\#}$ in \eqref{estimate2} we get 
$(\e_n^{\#})^2 Q^2(U_n^{\#}) \leq  ( 4 r_0 ^2 \cs ^3  \bS_{\rm min} )^2 $.
We conclude that $ \e_n^{\#} \in I_n$. 

\medskip

Next, for any $ \e' \in ( \e_n^{\#}, \cs )$ and any minimizer $  U_{\e'} $ of 
$E_{c(\e')}$  on $ \Co _{c(\e')} $ that solves  (TW$_{c(\e')}$), inequality \eqref{estimate1} 
holds with $ \e_n^{\#} $ and $ U_n^{\#}$ instead of $ \e_n$ and $ U_n$, respectively. 
The limit of the right-hand side as $ \e' \searrow \e _n^{\#} $ is $(\e_n^{\#})^2 Q^2(U_n^{\#})$. 
If $\e_n^{\#} | Q(U_n^{\#} | <  4 r_0 ^2 \cs ^3  \bS_{\rm min} $, 
as above we infer that there is $ \de_n > 0 $ such that $[ \e_n^{\#}, \e_n^{\#} + \de _n ] \subset I_n$, 
contradicting the fact that $ \e_n^{\#} = \sup I_n$. 
The claim \eqref{claim1} is thus proved. 

\medskip

Now we turn our attention to the sequence $(\e_n^{\#}, U_n^{\#})_{n\geq 1}$. 
It is clear that $ \e_n ^{\#} \to 0 $ (because $ \e_n^{\#} \in (\e_n, \e_{n-1})$). 
By Proposition \ref{asympto} ($ii$) there is $ K > 0 $ such that 
$$
E ( U_n^{\#} ) + c(\e_n^{\#}) Q ( U_n^{\#} ) =  E_{c(\e_n^{\#})} ( U_n^{\#} )= T_{c(\e_n^{\#})} \leq K \e_n^{\#} 
$$
and using \eqref{claim1} we find $ | E ( U_n^{\#} )| \leq \frac{K'}{\e_n^{\#} }$ 
for some constant $ K'> 0 $ and for all $n$ sufficiently large. 
Hence we may use Proposition \ref{convergence} and we infer that there is a subsequence 
$(\e_{n_k}^{\#}, U_{n_k}^{\#})_{k\geq 1}$ which satisfies the conclusion of Theorem \ref{res1}. 
In particular, we have 
$$
\lim_{k \ra \infty} \e_{n_k}^{\#} | Q( U_{n_k}^{\#}) | = r_0 ^2 \cs ^3 \bS_{\rm min} 
$$
and this contradicts the fact that $U_{n_k}^{\#} $ satisfies \eqref{claim1}. 
Proposition \ref{global3} is thus proven.
\carre

\subsection{Proof of Proposition \ref{lifting}}
\label{preuvelifting}

($i$) Since $ U \in \BE$, we have $ |U| - r_0 \in H^1( \R^N)$ (see the Introduction of \cite{CM1})
and then $ \Big| \ds \frac{ \p }{\p x_i} ( |U | - r_0 ) \Big| \leq \Big| \frac{ \p U}{\p x_i} \Big| $ a.e. in $ \R^N$. 
It is well-known (see, for instance, \cite{brezis} p. 164) that for any $ \phi \in H^1( \R^N)$ there holds
$$
\| \phi \|_{L^{2^*}(\R^N)} \leq C_S \prod _{i =1}^N \Big\| \frac{\p \phi}{\p x_i} \Big\| _{L^2( \R^N)}^{\frac  1N}. 
$$
We infer that 
\be
\label{sobo}
\| \, |U| - r_0 \|_{L^{2^*}(\R^N)} \leq  C_S \prod _{i =1}^N \Big\| \frac{\p U}{\p x_i} \Big\| _{L^2( \R^N)}^{\frac  1N}
\leq C_S  \Big\| \frac{\p U}{\p x_1} \Big\| _{L^2( \R^N)}^{\frac  1N} \cdot 
 \| \nabla _{x_{\perp}} U \| _{L^2( \R^N)}^{\frac{N}{N-1}}.
\ee
Assume first that (A2) holds. 
If $\Big\| \frac{\p U}{\p x_1} \Big\| _{L^2( \R^N)} \cdot \| \nabla _{x_{\perp}} U \| _{L^2( \R^N)}^{N-1} \leq 1$, 
from \eqref{sobo} we get $ \| \, |U| - r_0 \|_{L^{2^*}(\R^N)} \leq C_S$. 
Let $ \tilde{U} (x) = e^{-\frac{i c x_1}{2}} U(x) $. 
Then $ \tilde{U} \in H_{loc}^1( \R^N)$ and $ \tilde{U}$ solves the equation 
$$
\Delta \tilde{U} + \left( \frac{c^2}{4} + F(|\tilde{U} |^2) \right)  \tilde{U} = 0 \qquad \mbox{ in }  \R^N. 
$$
Since $ \| \tilde{U} \|_{L^{2^*} (B(x, 1))} \leq C$ for any $ x \in \R^N$ and for some constant $C>0$, 
 using the above equation and a standard bootstrap argument
(which works thanks to (A2)), we infer that $ \| \tilde{U} \|_{W^{2,p} (B(x, \frac{1}{2^{n_0}}))} \leq \tilde{C}_p$
for some $ n_0 \in \N$, $ \tilde{C}_p  > 0$ and for any $ x \in \R^N$ and any $ p \in [2, \infty)$. This clearly  implies 
  $ \| {U} \|_{W^{2,p} (B(x, \frac{1}{2^{n_0}}))} \leq {C_p}$ for any $ x \in \R^N$ and any $ p \in [2, \infty)$.
In particular, using the Sobolev embedding we see that there is $ L>0 $ (independent on $U$) 
such that $ \| \nabla U \|_{L^{\infty}( \R^N)} \leq L$.

Fix $ \d > 0$. If there is $ x_0 \in \R^N$ such that $ | \, |U(x_0)|  - r_0 | \geq \d$, 
we infer that $ \| \, |U(x)|  - r_0 | \geq \frac{\d}{2}$ for any $ x \in B( x_0,  \frac{ \d}{2L})$ 
and consequently 
\be
\label{low}
\| \, |U| - r_0 \|_{L^{2^*}(\R^N)}  \geq \frac{\d}{2} \left(\Lo ^N \left(  B( x_0 \frac{ \d}{2L}) \right) \right)^{\frac{1}{2^*}}
= \frac{\d}{2}  \left(\frac{\d}{2L} \right)^{\frac{N}{2^*}}  \left(\Lo ^N(  B(0,1)) \right)^{\frac{1}{2^*}} . 
\ee
Let $ \mu ( \d) = \min \left(
1, \frac{\d}{2}  \left(\frac{\d}{2L} \right)^{\frac{N}{2^*}}  \left( \Lo ^N(  B(0,1)) \right)^{\frac{1}{2^*}}
\right).$
From \eqref{sobo} and \eqref{low} we infer that $|\, |U(x) | - r_0 | < \de $ for any solution $ U \in \BE$ of (TW$_c$) 
satisfying 
$ \Big\| \frac{\p U}{\p x_1} \Big\| _{L^2( \R^N)} \cdot \| \nabla _{x_{\perp}} U \| _{L^2( \R^N)}^{N-1} \leq  \mu ( \de).$

If (A3) holds, it follows from the proof of Proposition 2.2 p. 1078-1080 in \cite{M2} thet there is $ L >0$, independent on $U$, 
such that  $ \| \nabla U \|_{L^{\infty}( \R^N)} \leq L$. 
The rest of the proof is as above. 

\medskip

($ii$) By Proposition 2.2 p. 1078 in \cite{M2} we know that $ U \in W_{loc}^{2,p}(\R^N)$ for any $ p \in [2, \infty)$. 
In particular, $U \in C^1( \R^N)$ . 
As in the proof of ($i$) we see that there is $ L > 0$, independent on $U$, such that $ || \nabla U \|_{L^{\infty}( \R^N)} \leq L$. 

Fix $ \d > 0$ and assume that there is $ x^0 = ( x_1^0, \dots, x_N^0)$ such that $ | \, |U(x^0)|  - r_0 | \geq \d$.
Then we have  $ | \, |U(x)|  - r_0 | \geq \frac{\d}{2}$ for any $ x \in B( x^0 , \frac{ \d}{2L})$ 
and, in particular, 
$|\, | U(x_1, x_2^0, \dots, x_N^0)| - r_0 | 
\geq \frac{ \d}{2}$ for any $ x_1 \in [x_1^0 - \frac{ \d}{2L}, x_1^0 + \frac{ \d}{2L}]$.
We infer that 
$|\, | U(x_1, x_{\perp})| - r_0 | \geq \frac{\d}{4} $ for any $ x_1 \in [x_1^0 - \frac{ \d}{2L}, x_1^0 + \frac{ \d}{2L}]$ 
and any $ x_{\perp} \in B_{\R^{N-1}} ( x_{\perp}^0, \frac{ \d}{4L}).$ 
Consequently
$$
\begin{array}{l}
\| \, |U( x_1, \cdot ) | - r_0 \|_{L^{\frac{2(N-1)}{N-3}}(\R^{N-1})}   \geq 
\frac{ \d}{4} \left( \Lo^{N-1} \left( B_{\R^{N-1}}\left(x_{\perp}^0, \frac{\d}{4L}\right) \right) \right)^{\frac{N-3}{2(N-1)}}
\\ \\
 \geq  \frac{ \d}{4}  \left(\frac{\d }{4L} \right)^{\frac{N-3}{2}} \left( \Lo^{N-1} ( B_{\R^{N-1}}(0, 1)) \right)^{\frac{N-3}{2(N-1)}} 
= K \d ^{ \frac{N-1}{2}}
\end{array}
$$
for all $ x_1 \in [x_1^0 - \frac{ \d}{2L}, x_1^0 + \frac{ \d}{2L}]$.
Using the Sobolev inequality in $ \R^{N-1}$ we get for  $  x_1 \in \left[x_1^0 - \frac{ \d}{2L}, x_1^0 + \frac{ \d}{2L}\right], $
$$
\int_{\R^{N-1}} |\nabla_{x_{\perp}} U( x_1, x_{\perp}) |^2 \, dx _{\perp}
\geq \frac{1}{\tilde{C}_S^2} \| \, |U( x_1, \cdot ) | - r_0 \|_{L^{\frac{2(N-1)}{N-3}}(\R^{N-1})} ^2
\geq \frac{K^2}{\tilde{C}_S^2} \de ^{N-1} .
$$
Integrating the above inequality on $[x_1^0 - \frac{ \d}{2L}, x_1^0 + \frac{ \d}{2L}]$ we obtain 
$ \| \nabla_{x_{\perp} } U \|_{ L^2 (\R^N )} ^2 \geq \frac{K^2}{L \tilde{C}_S^2} \d ^{N} = K_1  \d ^{N} .$
We conclude that if $\| \nabla_{x_{\perp} } U \|_{ L^2 (\R^N )} ^2  <\min(1,  K_1  \d ^{N}) $, then necessarily 
$  | \, |U|  - r_0 | < \d$ in $ \R^N$. 
\carre

\subsection{Proof of Proposition \ref{prop2d}}


It follows from Lemma 4.1 in \cite{CM1} 
that there are $ k_0 > 0$, $ C_1, C_2 > 0$ such that for all $ \psi \in \BE$ with $\ds \int_{\R^2} | \nabla \psi |^2 \, dx \leq k_0$
we have
\be
\label{kifkifpot}
C_1 \int_{\R^2} ( \chi^2(|\psi|) - r_0^2 )^2 \ dx \leq 
\int_{\R^2} V(|\psi|^2) \ dx  \leq C_2 \int_{\R^2} ( \chi^2(|\psi|) - r_0^2 )^2 \ dx .
\ee

We recall that in space dimension two,  nontrivial solutions $U_k$ to 
(TW$_c$) have been constructed in Theorem  \ref{th2d}  by considering the minimization problem 
$$
\mbox{ minimize } I(\psi) = Q( \psi ) + \int_{\R^2} V(| \psi |^2) \, dx \quad \mbox{ in } \BE \; \mbox{ under the constraint } 
\int_{\R^2} | \nabla \psi |^2 \, dx = k.
\eqno{(\Io _k)}
$$
If $\BU_k$ is a minimizer for $(\Io _k)$, there is $ c_k > 0$ such that $U_k = (\BU_k)_{c_k, c_k}$  
solves  (TW$_{c_k}$) and minimizes $E_{c_k} = E + c_k Q$ in the set 
$ \Big\{ \psi \in \BE \; \Big| \; \ds \int_{\R^2} | \nabla \psi |^2 \, dx = k \Big\}$.
Moreover, we have $ c_k \to \cs $ as $ k \to 0$. 
Lemma \ref{liftingfacile} implies that $|U_k| \to r_0$ uniformly on $ \R^2$ as $ k \to 0$; 
in particular, there is $ k_1 > 0 $ such that if $ k \in ( 0, k_1) $, 
we have $|U_k| \geq \frac{r_0}{2} $ in $ \R^2$. 
From the Pohozaev identities \eqref{Pohozaev} we get $c_k Q( U_k) + 2 \ds \int_{\R^2} V(|U_k|^2 )\, dx = 0$, and this gives 
\be
\label{scaling2}
I_{\rm min}(k) = I( \BU _k) = \frac{ 1}{c_k} Q(U_k) + \frac{ 1}{c_k^2 } \int_{\R^2} V(|U_k|^2 )\, dx 
= \frac{ 1}{2c_k} Q(U_k) = - \frac{ 1}{c_k^2 } \int_{\R^2} V(|U_k|^2 )\, dx.
\ee
By Lemma 5.2 in \cite{CM1} there is $ k_2 > 0 $ such that  $  - \frac{2k}{\cs ^2} \leq I_{\rm min}(k) \leq - \frac{k}{\cs ^2} $ 
for all $ k \in (0, k_2)$. 
Since $ c_k \to \cs $ as $ k \to 0$, the estimates \eqref{estim2d} follow directly from \eqref{kifkifpot} and \eqref{scaling2}.

\medskip

It remains to prove \eqref{kifkif}.
By Proposition \ref{asympto}, there is $ \mu _0 > 0 $ such that for $k$ sufficiently small we have 
$ I_{\rm min} (k) \leq  - \frac{k}{\cs^2} - \mu_0 k^3 .$
 By  scaling we  have
$$ \frac{1}{c_k^2} \Big( E_{c_k}(U_k) - \int_{\R^2} |\nabla U_k|^2 \ dx  \Big) 
= \frac{1}{c_k^2} \Big( c_k Q(U_k) + \int_{\R^2} V(|U_k|^2) \ dx  \Big) 
= I(\BU_k) = I_{\rm min}(k) \leq - \frac{k}{\cs^2} - \mu_0 k^3 . $$
Since $\cs^2 - c_k^2 = \e _k^2 $ and $\ds \int_{\R^2} |\nabla U_k|^2 \ dx  = k$, we get
\be
\label{ecusson}
E_{c_k}(U_k) \leq k \Big( 1 - \frac{c _k^2}{\cs^2} \Big) - \mu_0 c_k^2 k^3 
= \frac{k\e_k^2}{\cs^2} - \mu_0 c_k ^2 k^3 .
\ee
The second Pohozaev identity \eqref{Pohozaev} yields $E_{c_k}(U_k) = 
2 \ds \int_{\R^2} |\p_{2} U_k |^2 \ dx \geq 0$, thus
$ 0 \leq k \Big( \frac{\e_k^2}{\cs^2} - \mu_0 c_k ^2 k^2 \Big) $
and this implies
$$ \frac{\e_k^2}{\cs^2} \geq  \mu_0 c^2 k^2 . $$
Since $c \geq \cs/2$ for $k$ small, 
the left-hand side inequality in \eqref{kifkif} follows.

\medskip

In order to prove the second inequality in \eqref{kifkif}, we need the next Lemma. 
In the case of the Gross-Pitaevskii nonlinearity, this result follows from  Lemma 2.12 p. 597 in \cite{BGS1}.
In the case of general nonlinearities, it was proved in \cite{CM1}.

\begin{lem} [\cite{BGS1, CM1}]
\label{tools}
Let $N\geq 2$. 
There is $ \beta _* > 0$ such that  any solution $ U = \rho e^{ i \phi} \in \Eo$ of (TW$_c$)  verifying
$ r_0 - \beta_* \leq \rho \leq r_0 + \beta _*$ 
satisfies the identities 
\be
\label{blancheneige}
 E(U) + c Q(U) = \frac{2}{N} \int_{\R^N} |\nabla \rho |^2 \ dx \qquad \mbox{ and } 
\ee
\be
\label{grincheux}
2 \int_{\R^N} \rho^2 |\nabla \phi |^2 \ dx = 
c \int_{\R^N} ( \rho^2 - r_0^2 ) \p_1 \phi \ dx = - c Q(U) .
\ee
Furthermore, there exist $a_1, a_2 > 0$ such that
\be
\label{balai}
a_1 \|  \rho ^2 - r_0 ^2 \|_{L^2(\R^N)} \leq \| \nabla U \|_{L^2( \R^N)} \leq a_2 \|  \rho ^2 - r_0 ^2 \|_{L^2(\R^N)}.
\ee
\end{lem}

\noindent {\it Proof.}
Identity \eqref{grincheux} is Lemma 7.3 ($i$) in \cite{CM1}. 
Formally, it follows by 
multiplying the first equation in \eqref{phasemod} by $ \phi $ and 
integrating by parts over $\R^N$; see \cite{CM1} for a rigorous justification.

Combining the two Pohozaev identities in \eqref{Pohozaev}, we have 
$$ (N-2) \int_{\R^N} |\nabla U|^2 \ dx + N \int_{\R^N} V(|U|^2) \ dx 
+ c ( N-1 ) Q(U) = 0 . $$
Using that $|\nabla U|^2 = |\nabla \rho|^2 + \rho^2 |\nabla \phi|^2 $, 
we infer from \eqref{grincheux}
\begin{align*}
N(E(U) + c Q(U) ) = 2 \int_{\R^N} |\nabla U|^2 \ dx + c Q(U) 
= & \ 2 \int_{\R^N} |\nabla \rho|^2 \ dx 
+ \Big( 2 \int_{\R^N} \rho^2 |\nabla \phi|^2 \ dx + c Q(U) \Big) 
\\ = & \ 2 \int_{\R^N} |\nabla \rho|^2 \ dx ,
\end{align*}
and this establishes  \eqref{blancheneige}.
The estimate \eqref{balai} has been proven in \cite{CM1} (see inequality (7.17) there).
\carre

\medskip

We come back to the proof of Proposition \ref{prop2d}.
We write $ U_k = \rho e^{ i \phi}$ and we denote $ \eta = \rho ^2 - r_0 ^2$, so that $ \rho$, $\phi $ and $ \eta$ satisfy 
\eqref{phasemod}$-$\eqref{fond} (with $c_k$ instead of $c$). 
Taking the Fourier transform of \eqref{fond} we get
\be 
\label{fondfou}
\begin{array}{rcl}
\wh{\eta} ( \xi) & = & \ds \frac{  |\xi|^2}{|\xi |^4 + \cs ^2 |\xi |^2 - c_k ^2 \xi _1^2 } 
\Fou \left( -2 |\nabla U_k|^2 + 2 c_k \eta \frac{ \p \phi}{\p x_1} + 2 \rho ^2 F( \rho ^2) + \cs ^2 \eta  \right) 
\\
\\
& & 
\ds - 2 c_k \sum_{j=1}^N \frac{ \xi _1 \xi _j }{|\xi |^4 + \cs ^2 |\xi |^2 - c_k ^2 \xi _1^2 }
 \Fou \left( \eta \frac{ \p \phi}{\p x_j} \right).
\end{array}
\ee
It is easy to see that $ 2 \rho ^2 F( \rho ^2) + \cs ^2 \eta  = \BO ( (\rho ^2 - r_0 ^2) ^2) = \BO ( \eta^2)$, hence
$$
\| \Fou \left( 2 \rho ^2 F( \rho ^2) + \cs ^2 \eta \right) \|_{L^{\infty} ( \R^N)} \leq \| 2 \rho ^2 F( \rho ^2) + \cs ^2 \eta \|_{L^1( \R^N)}
\leq C \| \eta \|_{L^2( \R^N)} ^2. 
$$
Since $ r_0 - \beta_* < |U_k| < r_0 + \beta_* $ if $ k $ is sufficiently small and $ |\nabla U_k|^2 = |\nabla \rho |^2 + \rho^2 |\nabla \phi |^2$, 
 using \eqref{balai} we get 
$$
\Big\| \Fou \left( \eta \frac{ \p \phi}{\p x_j} \right) \Big\| _{L^{\infty}(\R^N)} 
\leq \Big\| \eta \frac{ \p \phi}{\p x_j} \Big\|  _{L^{1}(\R^N)} 
\leq \| \eta \|_{L^2(\R^N)} \Big\| \frac{ \p \phi}{\p x_j} \Big\|  _{L^{2}(\R^N)} 
\leq C  \| \eta \|_{L^2(\R^N)} ^2
$$
and $ \| \Fou ( |\nabla U_k| ^2) \|_{L^{\infty}(\R^N)} \leq \| \nabla U_k \|_{L^2(\R^N)}^2 \leq C \| \eta \|_{L^2(\R^N)} ^2 .$
Coming back to \eqref{fondfou} we discover
$$
| \wh{ \eta } ( \xi ) | \leq C  \| \eta \|_{L^2(\R^N)} ^2 \cdot \frac{  |\xi|^2}{|\xi |^4 + \cs ^2 |\xi |^2 - c_k ^2 \xi _1^2 } .
$$
Using Plancherel's formula and the above estimate we find 
\be
\label{planchereta}
\| \eta \|_{L^2(\R^N)} ^2 = \frac{1}{(2 \pi )^N} \int_{\R^N} |\wh{\eta} ( \xi )|^2\, d \xi 
\leq C \| \eta \|_{L^2(\R^N)} ^4 \int_{\R^N}  \frac{  |\xi|^4}{ (|\xi |^4 + \cs ^2 |\xi |^2 - c_k ^2 \xi _1^2) ^2 } \, d \xi.
\ee
If $N=2$, a straightforward computation using polar coordinates gives (see the proof of (2.59) p. 598 in \cite{BGS2}):
$$
\int_{\R^2}  \frac{  |\xi|^4}{ (|\xi |^4 + \cs ^2 |\xi |^2 - c_k ^2 \xi _1^2) ^2 } \, d \xi
= \frac{ \pi}{\cs \sqrt{ \cs ^2 - c_k^2}}  = \frac{ \pi}{\cs \e _k}.
$$
From to \eqref{planchereta} we get $ \| \eta \|_{L^2(\R^2)} ^2 \leq \frac{C}{\e _k}   \| \eta \|_{L^2(\R^2)} ^4$
and taking into account \eqref{balai} we infer that 
$\e _k \leq C  \| \eta \|_{L^2(\R^2)} ^2 \leq \tilde{C} \| \nabla U_k \|_{L^2( \R^2)} ^2 = \tilde{C} k. $
\carre

\medskip

Notice that at this stage, we have only upper bounds on the energy of travelling waves,  
and  we will have to prevent convergence towards the trivial solution to (SW). 
This will be done with the help of the following result. It was proven  in \cite{BGS2} in the case 
of the Gross-Pitaevskii nonlinearity (see Proposition 2.4 p. 595 there). We extend the proof 
to general nonlinearities. 

\begin{lem}
\label{minoinf} 
Let $N \geq 2$ and assume that (A1) holds and $F$ is twice differentiable at $ r_0^2$. 
There is $C>0$, depending only on $N$ and on $F$, such that any travelling wave $ U \in \BE$ of {\rm (NLS)} 
of speed $ c \in [0, \cs]$ such that $ \frac{r_0}{2} \leq|U| \leq \frac{ 3 r_0}{2}$ satisfies 
$$
\| \, |U| - r_0 \|_{L^{\infty}(\R^N)} \geq C( \cs ^2 - c^2) = C \e^2 (U).
$$
\end{lem}

\noindent {\it Proof.}
Let $U \in \BE$ be a travelling wave such that $ \frac{r_0}{2} \leq|U| \leq \frac{ 3 r_0}{2}$ in $ \R^N$. 
Then $U \in W_{loc}^{2,p}(\R^N)$, $ \nabla U \in W^{1,p}(\R^N)$ for all $ p \in [2, \infty)$ 
(see Proposition 2.2 p. 1078-1079 in \cite{M2}), and $U$ admits a lifting $U = \rho e^{i \phi}$, where 
$ \rho$ and $ \phi$ satisfy \eqref{phasemod}. 
Since $ U \in \BE$ we have $ \rho ^2 - r_0^2 \in H^1( \R^N)$ and then it is easy to see that 
$ \frac{ \rho ^2 - r_0^2}{\rho } \in H^1( \R^N)$. 
Multiplying the second equation in \eqref{phasemod} by $ \frac{ \rho ^2 - r_0^2}{\rho } $ 
and integrating by parts we get 
\be 
\label{ident1}
\int_{\R^N} \left( 1 + \frac{r_0^2}{\rho^2} \right) |\nabla \rho |^2 \, dx 
+ \int_{\R^N} ( \rho ^2 - r_0^2) |\nabla \phi |^2 - ( \rho ^2 - r_0^2) F( \rho ^2) - c ( \rho ^2 - r_0^2) \frac{ \p \phi }{\p x_1} \, dx = 0.
\ee
Denote $ \d = \| \, |U| - r_0 \|_{L^{\infty}(\R^N)} =  \| \rho  - r_0 \|_{L^{\infty}(\R^N)} .$ We have 
\be
\label{inegradrho}
\int_{\R^N} \left( 1 + \frac{r_0^2}{\rho^2} \right) |\nabla \rho |^2 \, dx  
\geq  \left( 1 + \frac{r_0^2}{(r_0 + \d )^2} \right) \int_{\R^N} |\nabla \rho |^2 \, dx   
\qquad \mbox{ and } 
\ee
\be
\label{majophase}
\Big|  \int_{\R^N} ( \rho ^2 - r_0^2) |\nabla \phi |^2 \, dx  \Big|
\leq   \int_{\R^N} \frac{| \rho ^2 - r_0^2 |}{\rho^2} \rho ^2 |\nabla \phi |^2 \, dx
\leq  \frac{ 2 r_0 \d + \d ^2 }{(r_0 - \d)^2}  \int_{\R^N} |\nabla U |^2 \, dx.
\ee
There is $ \tilde{C} > 0 $ such that $| F( s^2) - F'( r_0 ^2) ( s^2 - r_0 ^2) | \leq \tilde{C} ( s^2 - r_0 ^2)^2 $ 
for all $ s \in [\frac{r_0}{2}, \frac{ 3r_0}{2}]$. 
Remember  that $ - F'( r_0 ^2) = 2 a^2 $ and $ \cs = 2 a r_0$, thus 
\be
\label{60}
- ( \rho ^2 - r_0 ^2) F( \rho ^2) \geq - F'( r_0 ^2) ( \rho ^2 - r_0 ^2) ^2 - \tilde{C} | \rho ^2 - r_0 ^2 |^3 
\geq \left( 2 a^2 - \tilde{C} ( 2 r _0 \de + \de ^2 ) \right) (\rho ^2 - r_0 ^2) ^2. 
\ee
Using \eqref{grincheux} and \eqref{momentlift}, then \eqref{ident1} and \eqref{inegradrho}-\eqref{60} we get 
$$
\begin{array}{l}
\ds - 2 c Q( U) = 2 \int_{\R^N} \rho ^2 |\nabla \phi |^2 \, dx + c \int_{\R^N} ( \rho ^2 - r _0 ^2) \frac{ \p \phi}{\p x_1 } \, dx 
\\
\\
\ds =  2 \int_{\R^N} \rho ^2 |\nabla \phi |^2 \, dx + 
\int_{\R^N} \left( 1 + \frac{r_0^2}{\rho^2} \right) |\nabla \rho |^2 \, dx 
+ \int_{\R^N} ( \rho ^2 - r_0^2) |\nabla \phi |^2 - ( \rho ^2 - r_0^2) F( \rho ^2) \, dx 
\\
\\
\ds \geq 
 2 \int_{\R^N} \rho ^2 |\nabla \phi |^2 \, dx  
+  \int_{\R^N} \left( 1 + \frac{r_0^2}{(r_0 + \d )^2} \right) |\nabla \rho |^2   
- \frac{ 2 r_0 \d + \d ^2 }{(r_0 - \d)^2}  |\nabla U |^2 
+  \left( 2 a^2 - \tilde{C} ( 2 r _0 \de + \de ^2 ) \right)  (\rho ^2 - r_0 ^2) ^2 \, dx
\end{array}
$$
and we infer that there exists $K>0$, depending only on $F$, such that 
\be
\label{61}
 -2 c Q(U) \geq 2 ( 1 - K \d ) \int_{\R^N} |\nabla U |^2 + a^2 ( \rho ^2 - r_0 ^2)^2\, dx.
\ee
On the other hand, using \eqref{momentlift} we have 
\be 
\label{62}
\begin{array}{l}
\ds - Q( U) = \frac{ 2 a r_0}{ \cs} \int_{\R^N} ( \rho ^2 - r _0 ^2) \frac{ \p \phi}{\p x_1 } \, dx 
\leq \frac{1}{\cs}  \int_{\R^N}  r_0 ^2  \Big| \frac{ \p \phi}{\p x_1 } \Big| ^2 + a ^2 ( \rho ^2 - r_0 ^2)^2\, dx
\\
\\
\ds \leq \frac{1}{\cs} \int_{\R^N} \frac{ r_0 ^2}{(r_0 - \d )^2}   \rho ^2  \Big| \frac{ \p \phi}{\p x_1 } \Big| ^2  
+ a^2  ( \rho ^2 - r_0 ^2)^2\, dx 
\leq  \frac{1}{\cs}   \frac{ r_0 ^2}{(r_0 - \d )^2}  \int_{\R^N} |\nabla U |^2 + a^2 ( \rho ^2 - r_0 ^2)^2\, dx.
\end{array}
\ee
Since $U$ is not constant we have $ \ds  \int_{\R^N} |\nabla U |^2 + a^2 ( \rho ^2 - r_0 ^2)^2\, dx > 0$ 
and comparing \eqref{61} and \eqref{62}   we get 
$$
\frac{c}{\cs }   \frac{ r_0 ^2}{(r_0 - \d )^2}  \geq 1 - K \d. 
$$
If $ \d > \frac{1}{2K}$ the conclusion of Lemma \ref{minoinf} holds because $ \e (U)$ is bounded. 
Otherwise, the previous inequality is equivalent to 
$ \frac{ r_0 ^2}{(r_0 - \d )^2}  \frac{1}{ 1 - K \d} \geq \frac{ \cs}{\sqrt{ \cs ^2 - \e^2(U) }}. $
There are $ K_1, \; K_2 > 0 $ such that $  \frac{ r_0 ^2}{(r_0 - \d )^2}  \frac{1}{ 1 - K \d}  \leq 1 + K_1 \d $ and 
$ \frac{ \cs}{\sqrt{ \cs ^2 - \e^2 }} \geq 1 + K_2 \e^2 $ for all $ \d \in [0, \frac{1}{2K} ] $ and all $ \e \in [0, \cs )$ and we 
infer that $ 1 + K_1 \d \geq 1 + K_2 \e ^2(U)$, that is 
$ \d = \| \, | U| - r_0 \|_{L^{\infty} ( \R^N )} \geq \frac{K_2}{K_1} \e^2(U)$.
\carre

\subsection{Initial bounds for  $\bs{\BA_{\e}}$}

Let $U_c \in \BE$ be a travelling wave to (NLS) of speed $ c$ 
provided by Theorems \ref{th2dposit} or \ref{th2d} if $N=2$, respectively by Theorem \ref{thM} if $N=3$, 
such that $\frac{ r_0}{2} \leq |U| \leq \frac{ 3r_0}{2} $ in $ \R^N$. 
As in \eqref{ansatz}, we write $U_c(x) = \rho(x) e^{ i \phi(x)} = r_0 \sqrt{1+\e^2 \BA_{\e}(z) }\ \ex^{i\e \vp_{\e} (z)}, $ where 
$ \e = \sqrt{ \cs ^2 - c^2}, $ $z_1 = \e x_1 ,$  $ \ z_\perp = \e^2 x_\perp .$
According to Proposition 2.2 p. 1078-1079 in \cite{M2} we have 
$$
\| U_c \|_{C_b^1( \R^N)} \leq C \qquad \mbox{ and } \qquad \| \nabla U_c \|_{W^{1,p}(\R^N)} \leq C_p 
\quad \mbox{ for } p \in [2, \infty) .
$$
By scaling, we obtain the initial (rough) estimates
\be
\label{bourrinska}
\| \BA _{\e} \|_{L^{\infty }} \leq \frac{C}{\e ^2}, \quad 
\|  \p _{z_1} \BA _{\e} \|_{L^{\infty }} \leq \frac{C}{\e ^3}, \quad 
\| \nabla _{z_{\perp}} \BA _{\e} \|_{L^{\infty }} \leq \frac{C}{\e ^4}, \quad 
\| \p_{z_1} \vp _{\e} \|_{L^{\infty }} \leq \frac{C}{\e ^2}, \quad 
\| \nabla _{z_{\perp}} \vp_{\e}  \|_{L^{\infty }} \leq \frac{C}{\e ^3}
\ee
and 
\be 
\label{bourrinSKF}
\Big\| \frac{ \p ^2 \BA_{\e}}{\p  z_1 ^2} \Big\| _{L^p} \leq C_p \e^{-4 + \frac{2N-1}{p}}, \qquad
\Big\| \frac{ \p ^2 \BA_{\e}}{\p  z_1 \p z_j } \Big\| _{L^p} \leq C_p \e^{-5 + \frac{2N-1}{p}}, \qquad
\Big\| \frac{ \p ^2 \BA_{\e}}{\p  z_j \p z_k} \Big\| _{L^p} \leq C_p \e^{-6 + \frac{2N-1}{p}}
\ee
for any $ p \in [2, \infty) $  and all $ j, k \in \{ 2, \dots, N \}.$
We have: 

\begin{lem}
\label{BornEnergy} 
Assume that (A2) and (A4) are satisfied and $ \G \neq 0$. 
Let $U_c$ be a solution to {\rm (TW$_{c}$)} provided by Theorem \ref{th2d} if $N=2$, 
respectively by Theorem \ref{thM} if $N=3$ and let $ \e = \sqrt{\cs ^2 - c^2}$. 
If $N=3$ we assume moreover that $E(U_c) \leq \frac{K}{\e}$, where $K$ does not depend on $ \e$. 

There exist $ \e _0 > 0 $ and $ C > 0 $ (depending only on $F$, $N$, $K$) such that  
$U_c$ admits a lifting as in \eqref{ansatz} whenever $ \e \in (0, \e _0)$ and the following estimate holds:  
$$
 \int_{\R^N} | \p_{z_1} \vp_{\e} |^2 + |\nabla_{z_\perp} \vp_{\e}|^2 
+ \BA_{\e} ^2 + | \p_{z_1} \BA_{\e }  |^2 
+ \e^2 | \nabla_{z_\perp} \BA_{\e} |^2 \ dz \leq C . $$
\end{lem}

\noindent {\it Proof.} 
If $ N =2$ it follows from Theorem \ref{th2d} that 
$ k = \ds \int_{\R^2} |\nabla U_c |^2 \, dx $ is small if $ \e $ is small.  
Using Lemma \ref{liftingfacile} in the case $N=2$, respectively Corollary \ref{sanszero} if $N=3$, 
we infer that $|U_c|$ is arbitrarily close to $ r_0 $ if $ \e $ is sufficiently small 
and then it is clear that we have a lifting as in \eqref{ansatz}.

\medskip

 We will  repeatedly use the fact that there is a constant $C$ depending only on $F$ such that 
$$
 C |\p_j U_c|^2 \geq |\p_j (\rho^2)|^2 + |\p_j \phi|^2  \qquad \mbox{ for } 1 \leq j \leq N.
 $$
In view of the Taylor expansion of $V$ near $r_0^2$, for $\e$ sufficiently close to $ 0$ (so that $|U_c|$ is sufficiently close to $r_0$) we have
$$ V(|U_c|^2) \geq C (|U_c| - r_0)^2 . $$
By scaling, we infer that  for some $\d_1 >0$ depending only on $F$ there holds  
$$ E(U_c) = \int_{\R^N} |\nabla U_c|^2 + V(|U_c|^2) \ dx 
\geq \d_1 \e^{5-2N} \int_{\R^N} | \p_{z_1} \vp_{\e} |^2 + \BA_{\e}^2 \ dz . $$

In the case $ N =2$ it follows from  Proposition \ref{prop2d}  
that $E(U_c) \leq C \e $ for some $ C $ independent of $ \e $. 
In the case $ N =3$ we use the assumption  $E(U_c) \leq \frac{K}{\e}$.  
In both cases the previous inequality implies that 
\be
\label{bourrin2}
\int_{\R^N} | \p_{z_1} \vp_{\e} |^2 + \BA_{\e}^2 \ dz 
\leq C .
\ee
We have   $E_{c} (U_c) = T_{c} = \BO(\e )$ if $ N=3 $ by  Proposition \ref{asympto} $(ii)$, 
respectively  $E_{c} (U_c) = \BO( k  \e ^2 ) = \BO(\e ^3)$ by  \eqref{ecusson} and \eqref{kifkif} in the case $ N =2$. 
From the Pohozaev identity $P_c(U_c) = 0$ (see \eqref{Pohozaev})  we deduce
$$
 \frac{2r_0^2 \e^{7-2N}}{N-1} \int_{\R^N} 
|\nabla_{z_\perp} \vp_{\e}|^2 + \e ^2 |\nabla_{z_\perp} \BA_{\e}|^2 \ dz 
\leq C \frac{2}{N-1} \int_{\R^N} |\nabla_{\perp} U_c|^2 \ dx 
= C E_{c} (U_c) = \BO(\e ^{7-2N}) .
 $$
Thus  we get 
\be
\label{bourrin3}
\int_{\R^N} |\nabla_{z_\perp} \vp_{\e}|^2 + \e ^2 |\nabla_{z_\perp} \BA_{\e}|^2\ dz 
\leq C . 
\ee
Furthermore, by scaling the identity \eqref{blancheneige} in 
Lemma \ref{tools}  we obtain
$$ r_0^2 \e ^{7-2N} \int_{\R^N} |\p_{z_1} \BA_{\e } |^2 \ dz
\leq C \int_{\R^N} |\p_{x_1}  \rho  |^2 \ dx \leq C 
\int_{\R^N} |\nabla \rho  |^2 \ dx = C \frac{N}{2} E_{c } (U_c) = 
\BO (\e ^{7-2N} ) , $$
so that
\be
\label{bourrin4}
\int_{\R^N} |\p_{z_1} \BA_ {\e} |^2 \ dz \leq C .
\ee
Gathering \eqref{bourrin2}, \eqref{bourrin3} and \eqref{bourrin4} 
yields the desired inequality. \carre \\

\medskip

Using the above estimates, we shall find  $L^q$ bounds 
for $\BA_{\e}$. The proof is based 
on equation \eqref{Fonda}, that is
$$
\Big\{ \p_{z_1}^4 - \p_{z_1}^2 - \cs^2 \Delta_{z_\perp} 
+ 2 \e ^2 \p_{z_1}^2 \Delta_{z_\perp} + \e ^4 \Delta^2_{z_\perp} \Big\} 
\BA_{\e} =   \BR_{\e} ,
\eqno{\mbox{\eqref{Fonda}}}
$$
where
\begin{align*}
\BR _{\e} = & \ 
\{ \p_{z_1}^2 + \e ^2 \Delta_{z_\perp} \} \Big[ 
2(1 + \e^2 \BA _{\e}) \Big( (\p_{z_1} \vp {\e} )^2 
+ \e ^2 |\nabla_{z_\perp} \vp _{\e} |^2 \Big) + \e^2 \frac{(\p_{z_1} \BA _{\e} )^2 
+ \e ^2 |\nabla_{z_\perp} \BA _{\e} |^2}{2(1+ \e ^2 \BA _{\e} )} \Big]
\\ & \ 
- 2 c   \e ^2 \Delta_{z_\perp} ( \BA _{\e}  \p_{z_1} \vp _{\e}) 
+ 2 c   \e ^2 \ds \sum_{j=2}^N \p_{z_1} \p_{z_j} ( \BA _{\e} \p_{z_j} \vp _{\e} )
\\ & \ 
+ \{ \p_{z_1}^2 + \e ^2 \Delta_{z_\perp} \} \Big[ 
\cs^2 \Big( 1 - \frac{r_0^4F''(r_0^2)}{\cs^2} \Big) \BA _{\e}^2 
- \frac{1}{\e ^4} \tilde{F}_3(r_0^2 \e ^2 \BA _{\e}) \Big] 
\end{align*}
and  we   recall that $ \tilde{F}_3(\a) = \BO(\a^3)$ as $\a \to 0$. 

Let 
$$ 
D_{\e} (\xi) = \xi_1^4 + \xi_1^2 + \cs^2 |\xi_\perp|^2 
+ 2 \e^2 \xi_1^2 |\xi_\perp|^2 + \e^4 |\xi_\perp|^4 
= ( \xi _1 ^2 + \e ^2|\xi _{\perp} |^2)^2 + \xi _1 ^2 + \cs ^2 |\xi _{\perp}|^2.
 $$
We will consider  the following kernels:
$$ \BK^1_{\e }  (z) = \Fou^{-1} \Big( \frac{\xi_1^2}{D_{\e} (\xi)} \Big) , 
\quad \quad \BK^\perp_{\e} (z) = 
\Fou^{-1} \Big( \frac{|\xi_\perp|^2}{D_{\e}(\xi)} \Big) 
\quad \quad  {\rm and} \quad \quad 
\BK^{1,j}_{\e} (z) = \Fou^{-1} \Big( \frac{\xi_1 \xi_j}{D_{\e}(\xi)} \Big) , \quad j = 2, \dots, N.$$
Then we may rewrite \eqref{Fonda}  as a convolution equation
\be
\label{Henry}
\BA_{\e} = \Big( \BK^1_{\e} + \e ^2 \BK^\perp_{\e} \Big) * G_{\e} 
+ 2 c \e ^2 \BK_{\e}^{\perp} * (\BA_{\e} \p_{z_1} \vp_{\e}) 
- 2 c(\e) \e ^2 \sum_{j=2}^N \BK^{1,j}_{\e} * 
(\BA_{\e} \p_{z_j} \vp_{\e}) ,
\ee
where
\begin{align*}
G_{\e}  = & \ 
(1 + \e^2 \BA_{\e}) \Big( (\p_{z_1} \vp_{\e})^2 
+ \e ^2 |\nabla_{z_\perp} \vp_{\e}|^2 \Big) 
+ \e ^2 \frac{(\p_{z_1} \BA_{\e} )^2 + \e^2 |\nabla_{z_\perp} \BA_{\e}|^2}
{4(1+ \e ^2 \BA_{\e}  )} \\ & \ 
+ \frac{\cs^2}{4} ( \G - 2 ) \BA_{\e} ^2 
- \frac{1}{\e ^4} \tilde{F}_3(r_0^2 \e ^2 \BA_{\e} ) .
\end{align*}

\begin{lem}
\label{Grenouille} 
The following estimates hold for $N=2$, $3$ and $\e$ small   enough:

\medskip

 (i) For all $ 2 \leq p \leq \ii $ we have $ \n \p_{z_1} \BA_{\e}  \n_{L^p} 
+ \e  \n \nabla_{z_\perp} \BA_{\e}  \n_{L^p} 
\leq C \e ^{\frac{6}{p}-3} $.

\medskip

 (ii)  There exists $C>0 $ such that  
$ \n \BA_{\e} \n_{L^{3q}} \leq C \e ^{- \frac23 } 
\n \BA_{\e}  \n^{\frac23}_{L^{2q}}$  for any $ 1 \leq q \leq \infty$.

\medskip

 (iii)  If $N=3$,  for any $ 2 \leq p < 8/3 $ there is  $ C_p > 0 $ such that 
$ \n \BA_{\e}  \n_{L^p(\R^3)} \leq C_p $. 

\medskip

(iv)  If $N=2$,  for any 
$ 2 \leq p < 4 $ there is $ C_p > 0 $ such that $ \n \BA_{\e} \n_{L^p(\R^2)} \leq C_p $.
\end{lem}

\noindent {\it Proof.} For $(i)$, it suffices to notice that the 
estimate is true for $p=2$ by Lemma \ref{BornEnergy} and for $p=\ii$ 
by \eqref{bourrinska}, therefore it holds  for any $2 \leq p \leq \ii$ by interpolation. 
For $(ii)$ we just interpolate the exponent $3q$ between $2q$ and $\infty$ and we use \eqref{bourrinska}:
$$ 
\n \BA_{\e}  \n_{L^{3q}} \leq 
\n \BA_{\e} \n_{L^{2q}}^{ \frac 23 } 
\n \BA_{\e}  \n_{L^{\infty}}^{ \frac 13 } \leq C \e ^{ -\frac 23 } 
\n \BA_{\e}  \n_{L^{2q}}^{ \frac 23 } .
$$

Next we prove  $(iii)$. As already 
mentioned, a uniform $L^p$ bound (for $2\leq p \leq 8/3$) on 
the kernels $\BK^1_{\e}$, $ \e ^2 \BK^\perp_{\e}$ and $ \e ^2 \BK^{1,j}_ {\e} $ 
is established in \cite{BGS1} by using  a Sobolev estimate. Unfortunately  this is no longer possible 
in dimension $N=3$. We thus rely on a suitable decomposition of $ \BA_{\e} $ in the Fourier space. 
Some terms are controlled  by using  the energy bounds in Lemma  \ref{BornEnergy},  the others by using \eqref{Henry}.

We consider a set  of parameters  $\a$, $\beta$, $\g \in (1,2)$ and $\nu > 5/2 $ 
with $ \a \geq \beta $ and $\a \geq \g $ (to  be fixed later).
For $ \e \in (0,1)$, let 
$$
\begin{array}{c}
E^I = \{ \xi \in \R^N \; \big| \; |\xi _{\perp} | < 1 \}, \quad
E^{II} = \{ \xi \in \R^N \; \big| \;  |\xi_\perp| > \e ^{-\a} \}, \quad 
E^{III} = \{ \xi \in \R^N \; \big| \; \e ^{-\beta} \leq |\xi_\perp|  \leq \e ^{-\a},\, |\xi_1| < 1 \}, 
\\
\\
E^{IV} =  \{ \xi \in \R^N \; \big| \; \e ^{-\g} \leq |\xi_\perp| 
\leq \e ^{-\a}, \, 1 \leq |\xi_1|^\nu \leq |\xi_\perp| \},  \qquad 
E^{V} =  \{ \xi \in \R^N \; \big| \; 1 \leq |\xi_\perp| \leq \e ^{-\a}, \, |\xi_1|^\nu > |\xi_\perp| \}, 
\\
\\
E^{VI} =  \{ \xi \in \R^N \; \big| \; 1 \leq |\xi_\perp| < \e ^{-\beta},  \, |\xi_1| < 1 \}, \qquad 
E^{VII} =  \{ \xi \in \R^N \; \big| \; 1 \leq |\xi_\perp| < \e ^{-\g},  \, 1 \leq |\xi_1|^\nu \leq |\xi_\perp| \}.
\end{array}
$$
It is easy to see  that the sets $E^I, \dots , E^{VII}$ are disjoint and cover $ \R^N$. 
For $J \in \{ I, \dots , VII \} $ we denote $ \BA_{\e} ^J = \Fou ^{-1} (\wh{\BA } _{\e} \mathbf{1}_{E^J})$, 
so that $ \BA_{\e} = \BA_{\e} ^I + \dots + \BA_{\e}^{VII} $, and we estimate each term separately.  

\medskip

For $ \BA_{\e} ^{I} $ we use
$$ \n \nabla_{z_\perp} \BA_{\e }^{I} \n_{L^2} 
= \n \xi_{\perp} \widehat{\BA}_{\e} {\bf 1}_{\{ |\xi_\perp| < 1 \} } \n_{L^2} 
\leq \n \widehat{\BA}_{\e }  {\bf 1}_{ \{ |\xi_\perp|\leq 1 \} } \n_{L^2} 
\leq \n \widehat{\BA}_{\E }  \n_{L^2} = \n \BA_{\E }  \n_{L^2} \leq C . $$
By Lemma \ref{BornEnergy}, $\BA_{\e} $ and $ \p_{z_1} \BA_{\e } $ 
are uniformly bounded in $L^2$, thus we have
$$ \n \BA_ {\e} ^{I} \n_{L^2} + \n \p_{z_1} \BA_{\e } ^{I} \n_{L^2} \leq C . $$
Hence $ \BA_{\e } ^{I} $ is uniformly bounded in $H^1$, and using  the Sobolev 
embedding  we deduce
\be
\label{Timide1}
\forall \; 2 \leq p \leq 6, \quad \quad \quad 
\n \BA_{\e } ^{I} \n_{L^p} \leq C . 
\ee
We will use the  Riesz-Thorin theorem to bound  $\BA_{\e }^{II}$: if $1<q = \frac{p}{p-1} <2$ is the conjugate exponent 
of $ p \in (2,  \ii) $, there holds 
$$ \n \BA_{\e } ^{II} \n_{L^p} \leq C \n \widehat{\BA}_{\e } ^{II} \n_{L^q} . $$
Thus it suffices to bound $  \n \widehat{\BA}_{\e } ^{II} \n_{L^q} $. 
Using the  H\"older inequality with exponents $\frac{2}{q} $ and  $\frac{2}{2-q} $, 
we have 
\begin{align*} 
\n \widehat{\BA}_{\e } ^{II} \n_{L^q}^q 
= & \ \int_{\R^3} 
\Big( (|\xi_1| + \e  |\xi_\perp|)  |\widehat{\BA}_{\e } | \Big)^q 
\times \frac{{\bf 1}_{\{ |\xi_\perp| > \e ^{-\a} \} } }{(|\xi_1| + \e  |\xi_\perp|)^q } \,  d \xi \\ 
\leq & \ \n (|\xi_1| + \e  |\xi_\perp|) \widehat{\BA}_{\e }  \n_{L^2}^q 
\left( \int_{\R^3} \frac{ {\bf 1}_{\{ |\xi_\perp| \geq \e ^{-\a} \} } }{(|\xi_1| + \e |\xi_\perp|)^{\frac{2q}{2-q} }} \,  d \xi \right)^{\frac{2-q}{q} } \\
\leq & \ C_q ( \n \p_{z_1} \BA_{\e }  \n_{L^2} 
+ \e  \n \nabla_{z_\perp} \BA_{\e }  \n_{L^2})^q 
\left( \int_{\e ^{-\a}}^{\ii} 
\frac{R \, dR}{ (\e  R)^{\frac{3q-2}{2-q}}} \right)^{\frac{2-q}{q} }. 
\end{align*}
(We have computed the integral in $\xi_1$ and we 
used cylindrical coordinates for the third line.)  Provided that $\frac{3q-2}{2-q}> 2$ (or, equivalently, 
$ q > 6/5 $), the last  integral in $R$ is
$$ C (q) \e ^{-\frac{3q-2}{2-q}} \times \e ^{\a \frac{5q-6}{2-q}} 
\leq C_q $$
as soon as $\a \geq \frac{3q-2}{5q-6} = \frac{2+p}{6-p}$, that is $ p \leq 6 - \frac{8}{\a + 1 }$.  
Notice that $2<6- \frac{8}{\a+1} <6 $ because $\a > 1 $. 
By Lemma   \ref{BornEnergy} we get
\be
\label{Timide2}
\forall \; 2 \leq p \leq 6 - \frac{8}{\a+1} , 
\quad \quad \quad \n \BA_{ \e } ^{II} \n_{L^p} \leq C(\a) . 
\ee

Using  similar arguments, we have
\begin{align*} 
\n \BA_{\e } ^{III} \n_{L^p}^q 
\leq & \ C \n \widehat{\BA}_{\e} ^{III} \n_{L^q}^q \\
= & \ C \int_{\R^3} 
\Big( \e  |\xi_\perp| \cdot |\wh{\BA}_{\e } | \Big)^q 
\times \frac{{\bf 1}_{ \{ \e ^{-\beta} \leq |\xi_\perp| \leq \e ^{-\a} , \, |\xi_1| <  1 \} } }
{(\e |\xi_\perp|)^q }  \, d \xi \\ 
\leq & \ C ( \e \n \nabla_{z_\perp} \BA_{\e }  \n_{L^2} )^q 
\left( \int_{\R^3} \frac{ {\bf 1}_{ \{ \e ^{-\beta} \leq |\xi_\perp| \leq \e ^{-\a} , \, |\xi_1| \leq 1 \} } }{(\e  |\xi_\perp|)^{\frac{2q}{2-q} }} 
\, d \xi \right)^{\frac{2-q}{q} }
\\ \leq & \ 
C_q \left( \e ^{-\frac{2q}{2-q} } \int_{\e ^{-\beta}}^{\e ^{-\a}} 
\frac{dR}{R^{\frac{4q-4}{2-q} + 1}} \right)^{\frac{2-q}{q} } \leq C_q 
\end{align*}
 if $\beta \frac{4q-4}{2-q} - \frac{2q}{2-q} \geq 0 $, 
that is $ 2 \beta \geq \frac{q}{(q-1)} = p $. Consequently,
\be
\label{Timide3}
\forall \; 2 \leq p \leq 2 \beta , 
\quad \quad \quad \n \BA_{\e }^{III} \n_{L^p} \leq C(\beta) . 
\ee
Similarly we get a bound for $ \BA_{\e } ^{IV}$:
\begin{align*} 
\n \BA_{\e } ^{IV} \n_{L^p}^q 
\leq & \ C ( \e  \n \nabla_{z_\perp} \BA_{\e }  \n_{L^2} )^q 
\left( \int_{\R^3} \frac{ {\bf 1}_{ \{ \e ^{-\g} \leq |\xi_\perp| \leq \e ^{-\a} ,\, 1 \leq |\xi_1|^\nu \leq |\xi_\perp| \} } }{(\e |\xi_\perp|)^{\frac{2q}{2-q} }} 
\, d  \xi \right)^{\frac{2-q}{q} }
\\ \leq & \ 
C_q \left( \e ^{-\frac{2q}{2-q} } \int_{\e ^{-\g}}^{\e ^{-\a}} 
\frac{R^{\frac{1}{\nu}} \, dR}{R^{\frac{4q-4}{2-q} + 1}} \right)^{\frac{2-q}{q} } 
\leq C_q
\end{align*}
provided that   $\g \frac{4q-4}{2-q} - \frac{2q}{2-q} - \frac{\g}{\nu} 
\geq 0 $, which is equivalent to  $ p \leq \frac{2\g (2 \nu + 1)}{2\nu+\g} $ (notice that
 $\frac{2\g (2 \nu + 1)}{2\nu+\g} > 2$  because $\g > 1$). Therefore,
\be
\label{Timide4}
\forall \;  2 \leq p \leq \frac{2\g (2 \nu + 1)}{2\nu+\g} , 
\quad \quad \quad \n \BA_{\e } ^{IV} \n_{L^p} \leq C(\nu) . 
\ee

We use the fact that $\| \p_{z_1} \BA_{\e} \|_{L^2}$ 
is bounded independently of $ \e$ (see part $(i)$) in order to estimate $\BA_{\e }^V$:
\begin{align*}
 \n \BA_{\e } ^{V} \n^q_{L^p} \leq & \ C \n \widehat{\BA}_{\e } ^{V} \n^q_{L^q} \\
= & \ C \int_{\R^3 } |\xi_1 \widehat{\BA}_{\e } |^q \times 
\frac{{\bf 1}_{\{ 1 \leq |\xi_\perp | \leq \e ^{-\a}, \, |\xi_\perp| < |\xi_1|^\nu \}}}{|\xi_1|^q} \, d \xi \\ 
\leq & \ C \n \p_{z_1} \BA_{\e }  \n_{L^2}^{q} \left( 
\int_{\R^3} \frac{{\bf 1}_{\{ 1 \leq |\xi_\perp | \leq \e ^{-\a}, \, |\xi_\perp| \leq |\xi_1|^\nu \}}}{|\xi_1|^{\frac{2q}{2-q}}} \, d \xi \right)^{\frac{2-q}{2}} \\ 
\leq & \ C \left( \int_1^{\e ^{-\a} } \frac{ R \, dR}{R^{(\frac{2q}{2-q}-1)/\nu}} 
\right)^{\frac{2-q}{2}} ,
\end{align*}
by using cylindrical coordinates in the fourth line.
We have 
$\frac{2q}{2-q} > 1$ for $q \in  [1, 2)$ and  the last integral is bounded 
independently of $\e $ as soon as $\frac{1}{\nu} \left(  \frac{2q}{2-q} - 1 \right)  > 2 $, 
that is $ p < \frac{4\nu +2}{2\nu-1} $. 
It is obvious that $ \frac{4\nu +2}{2\nu-1} > 2$ for $\nu > 1/2$. As a consequence, we get 
\be
\label{Timide5}
\forall \;  2 \leq p < \frac{4\nu +2}{2\nu-1} , 
\quad \quad \quad \n \BA_{\e } ^{V} \n_{L^p} \leq C(p) . 
\ee

We  use the convolution  equation \eqref{Henry} to estimate  $\BA_{\e } ^{VI}$ and $\BA_{\e } ^{VII}$.
Applying  the Fourier transform  to \eqref{Henry}  we obtain  the pointwise bound
\begin{align*}
|\widehat{\BA}_{\e }  (\xi) | = & \ \Big| 
\Big( \widehat{\BK}^1_{\e } + \e ^2 \widehat{\BK}^\perp_{\e } \Big) \widehat{G}_{\e }
+ 2 c(\e ) \widehat{\BK}^{\perp }_{\e }  \Fou ( \BA_{\e }  \p_{z_1} \vp_{\e }  ) 
- 2 c(\e ) \e ^2 \sum_{j=2}^N \widehat{\BK}^{1,j}_{\e } 
\Fou (\BA_{\e }  \p_{z_j} \vp_{\e }) \Big| \\ 
\leq & \ C \Big( |\widehat{\BK}^1_{\e }| + \e ^2 |\widehat{\BK}^\perp_{\e }| 
+ \e ^2 \sum_{j=2}^N |\widehat{\BK}^{1,j}_{\e }| \Big) 
\Big( \n \widehat{G}_{\e } \n_{L^\ii} + \n \Fou ( \BA_{\e } \p_{z_1} \vp_{\e } ) \n_{L^\ii} 
+ \sum_{j=2}^N \n \Fou ( \BA_{\e } \p_{z_j} \vp_{\e } ) \n_{L^\ii} \Big) .
\end{align*}
The estimates in Lemma \ref{BornEnergy} and the boundedness of 
$\Fou : L^1 \to L^\ii$ imply that the second factor is bounded independently 
of $\e$. Therefore 
\be
\label{Atchoum}
|\widehat{\BA}_{\e } (\xi) | \leq C \Big( |\widehat{\BK}^1_{\e }| 
+ \e ^2 |\widehat{\BK}^\perp_{\e }| 
+ \e ^2 \sum_{j=2}^N |\widehat{\BK}^{1,j}_{\e }| \Big) 
\leq C \frac{\xi_1^2 + \e ^2 |\xi_\perp|^2 
+ \e ^2 |\xi_1| \cdot |\xi_\perp|}{D_{\e }(\xi)} 
\leq C \frac{\xi_1^2 + \e ^2 |\xi_\perp|^2 }{D_{\e }(\xi)} 
\ee
because  $ 2 \e ^2 |\xi_1| \cdot |\xi_\perp| \leq \xi_1^2 + \e ^4 |\xi_\perp|^2$. 
If  $ \xi \in E^{VI}$ we have $ |\xi_1| \leq 1 $ and 
$ 1 \leq |\xi_\perp| \leq \e^{-\beta} \leq \e^{-2}$ (because
$\beta < 2$), hence there is   some constant $C$ depending only on 
$\cs$ such that
$$ C |\xi_\perp|^2 \geq D_{\e }(\xi) = \xi_1^4 + \xi_1^2 + \cs^2 |\xi_\perp|^2 
+ 2 \e ^2 \xi_1^2 |\xi_\perp|^2 + \e ^4 |\xi_\perp|^4 
\geq \frac{|\xi_\perp|^2}{C} . $$
Using the Riesz-Thorin theorem 
with exponents $2 < p < \ii$ and $q=p/(p-1) \in (1,2)$ as well as  \eqref{Atchoum} we find 
\begin{align*} 
\n \BA_{\e } ^{VI} \n_{L^p}^q 
\leq & \ 
C \n \widehat{\BA}_{\e } ^{VI} \n_{L^q}^q \\ 
\leq & \ C \int_{\R^3} 
{\bf 1}_{ \{ 1 \leq |\xi_\perp| \leq \e ^{-\beta}, \, |\xi_1| \leq 1 \}} 
\frac{( \xi_1^2 + \e ^2 |\xi_\perp|^2 )^q}{ |\xi_\perp|^{2q} } \ d \xi \\ 
\leq & \ C \int_{\R^3} 
{\bf 1}_{ \{ 1 \leq |\xi_\perp| \leq \e ^{-\beta}, \, |\xi_1| \leq 1 \} } 
\left( \frac{\xi_1^{2q}}{ |\xi_\perp|^{2q} } + \e ^{2q} \right) \ d \xi \\ 
\leq & \ C \int_{|\xi_\perp|\geq 1} \frac{d\xi_\perp}{|\xi_\perp|^{2q} } 
+ C \e ^{2q -2\beta} \leq C_q 
\end{align*}
provided that $q> 1$ 
and $ q \geq \beta $. 
We have $ q \geq \beta $ if and only if  $ p \leq \frac{\beta}{\beta-1}$. It is obvious that 
$\frac{\beta}{\beta-1}>2$ because $ 1 < \beta < 2$. Hence we obtain 
\be
\label{Timide6}
\forall \;  2 \leq p \leq \frac{\beta}{\beta-1} , 
\quad \quad \quad \n \BA_{\e }^{VI} \n_{L^p} \leq C(\beta) . 
\ee

In order to estimate  $\BA_{\e } ^{VII}$ we notice that for $ \xi \in E^{VII} $ we have
$ 1 \leq |\xi_\perp| \leq \e ^{-\g}$ and 
$ 1 \leq |\xi_1|^\nu \leq |\xi_\perp|$, thus
$|\xi_1|^2 \leq |\xi_\perp| \leq \e ^{-2}$ because 
$\nu \geq 5/2 > 2  $ and $ \g \leq 2$. Hence there exists 
$C > 0$ depending only on $\cs$ such that 
$$
 C |\xi_\perp|^2 \geq D_{\e} (\xi) = \xi_1^4 + \xi_1^2 + \cs^2 |\xi_\perp|^2 
+ 2 \e ^2 \xi_1^2 |\xi_\perp|^2 + \e ^4 |\xi_\perp|^4 
\geq \frac{|\xi_\perp|^2}{C} .
$$
Using \eqref{Atchoum} we get 
\begin{align*} 
\n \BA_{\e }^{VII} \n_{L^p}^q 
\leq & \ C \int_{\R^3} 
{\bf 1}_{ \{ 1 \leq |\xi_\perp| \leq \e ^{-\g}, \, 
1 \leq |\xi_1|^\nu \leq |\xi_\perp| \}} 
\left( \frac{\xi_1^{2q}}{ |\xi_\perp|^{2q} } + \e ^{2q} \right) \, d \xi \\ 
\leq & \ C \int_{|\xi_\perp|\geq 1} \frac{|\xi_\perp|^{\frac{2q+1}{\nu}}}{|\xi_\perp|^{2q} } \, d\xi_\perp
+ C \e ^{2q} \int_1^{\e ^{-\g}} R^{1+\frac{1}{\nu}} \, dR \leq C_q 
\end{align*}
provided that $2q - \frac{2q+1}{\nu} > 2 $ 
and  $2q - \g (2 + \frac{1}{\nu} ) \geq 0 $. 
These inequalities are equivalent to  $ p < \frac{2\nu +1}{3}$ and 
$ p \leq \frac{\g(2\nu +1)}{\g(2\nu +1) - 2 \nu}$, respectively.
 Since $\nu > 5/2$, we  have $ \frac{2\nu +1}{3} > 2 $ and $ \frac{4\nu}{2\nu+1} > 5/3 $ and $ \frac{\nu}{\nu-1} < 5/3 $. 
It is easy to see that $ \frac{\g(2\nu +1)}{\g(2\nu +1)-2\nu} > 2$ if and only if $ \g < \frac{4\nu}{2\nu+1} $,  
and that  $ \frac{\g(2\nu +1)}{\g(2\nu +1)-2\nu} > \frac{2\nu +1}{3}$ if and only if $ \g < \frac{\nu}{\nu-1} $.  
Hence
\be
\label{Timide7}
\left\{ \begin{array}{ll}
\ds{\forall \; 1 \leq \g \leq \frac{\nu}{\nu-1}, \quad 
\forall \;  2 \leq p < \frac{2\nu +1}{3}} , 
& \quad \quad \quad \n \BA_{\e }^{VII} \n_{L^p} \leq C(p,\nu) \\ \\
\ds{ \forall  \; \frac{\nu}{\nu-1} < \g \leq \frac53, \quad 
\forall \; 2 \leq p \leq \frac{\g(2\nu +1)}{\g(2\nu +1)-2\nu} ,} 
& \quad \quad \quad \n \BA_{\e }^{VII} \n_{L^p} \leq C(\g,\nu) .
\end{array}
\right.
\ee

We now choose the parameters $\a$, $\beta$, $\g$ and $\nu$. In view of 
\eqref{Timide3} and \eqref{Timide6}, we fix $\beta = 3/2$, so that 
$ 2 \beta = \beta/(\beta-1) = 3 $. We  set $\a = 5/3 > 3/2 = \beta $. 
Then  by \eqref{Timide1}, \eqref{Timide2}, \eqref{Timide3} and 
\eqref{Timide6}  it follows that 
$$ \forall \;  2 \leq p \leq 3, \quad \quad \quad 
\n \BA_{\e }^{I} \n_{L^p} + \n \BA_{\e }^{II} \n_{L^p} + \n \BA_{\e }^{III} \n_{L^p} 
+ \n \BA_{\e }^{VI} \n_{L^p} \leq C . $$
For the other terms,  
we notice that in the case $1 \leq \g \leq \frac{\nu}{\nu-1}$ we have 
$$ \frac{2\g (2 \nu + 1)}{2\nu+\g} \leq \frac{4 \nu + 2}{2\nu- 1} , $$
with equality if $\g = \frac{\nu}{\nu-1}$.
We also observe  that
$$ \frac{2\nu +1}{3} < \frac{4 \nu + 2}{2\nu-1} < \frac83 
\quad {\rm if} \quad \nu < \frac72, 
\quad \quad \quad {\rm respectively} \quad \quad \quad 
\frac83 < \frac{4 \nu + 2}{2\nu-1} < \frac{2\nu +1}{3} \quad 
{\rm if} \quad \nu > \frac72 . $$
Then we fix 
$\nu = 7 / 2$ and $\g = \frac{\nu}{\nu-1} = 7 / 5 < 5 / 3 $ 
and using \eqref{Timide4}, \eqref{Timide5} and \eqref{Timide7} 
we obtain
$$ \forall \; 2 \leq p < \frac83, \quad \quad \quad 
\n \BA_{\e }^{IV} \n_{L^p} + \n \BA_{\e }^{V} \n_{L^p} 
+ \n \BA_{\e }^{VII} \n_{L^p} \leq C . $$
This concludes the proof of $(iii)$.

\medskip

$(iv)$ We use 
the same inequalities as in the three-dimensional case with 
$ 1 < \nu < 3$ and   $\a$, $\beta$, $\g \in (1,2)$ satisfying 
$\beta \leq \a$ and $ \g \leq \a$.  We  get
$$ \begin{array}{llll}
\ds{ \forall \; 2 \leq p < \ii}, 
& \quad  \n \BA_{\e }^{I} \n_{L^p} \leq C_p; 
& \quad \quad \quad \ds{ \forall \; 2 \leq p \leq 4 \a - 2}, 
&  \quad \n \BA_{\e } ^{II} \n_{L^p} \leq C_p; \\ \\
\ds{ \forall \;  2 \leq p \leq \frac{2\beta}{2 - \beta},} 
&  \quad \n \BA_{\e }^{III} \n_{L^p} \leq C(\beta); 
& \quad \quad \quad \ds{ \forall \; 2 \leq p \leq 
\frac{2\g(\nu+1)}{\g + \nu(2-\g)},} 
&  \quad \n \BA_{\e }^{IV} \n_{L^p} \leq C(\beta) ; \\  \\
\ds{ \forall \;  2 \leq p < 2\frac{\nu+1}{\nu-1},} 
& \quad  \n \BA_{\e }^{V} \n_{L^p} \leq C_p ; 
& \quad \quad \quad \ds{ \forall \; 2 \leq p < \ii,} 
& \quad   \n \BA_{\e } ^{VI} \n_{L^p} \leq C_p 
\end{array} $$
and
$$ \forall \;  1 \leq \g \leq \frac{\nu}{\nu-1}, \quad 
\forall \; 2 \leq p < \frac{\nu +1}{3-\nu} , \quad \quad \quad 
\n \BA_{\e }^{VII} \n _{L^p} \leq C_p . $$
Then we choose 
$$ \beta = \frac43, \quad \quad \quad 
\a = \frac53 ,
\quad \quad \quad \nu = 3^- , \quad \quad \quad 
\g = \frac{\nu}{\nu-1} = \frac32^+ ,$$ 
so that  $ \a > \beta $ and $ \a > \g$. We infer that 
$$ \forall \; 2 \leq p < 4 , \quad \quad \quad 
\n \BA_{\e }  \n_{L^p} \leq C_p . $$
This completes  the proof in the case $N=2$.\carre

\subsection{Proof of Proposition \ref{Born}}

We first recall the Fourier multiplier properties of the kernels  $\BK_{\e }^1, $ $ \BK_{\e }^{\perp}$ and $\BK_{\e}^{1,j }$.
We skip the proof since it is  the same as in section 
5.2 in \cite{BGS1} and does not depend on the space dimension $N$.

\begin{lem}
\label{Multiply} 
Let $1 < q < \ii$. There exists $C_q>0$ (depending also on $\cs$) 
such that for any $ \e \in (0, 1 )$, any  $2 \leq j \leq N$ and $h \in L^q$ 
we have
\begin{align*}
& \n \BK^1_{\e } \star h \n_{L^q} \\ 
& + \n \p_{z_1} \BK^1_{\e } \star h \n_{L^q} 
+ \n \nabla_{z_\perp} \BK^1_{\e } \star h \n_{L^q} \\ 
& + \n \p_{z_1}^2 \BK^1_{\e } \star h \n_{L^q} 
+ \e \n \p_{z_1} \nabla_{z_\perp} \BK^1_{\e } \star h \n_{L^q} 
+ \e ^2 \n \nabla_{z_\perp}^2 \BK^1_{\e} \star h \n_{L^q} 
\leq C_q \n h \n_{L^q} ,
\end{align*}
\begin{align*}
& \n \BK^\perp_{\e } \star h \n_{L^q} \\
& + \e \n \p_{z_1} \BK^\perp_{\e } \star h \n_{L^q} 
+ \e ^2 \n \nabla_{z_\perp} \BK^\perp_{\e } \star h \n_{L^q} \\ 
& + \e ^2 \n \p_{z_1}^2 \BK^\perp_{\e } \star h \n_{L^q} 
+ \e ^3 \n \p_{z_1} \nabla_{z_\perp} \BK^\perp_{\e } \star h \n_{L^q} 
+ \e ^4 \n \nabla_{z_\perp}^2 \BK^\perp_{\e } \star h \n_{L^q} 
\leq C_q \n h \n_{L^q}
\end{align*}
and
\begin{align*}
& \n \BK^{1,j}_ {\e } \star h \n_{L^q} \\
& + \n \p_{z_1} \BK^{1,j}_{\e } \star h \n_{L^q} 
+ \e  \n \nabla_{z_\perp} \BK^{1,j}_{\e } \star h \n_{L^q} \\ 
& + \e \n \p_{z_1}^2 \BK^{1,j}_{\e } \star h \n_{L^q} 
+ \e ^2 \n \p_{z_1} \nabla_{z_\perp} \BK^{1,j}_{\e } \star h \n_{L^q} 
+ \e ^3 \n \nabla_{z_\perp}^2 \BK^{1,j}_{\e } \star h \n_{L^q} 
\leq C_q \n h \n_{L^q}.
\end{align*}
\end{lem}

The proof of 
\eqref{goodestimate} is then divided into 5 Steps. 

\bigskip

\noindent {\bf Step 1.} There is $ \e _ 1 > 0 $ and for any $1 < q < \ii $ there exists $C_q$ 
(depending also on $F$) such that for all $\e \in (0, \e _1)$,
\begin{align*}
\n \BA_{\e } \n_{L^q} + & \ \n \nabla_z \BA_{\e } \n_{L^q} 
+ \n \p^2_{z_1} \BA_{\e } \n_{L^q} 
+ \e \n \p_{z_1} \nabla_{z_\perp} \BA_{\e } \n_{L^q} 
+ \e ^2 \n \nabla_{z_\perp}^2 \BA_{\e } \n_{L^q} \nonumber \\
& \leq C_q \Big( \n \BA_{\e } \n^2_{L^{2q}} 
+ \e ^2 \Big[ \n \p_{z_1} \BA_{\e } \n_{L^{2q}} 
+ \e  \n \nabla_{z_\perp} \BA_{\e } \n_{L^{2q}} \Big]^2 \Big) .
\end{align*}

The proof is very similar to that of  Lemma 6.2 p. 268 in \cite{BGS1} and thus  is  only sketched. 
Indeed, if $ U = \rho e ^{ i \phi }$ is a finite energy solution to (TW$_c$) such that $ \frac{ r_0 }{2 } \leq \rho \leq 2  r _0$ 
then the first equation  in \eqref{phasemod} 
can be written as 
$$ 
2 r_0 ^2 \Delta \phi = c \frac{\p }{\p x_1} ( \rho ^2 - r_0 ^2) - 2 \mbox{div}\left( (\rho ^2 - r_0 ^2) \nabla \phi  \right)
$$
and this gives
$$
2 r_0 ^2 \frac{ \p \phi}{\p x_j } = c R_j R_1 ( \rho ^2 - r_0 ^2) - 2 \sum_{k = 1}^{N} R_j R_k \left( ( \rho ^2 - r_0 ^2 ) \frac{ \p \phi}{\p x_k } \right), 
$$
where $R_k$ is the Riesz transform (defined by $ R_k f = \Fou ^{-1} \left( \frac{ i \xi _k}{|\xi |} \wh{f} \right)$).
It is well-known that the Riesz transform maps continuously $ L^p(\R^N)$ into $ L^p(\R^N)$ for $ 1 < p < \ii$. 
From the above we infer that for any $ q \in (1, \ii) $ and any $ j \in \{ 1, \dots, N \}$ we have
$$
\Big\|  \frac{ \p \phi}{\p x_j } \Big\|_{L^q} \leq C(q) \| \rho ^2 - r_0 ^2 \|_{L^q} 
+ C(q) \sum_{k=1}^N \Big\| (\rho ^2 - r_0 ^2  ) \frac{ \p \phi}{\p x_j } \Big\|_{L^q}
\leq C(q) \| \rho ^2 - r_0 ^2 \|_{L^q}  + C(q) \| \rho ^2 - r_0 ^2  \|_{L^{\ii } } \| \nabla \phi \|_{L^q}
$$
and this implies 
$$
\| \nabla \phi \|_{L^q} \leq  C(q) \| \rho ^2 - r_0 ^2 \|_{L^q}  + C(q) \| \rho ^2 - r_0 ^2  \|_{L^{\ii } } \| \nabla \phi \|_{L^q}.
$$
If $\| \rho ^2 - r_0 ^2  \|_{L^{\ii } } $ is sufficiently small we get 
 $ \n \nabla \phi \n _{L^q} \leq \tilde{C} (q)  \| \rho ^2 - r_0 ^2 \|_{L^q}   \leq K(q) \n \rho - r_0 \n _{L^q}.$
By scaling, this estimate implies  that for $1 < q < \ii$,
\beq
\label{phasestimate1}
 \n \p_{z_1} \vp_{\e } \n_{L^q} 
+ \e \n \nabla_{z_\perp} \vp_{\e } \n_{L^q} \leq C_q \n \BA_{\e } \n_{L^q} . 
\eeq
Hence, by H\"older's inequality  and Lemma \ref{Grenouille} $(ii)$, 
\begin{align*}
\n G_{\e }  \n_{L^q} \leq & \ C_q \Big( \n \BA_{\e }  \n^2_{L^{2q}} 
+ \e ^2 \n \BA_{\e }  \n^3_{L^{3q}} 
+ \e ^2 \n \p_{z_1} \BA_{\e } \n^2_{L^{2q}} 
+ \e ^4 \n \nabla_{z_\perp} \BA_{\e }  \n^2_{L^{2q}} \Big) \\
\leq & \ C_q \Big( \n \BA_{\e }  \n^2_{L^{2q}}  
+ \e ^2 \Big[ \n \p_{z_1} \BA_{\e }  \n_{L^{2q}} 
+ \e  \n \nabla_{z_\perp} \BA_ {\e }  \n_{L^{2q}} \Big]^2 \Big)  .
\end{align*}
 We take  the derivatives up to order $2$ 
of \eqref{Henry} and then  the conclusion follows from Lemma \ref{Multiply}.

\bigskip

\noindent {\bf Step 2.} Let $N=3$. There is $ \e _2 > 0$ and  for any $1 < p < 3/2 $ there exists 
$C_p$ (also depending  on $F$) such that for any $\e \in (0, \e _2)$ there holds
\begin{align*}
\n \BA_{\e }  \n_{L^p} + \n \nabla \BA_{\e }  \n_{L^p} 
+ \n \p^2_{z_1} \BA_{\e }  \n_{L^p} 
+ \e  \n \p_{z_1} \nabla_{z_\perp} \BA_{\e }  \n_{L^p} 
+ \e ^2 \n \nabla_{z_\perp}^2 \BA_{\e } \n_{L^p} \leq C_p .
\end{align*}

If $1 \leq q \leq 3/2$, we have by Lemma \ref{Grenouille} $(i)$
$$ \e  \Big[ \n \p_{z_1} \BA_{\e } \n_{L^{2q}} 
+ \e  \n \nabla_{z_\perp} \BA_{\e } \n_{L^{2q}} \Big] \leq C . $$
Thus for $ 1 < q \leq 3/2$  we infer from Step 1 that
\begin{align}
\label{Garulfo}
\n \BA_{\e } \n_{L^q} + \n \nabla_z \BA_{\e } \n_{L^q} 
+ & \ \n \p^2_{z_1} \BA_{\e } \n_{L^q} 
+ \e \n \p_{z_1} \nabla_{z_\perp} \BA_{\e } \n_{L^q} 
+ \e ^2 \n \nabla_{z_\perp}^2 \BA_{\e } \n_{L^q} 
\leq C_q + C_q \n \BA_{\e } \n^2_{L^{2q}} . 
\end{align}
If $ 1 < p < 4/3 $,  we use  \eqref{Garulfo} 
combined with Lemma \ref{Grenouille} $(iii)$ with exponent
$ 2p \in [2,8/3)$ to get
\begin{align}
\label{Prof}
\n \BA_{\e } \n_{L^p} + \n \nabla_z \BA_{\e } \n_{L^p} 
+ & \, \n \p^2_{z_1} \BA_{\e } \n_{L^p} 
+ \e \n \p_{z_1} \nabla_{z_\perp} \BA_{\e } \n_{L^p} 
+ \e ^2 \n \nabla_{z_\perp}^2 \BA_{\e } \n_{L^p} \leq C_p .
\end{align}
This proves Step 2 for $1 < p < 4/3$. In dimension $N=3$, the Sobolev 
inequality does not enable us to improve the $L^q$ integrability of $\BA_{\e }$ to some 
$q>8/3$. We thus rely on 
the decomposition of $\BA_{\e }$ as 
$ \BA_{\e } = \BA_{\e }^{I} + \BA_{\e }^{II} + \BA_{\e }^{III} + \BA_{\e }^{IV} + \BA_{\e }^{V} 
+ \BA_{\e }^{VI} + \BA_{\e }^{VII} ,  $
exactly as  in Lemma \ref{Grenouille}. We choose  $\a = 5/3$, $\beta = 3/2$. 
By the estimates in the proof of Lemma \ref{Grenouille} $(iii)$  we  have  then
$$ \forall \;  2 \leq p \leq 3, \quad \quad \quad 
\n \BA_{\e }^{I} \n_{L^p} + \n \BA_{\e }^{II} \n_{L^p} + \n \BA_{\e }^{III} \n_{L^p} 
+ \n \BA_{\e }^{VI} \n_{L^p} \leq C .  $$
It remains  to bound $\BA_{\e }^{IV}$, $\BA_{\e }^{V}$ and $\BA_{\e }^{VII}$ in 
$L^{3^-}$. In view of \eqref{Timide5}, we choose $\nu = 5/2$, so 
that $\frac{4\nu+2}{2\nu-1} = 3$, and thus
$$ \forall \; 2 \leq p < 3, \quad \quad \quad \n \BA_{\e }^V \n_{L^p} \leq C_p . $$
We cancel out $\BA_{\e }^{IV}$ by taking $\g = 5/3 = \a$. 
Next we turn our attention to  the "bad term" $\BA_{\e }^{VII}$. By   \eqref{Prof} we get
$$ \forall \;  1 < p < \frac43 , \quad \quad \quad 
\n \nabla_{z_\perp} \BA_{\e } \n_{L^p} \leq C_p , $$
hence, by the Riesz-Thorin theorem,
$$ \forall \; 4 < r < \ii , \quad \quad \quad 
\n \xi_{\perp} \widehat{\BA}_{\e } \n_{L^r} = 
\n \Fou( \nabla_{z_\perp} \BA_{\e } ) \n_{L^r} \leq C_r . $$
Consequently, for $4 < r < \ii$, $ 2 < p < \ii $ and $q= p / (p-1) \in (1,2)$,   using once again the Riesz-Thorin theorem 
and  the H\"older inequality with exponents $\frac{r}{q} $ and $ \frac{r}{r-q} $ we get
\begin{align*}
\n \BA_{\e }^{VII} \n_{L^p}^q \leq & \ C \n \widehat{\BA}_{\e }^{VII} \n_{L^q}^q \\
= & \ C \int_{\R^3} ( |\xi_\perp| \cdot |\widehat{\BA}_ {\e } | )^q 
\times \frac{{\bf 1}_{ \{ 1 \leq |\xi_\perp| \leq \e ^{-\g},\, 
1 \leq |\xi_1|^\nu \leq |\xi_\perp| \} }}{|\xi_\perp|^q} \ d \xi \\ 
\leq & \ C \n \xi_\perp \widehat{\BA}_{\e } \n_{L^r}^{q} 
\Big( \int_{\R^3} \frac{{\bf 1}_{ \{ 1 \leq |\xi_\perp| \leq \e ^{-\g},\, 
1 \leq |\xi_1|^\nu \leq |\xi_\perp| \} }}{|\xi_\perp|^{\frac{rq}{r-q}}} \ d \xi \Big)^{\frac{r-q}{r}} \\ 
\leq & \ C_{r,q} \Big( \int_1^{\e ^{-\g}} \frac{R^{1+\frac{1}{\nu}}}{R^{\frac{rq}{r-q}}} \ d R \Big)^{\frac{r-q}{r}} 
\leq C_{r,q} 
\end{align*}
provided that $ \frac{rq}{r-q} > 2 + \frac{1}{\nu} = 12 / 5 $. 
Now  let $ 2 \leq p < 3 $ be fixed, so that  $ 3/2 < q \leq 2 $. Since 
$ 3/2 < q \leq 2 $ and $q \longmapsto  \frac{4q}{4-q} $ is increasing on $ (3/2 , 2 ] $, 
we have $\frac{4q}{4-q} > 12 / 5$. Furthermore,  we have 
$ \frac{rq}{r-q} \to \frac{4q}{4-q} > 12 / 5$ as $r \to 4$. Hence  we may choose 
 $r > 4$ such that $ \frac{rq}{r-q} > 2 + \frac{1}{\nu} = 12 / 5 $. 
As a consequence, we have
$$ \forall \; 2 \leq p < 3, \quad \quad \quad 
\n \BA_{\e }^{VII} \n_{L^p} \leq C_p . $$
Collecting the above estimates for $\BA_{\e }^{I} $, ... , $\BA_{\e }^{VII}$  
we deduce
$$ \forall \; 2 \leq p < 3, \quad \quad \quad \n \BA_{\e } \n_{L^p} \leq C_p . $$
Then we use  once again \eqref{Garulfo} with exponent $ p/2 \in (1,3/2) $ 
to infer that Step 2 holds for $ 1 < p < 3/2 $.\\

In order to be able to use Step 1 with some $q > 3/2$, we need to prove that 
$\BA_{\e }$, $\e_n \p_{z_1} \BA_{\e }$ and $ \e_n^2 \nabla_{z_\perp} \BA_{\e }$ 
are uniformly bounded in $L^p$ for some $p>3$. This is what we will prove next.

\bigskip

\noindent {\bf Step 3.} If $N=3$, the following bounds hold:
$$ \left\{ \begin{array}{ll}
\ds{ \forall \;  2 \leq p < 15/4 = 3.75 }, & 
\quad \quad \quad \n \BA_{\e } \n_{L^p} \leq C_p ;
\\ \\
\ds{ \forall \; 2 \leq p < 18/5 = 3.6 }, & 
\quad \quad \quad \e  \n \p_{z_1} \BA_{\e } \n_{L^p} \leq C_p ;
\\ \\
\ds{ \forall \; 2 \leq p < 18/5 = 3.6 }, & 
\quad \quad \quad \e ^2 \n \nabla_{z_\perp} \BA_{\e } \n_{L^p} \leq C_p .
\end{array}
\right. $$

Fix $ r \in (3, \ii)$,  $ p \in (2, \ii )$ and let $q= p/(p-1) \in (1,2)$ be the conjugate exponent of $p$.
 By  the Riesz-Thorin theorem and the H\"older inequality with exponents ${\frac{r}{q}} $ and $ {\frac{r}{r-q}}$  we have
\begin{align}
\label{sorciere}
\n \BA_{\e } \n_{L^p}^q \leq & \ C \n \widehat{\BA}_{\e } \n_{L^q}^q \nonumber \\
= & \ C \int_{\R^3} \Big[ (1 +|\xi_1|^2 + |\xi_\perp|) \cdot 
|\widehat{\BA}_{\e }| \Big]^q 
\times \frac{d \xi}{(1 + |\xi_1|^2 + |\xi_\perp|)^q} \nonumber \\ 
\leq & \ C \Big( \n \widehat{\BA}_{\e } \n_{L^r} +\n \xi_1^2 \widehat{\BA}_{\e } \n_{L^r} 
+ \n \xi_\perp \widehat{\BA}_{\e } \n_{L^r} \Big)^{q} 
\Big( \int_{\R^3} \frac{ d \xi}{( 1+ |\xi_1|^2 + |\xi_\perp|)^{\frac{rq}{r-q}}} \Big)^{\frac{r-q}{r}} .
\end{align}
We bound the first parenthesis using  again the Riesz-Thorin theorem: 
since $ r \in(3, \infty) $, its conjugate exponent $ r/(r-1) $ belongs to 
$ (1,3/2) $  and then Step 2  holds for the exponent  $ r$ instead of $p$, hence
\begin{align*}
\n \widehat{\BA}_{\e } \n_{L^r} + \n \xi_1^2 \widehat{\BA}_{\e } \n_{L^r} 
+ \n \xi_\perp \widehat{\BA}_{\e } \n_{L^r} 
= & \ \n \Fou( \BA_{\e }) \n_{L^r} + \n \Fou( \p_{z_1}^2 \BA_{\e }) \n_{L^r} 
+ \n \Fou( \nabla_{z_\perp} \BA_{\e } ) \n_{L^r} \\
\leq & \ C \Big( \n \BA_{\e } \n_{L^{\frac{r}{r-1}}} 
+ \n \p_{z_1}^2 \BA_{\e } \n_{L^{\frac{r}{r-1}}} 
+ \n \nabla_{z_\perp} \BA_{\e } \n_{L^{\frac{r}{r-1}}} \Big) \leq C_r .
\end{align*}
Next, we compute using cylindrical coordinates
\begin{align*}
\int_{\R^3} & \ \frac{ d \xi }{( 1+ |\xi_1|^2 + |\xi_\perp|)^{\frac{rq}{r-q}}} \\
& \leq 4 \pi \Big[ 
\int_0^1 \int_0^{+\ii} \frac{ R dR }{(1+R)^{\frac{rq}{r-q}}} \ d \xi_1 
+ \int_1^{+\ii} \int_0^{\xi_1^2} \frac{ R dR }{\xi_1^{\frac{2rq}{r-q}}} 
\ d \xi_1 + \int_1^{+\ii} \int_{\xi_1^2}^{+\ii} 
\frac{ R dR }{R^{\frac{rq}{r-q}}} \ d \xi_1 \Big] \\
& \leq 4 \pi \Big[ \int_0^{+\ii} \frac{ R dR }{(1+R)^{\frac{rq}{r-q}}} 
+ \frac12 \int_1^{+\ii} \frac{\xi_1^4}{\xi_1^{\frac{2rq}{r-q}}} \ d \xi_1 
+ \frac{1}{\frac{rq}{r-q}-2} 
\int_1^{+\ii} \frac{d \xi_1}{\xi_1^{2(\frac{rq}{r-q}-2)} } \Big] .
\end{align*}
The integrals in the last line are finite 
provided that $\frac{rq}{r-q} > 2$ (for the first integral), $\frac{2rq}{r-q} > 5$ 
(for the second integral) and $ 2(\frac{rq}{r-q}-2) > 1$ (for the third integral), hence their sum is finite if
 $\frac{rq}{r-q} > 5/2 $.  Note that $\frac{rq}{r-q} \to \frac{3q}{3-q} $ as $r \to 3$  and $ \frac{3q}{3-q} > 5/2 $ 
for $ q  \in (\frac{15}{11}, 3)$.   If $ 2 < p < 15 / 4 = 3.75$  we have 
$ 15 / 11 < q < 2 $ and we may choose  $ r > 3 $ 
(and $r$  close to $3$) such that $\frac{rq}{r-q} > 5/2 $. Then it follows from 
the two estimates above that
$$ \forall \; 2  \leq p < \frac{15}{4} , \quad \quad \quad 
\n \BA_{\e } \n_{L^p} \leq C_p . $$

Now we  turn our attention to the bound  on $\e  \p_{z_1} \BA_{\e }$. 
Let $r \in ( 1, \frac 32)$, $ q \in [2, \infty)$ and $ s \in (r,q)$. 
We use the estimates in Step 2 for $ \Big\| \frac{ \p ^2 \BA_{\e }}{\p z_i \p z_j} \Big\|_{L^r}$ and \eqref{bourrinSKF} with $N=3$
for $ \Big\| \frac{ \p ^2 \BA_{\e }}{\p z_i \p z_j} \Big\|_{L^q}$, then we interpolate to get
\beq
\label{81a}
 \Big\| \frac{ \p ^2 \BA_{\e }}{\p z_1 ^2} \Big\|_{L^s}
+ \e  \Big\| \frac{ \p ^2 \BA_{\e }}{\p z_1  \p z_j} \Big\|_{L^s}
+ \e ^2 \Big\| \nabla_{\perp} ^2 \BA _{\e} \Big\|_{L^s} \leq C_{r, q} \e^{\left(-4 + \frac{2N-1}{q} \right) \frac{ 1 - \frac rs }{1 - \frac rq}}.
\eeq
If $ s \in (r,3)$, from the Sobolev inequality and the above estimate we obtain 
\beq
\label{81}
\| \p_{z_1} \BA _{\e} \|_{L^{\frac{3s}{3-s}}} \leq C_s \| \p_{z_1} ^2 \BA_{\e } \|_{L^s} ^{\frac 13} \| \p _{z_1} \nabla_{\perp} \BA_{\e } \|_{L^s} ^{\frac 23} 
\leq C_{s,r,q} \e ^{-\frac 23} \e^{\left(-4 + \frac{5}{q} \right) \frac{ 1 - \frac rs }{1 - \frac rq}}.
\eeq
We have $ - \frac 23 + \left(-4 + \frac{5}{q} \right) \frac{ 1 - \frac rs }{1 - \frac rq} \to - \frac{14}{3}  + \frac{4r}{s} $
as $ q \to \ii$ uniformly with respect to $ r \in [1, \frac 32]$ and $ s \in [1,3]$. 
If $ 1 < s < \frac{18}{11} \approx 1.636$ 
we have $  - \frac{14}{3}  + \frac{4r}{s} \to  - \frac{14}{3} +  \frac 6s > -1$ as $ r \to \frac 32$. 
For any fixed $ s \in (1, \frac{18}{11})$ we may choose $q$ sufficiently large and  $ r \in (1,\frac 32)$ 
sufficiently close to $ \frac 32$ such that 
$- \frac 23 + \left(-4 + \frac{5}{q} \right) \frac{ 1 - \frac rs }{1 - \frac rq} > -1$. 
Since $\frac{3s}{3-s} \nearrow \frac{18}{5} $ as $ s \nearrow \frac{18}{11}$, from  \eqref{81} we get 
$$
\forall \; p \in \left(1, \frac {18}{5} \right), \qquad \quad \| \p_{z_1} \BA_{\e } \|_{L^p} \leq C_p \e^{-1}. 
$$

Let $r \in ( 1, \frac 32)$, $ q \in [3, \infty)$ and $ s \in (r,3)$. 
Using the Sobolev inequality and \eqref{81a} we have 
$$ \n  \nabla_{z_\perp} \BA_{\e } \n_{L^{\frac{3s}{3-s}}} \leq 
C_p \n  \p_{z_1} \nabla_{z_\perp} \BA_{\e } \n_{L^s}^{\frac13} 
\n  \nabla^2_{z_\perp} \BA_{\e } \n_{L^s}^{\frac23} 
\leq C_{s,r,q} \e ^{-\frac 53} \e^{\left(-4 + \frac{5}{q} \right) \frac{ 1 - \frac rs }{1 - \frac rq}}.
 $$
Proceeding  as above we infer that
$$ \forall \; 1 < p < 18 / 5, \quad \quad \quad 
\e ^2 \n   \nabla_{z_\perp} \BA_{\e } \n_{L^p} \leq C_p . $$

\medskip

\noindent {\bf Step 4.} Conclusion in the case $N=3$.

Fix $ 1 < p < 9 / 5 = 1.8 $. 
Since $ 2 < 2p < 18/ 5 < 15/4 $, we may use Step 1 (with  $p$ instead of $q$) 
and  Step 3 to deduce that 
\begin{align}
\label{princecharmant}
\n \BA_{\e } \n_{L^p} + & \ \n \nabla_z \BA_{\e } \n_{L^p} 
+ \n \p^2_{z_1} \BA_{\e } \n_{L^p} 
+ \e_n \n \p_{z_1} \nabla_{z_\perp} \BA_{\e } \n_{L^p} 
+ \e_n^2 \n \nabla_{z_\perp}^2 \BA_{\e } \n_{L^p} \nonumber \\
& \leq C_p \Big( \n \BA_{\e } \n^2_{L^{2p}} 
+ \Big[ \e_n \n \p_{z_1} \BA_{\e } \n_{L^{2p}} 
+ \e_n^2 \n \nabla_{z_\perp} \BA_{\e } \n_{L^{2p}} \Big]^2 \Big) \leq C_p .
\end{align}
Hence \eqref{goodestimate} holds for $p \in (1, \; 9/5)$. 
In particular, by the Sobolev imbeddding 
$ W^{1,p} \hookrightarrow L^{\frac{3p}{3-p}}$ with $ 1 < p < 9/5 $   we have
$$ \forall \, 1 < q < 9/2 = 4.5 , \quad \quad \quad \n \BA_{\e } \n_{L^q} \leq C_q . $$
On the other hand, for any $1 < p < 9/5 $,
$$  \e \n  \p_{z_1} \BA_{\e } \n_{W^{1,p}} 
= \e \n  \p_{z_1} \BA_{\e } \n_{L^{p}} 
+  \e \n  \p^2_{z_1} \BA_{\e } \n_{L^{p}} 
+  \e \n  \nabla_{z_\perp} \p_{z_1} \BA_{\e } \n_{L^{p}} \leq C_p 
\quad \quad {\rm and} \quad \quad  \e ^2 \n  \nabla_{z_\perp} \BA_{\e } \n_{W^{1,p}} 
\leq C_p , $$
hence by the Sobolev embdding,
$$ \forall \; 1 < q < 9/2 = 4.5 , \quad \quad \quad 
\e \n  \p_{z_1} \BA_{\e } \n_{L^q} +  \e  \n  \nabla_{z_\perp} \BA_{\e } \n_{L^q} 
\leq C_q . $$
Thus we may apply Step 1  again  to infer 
that \eqref{princecharmant} holds now for $1 < p < 9/4 = 2.25 $. By the
Sobolev embedding $ W^{1,p} \hookrightarrow L^{\frac{3p}{3-p}}$, we 
deduce as before that
$$ \forall \; 1 < q < 9, \quad \quad \quad \n \BA_{\e } \n_{L^q} 
+ \e \n  \p_{z_1} \BA_{\e } \n_{L^q} 
+  \e ^2 \n  \nabla_{z_\perp} \BA_{\e } \n_{L^q} \leq C_q . $$
Applying  Step 1, we discover that \eqref{princecharmant} 
holds  for any $1 < p < 9/2 $. Since $ 9/2 > 3 $, the Sobolev 
embedding yields 
$$ \forall \; 1 < q \leq \ii, \quad \quad \quad 
\n \BA_{\e } \n_{L^p} + \e \n  \p_{z_1} \BA_{\e } \n_{L^p} 
+ \e ^2  \n  \nabla_{z_\perp} \BA_{\e } \n_{L^p} \leq C_p , $$
and  the  conclusion  follows  using  again  Step  1.

\bigskip

\noindent {\bf Step 5.} Conclusion in the case $N=2$. 
The proof of \eqref{goodestimate} in the two-dimensional case is much easier: for any $1 < p < \frac 32$, 
we have by Step 1 and Lemma \ref{Grenouille} $(i)$ and $(iv)$ 
$$ \n \BA_{\e } \n_{L^p} + \n \nabla_z \BA_{\e } \n_{L^p} 
+ \n \p^2_{z_1} \BA_{\e } \n_{L^p} 
+ \e \n \p_{z_1} \nabla_{z_\perp} \BA_{\e } \n_{L^p} 
+ \e ^2 \n \nabla_{z_\perp}^2 \BA_{\e } \n_{L^p} \leq C_p . $$
Thus, by the Sobolev embedding $W^{1,p} (\R^2 ) \hookrightarrow L^{\frac{2p}{2-p}} (\R^2 ) $,
\beq
\label{85}
 \forall \; 1 < q < 6, \quad \quad \quad \n \BA_{\e } \n_{L^q} \leq C_q 
\quad \quad \quad {\rm and} \quad \quad \quad 
\e_n \Big[ \n \p_{z_1} \BA_{\e } \n_{L^{q}} 
+ \e_n \n \nabla_{z_\perp} \BA_{\e } \n_{L^{q}} \Big] \leq C_q . 
\eeq
Applying Step 1 once again, we infer that \eqref{princecharmant} holds for any $ p \in (1, 3)$.  
Since $ 3 >2$, the Sobolev embedding implies that \eqref{85} holds for any $ q \in (1, \infty]$. 
Repeating the argument we get the desired conclusion. \\


Since 
$A_{\e } = \e^{-2} ( \sqrt{1+\e ^2 \BA_{\e }} - 1 )$, uniform 
bounds bounds on $ A_{\e }$ and its derivatives up to order 2 follow immediately from \eqref{goodestimate}.

\bigskip

It remains to prove \eqref{goodestimate2}.  The uniform bounds on  
$\p_{z_1} \vp_{\e }$ and $ \e  \nabla_{z_\perp} \vp_ {\e }$  follow from \eqref{phasestimate1} and \eqref{goodestimate}. 

Let  $ U = \rho \ex ^{i \phi } $   be a finite energy solution to (TW$_c$), from the first equation in \eqref{phasemod} 
we have 
$$
2 \rho ^2 \Delta \phi = c \frac{ \p }{ \p x_1 } ( \rho ^2 - r_0 ^2) - 2 \nabla ( \rho ^2) \cdot \nabla \phi. 
$$
If $ \rho \geq \frac{ r_0}{2} $ and  $ c \in (0, \cs)$, using the  properties of the Riesz transform we get 
for any $ j, k \in \{ 1, \dots, N \}$ and any $ q \in (1, \infty)$
$$
\Big\| \frac{ \p ^2 \phi}{\p x_j \p x_k } \Big\|_{L^q} = \| R_j R_k ( \Delta \phi ) \|_{L^q} 
\leq C \|  \Delta \phi  \|_{L^q}  
\leq C \Big\|  \frac{ \p }{ \p x_1 } ( \rho ^2 - r_0 ^2) \Big\| _{L^q} + C \| \nabla ( \rho ^2) \cdot \nabla \phi \|_{L^q}. 
$$
In the case $ U = U_{\e}$, $ \rho (x)  = r_0 \sqrt{ 1 + \e ^2 \BA_{\e} (z) }$, $\phi (x) = \e \vp _{\e} (z)$, 
using \eqref{goodestimate} and \eqref{phasestimate1} we get 
$$
\Big\| \frac{ \p ^2 \phi}{\p x_j \p x_k } \Big\|_{L^q}  
\leq \e^{ 3 - \frac{2N -1}{q} } \Big\| \frac{ \p \BA_{\e}}{\p z _1} \Big\| _{L^q} 
+ C \e^{ 5 - \frac{2N -1}{q}} \Big\| \frac{ \p \BA_{\e}}{\p z _1} \cdot \frac{ \p \vp_{\e }}{\p z_1} \Big\| _{L^q}  
+ C \e^{7 - \frac{2N -1}{q}}  \sum_{ j =2}^N \Big\| \frac{ \p \BA_{\e}}{\p z _j} \cdot \frac{ \p \vp_{\e }}{\p z_j} \Big\| _{L^q} 
\leq C_q \e^{ 3 - \frac{2N -1}{q} }.
$$
By scaling we find  for $ j,k \in \{2, \dots, N \}$, 
\beq
\label{phasestimate3}
\Big\| \frac{ \p ^2 \vp _{\e}}{\p z_1 ^2} \Big\|_{L^q} 
 + \e \Big\| \frac{ \p ^2 \vp _{\e}}{\p z_1  \p z_j } \Big\|_{L^q}  
+  \e ^2 \Big\| \frac{ \p ^2 \vp _{\e}}{\p z_j  \p z_k } \Big\|_{L^q}  \leq C_q. 
\eeq
By assumption (A4) there is $ \de > 0$ such that $F$ is $C^2$ on $ (\, ( r_0  -2 \de )^2, ( r_0 + 2 \de )^2) $.
Let $ U = \rho \ex ^{ i \phi }$  be a  solution to (TW$_c$) such that $ r_0 - \de \leq \rho \leq r_0 + \de$. 
Differentiating (TW$_c$) and using standard elliptic regularity theory it is not hard to see that $ U \in W_{loc}^{4, p} ( \R^N)$ 
and $ \nabla U \in W^{3, p}(\R^N)$ for any $ p \in (1, \infty) $ 
(see the proof Proposition 2.2 (ii) p. 1079 in \cite{M2}). We infer that $ \nabla \rho, \,  \nabla \phi \in  W^{3, p} ( \R^N)$ 
for $ p \in (1, \infty) $. 
Differentiating the first equation in \eqref{phasemod}  with respect to $ x_1 $ we find 
\beq
\label{deriveq}
c \frac{ \p ^2}{\p x_1 ^2} \left( \rho ^2 - r_0 ^2\right) = 2 \nabla \left( \frac{ \p ( \rho ^2)}{\p x_1 } \right) \cdot \nabla \phi 
+ 2 \nabla ( \rho ^2) \cdot \nabla \left( \frac{ \p \phi}{\p x_1} \right) + 2 \frac{ \p ( \rho ^2)}{\p x_1 } \Delta \phi 
+ 2 \rho ^2 \Delta \left( \frac{ \p \phi}{\p x_1} \right). 
\eeq
If $ U = U_{\e}$, $ \rho (x) = r_0 \sqrt{ 1 + \e ^2 \BA_{\e } (z) } $ and 
$ \phi (x) = \e \vp_{\e} (x)$, we perform a scaling and then we use 
\eqref{goodestimate},  \eqref{phasestimate1} and \eqref{phasestimate3} to get, for $ 1 < q < \ii$ and all $ \e $ sufficiently small, 
$$
\begin{array}{c}
\ds \Big\| \frac{ \p ^2}{\p x_1 ^2} \left( \rho ^2 - r_0 ^2 \right) \Big\|_{L^q} 
= \e^{ 4 + \frac{1 - 2N}{q}} \Big\| \frac{ \p ^2 \BA_{\e}}{\p z_1 ^2 } \Big\|_{L^q} \leq C_q  \e^{ 4 + \frac{1 - 2N}{q}}, 
\\
\\
\ds \Big\| \frac{ \p ^2 ( \rho ^2)}{\p x_1 ^2 } \cdot \frac{\p \phi}{\p x_1} \Big\|_{L^q} 
\leq \Big\|   \frac{ \p ^2 ( \rho ^2)}{\p x_1 ^2 }   \Big\|_{L^{2q}}   \Big\|   \frac{\p \phi}{\p x_1}  \Big\|_{L^{2q}}   
= \e^{ 6 + \frac{1 - 2N}{q}} \Big\| \frac{ \p^2 \BA_{\e}}{\p z_1 ^2}   \Big\|_{L^{2q}}  \Big\| \frac{\p \vp _{\e}}{\p z_1} \Big\|_{L^{2q}} 
\leq C_q  \e^{ 6 + \frac{1 - 2N}{q}}, 
\\
\\
\ds \Big\| \frac{ \p ^2 ( \rho ^2)}{\p x_1  \p x_{k}  } \cdot  \frac{\p \phi}{\p x_k} \Big\|_{L^q} 
\leq \Big\|   \frac{ \p ^2 ( \rho ^2)}{\p x_1  \p x_k }   \Big\|_{L^{2q}}   \Big\|   \frac{\p \phi}{\p x_k}  \Big\|_{L^{2q}}   
= \e^{ 8 + \frac{1 - 2N}{q}} \Big\| \frac{ \p^2 \BA_{\e}}{\p z_1  \p z_k}   \Big\|_{L^{2q}}  \Big\| \frac{\p \vp _{\e}}{\p z_k} \Big\|_{L^{2q}} 
\leq C_q  \e^{ 6 + \frac{1 - 2N}{q}}, 
\\
\\
\ds \Big\| \frac{ \p  ( \rho ^2)}{\p x_1  } \cdot  \frac{\p ^2 \phi}{\p x_1 ^2 } \Big\|_{L^q} 
\leq \Big\|   \frac{ \p  ( \rho ^2)}{\p x_1  }   \Big\|_{L^{2q}}   \Big\|   \frac{\p ^2 \phi}{\p x_1 ^2 }  \Big\|_{L^{2q}}   
= \e^{ 6 + \frac{1 - 2N}{q}} \Big\| \frac{ \p \BA_{\e}}{\p z_1 }   \Big\|_{L^{2q}}  \Big\| \frac{\p ^2 \vp _{\e}}{\p z_1 ^2 } \Big\|_{L^{2q}} 
\leq C_q  \e^{ 6 + \frac{1 - 2N}{q}}, 
\\
\\
\ds \Big\| \frac{ \p  ( \rho ^2)}{\p x_k  } \cdot  \frac{\p ^2 \phi}{\p x_1  \p x_k } \Big\|_{L^q} 
\leq \Big\|   \frac{ \p  ( \rho ^2)}{\p x_k  }   \Big\|_{L^{2q}}   \Big\|   \frac{\p ^2 \phi}{\p x_1 \p x_k  }  \Big\|_{L^{2q}}   
= \e^{ 8 + \frac{1 - 2N}{q}} \Big\| \frac{ \p \BA_{\e}}{\p z_k }   \Big\|_{L^{2q}}  \Big\| \frac{\p ^2 \vp _{\e}}{\p z_1 \p z_k } \Big\|_{L^{2q}} 
\leq C_q  \e^{ 7 + \frac{1 - 2N}{q}}, 
\\
\\
\ds \Big\| \frac{ \p  ( \rho ^2)}{\p x_1  } \Big\|_{L^q}  = \e^{ 3 + \frac{1 - 2N}{q}}  \Big\|  \frac{ \p  \BA_{\e}}{\p z_1}   \Big\|_{L^q} 
\leq  C_q  \e^{ 3 + \frac{1 - 2N}{q}}, 
\\
\\
\ds \Big\| \frac{ \p ^2  \phi }{\p x_1 ^2 } \Big\|_{L^q}  = \e^{ 3 + \frac{1 - 2N}{q}}  \Big\|  \frac{ \p ^2  \vp_{\e}}{\p z_1 ^2 }   \Big\|_{L^q} 
\leq  C_q  \e^{ 3 + \frac{1 - 2N}{q}}
\qquad \mbox{ and } \qquad
\Big\| \frac{ \p ^2  \phi  }{\p x_k ^2 } \Big\|_{L^q}  = \e^{ 5 + \frac{1 - 2N}{q}}  \Big\|  \frac{ \p ^2  \vp_{\e}}{\p z_k ^2 }   \Big\|_{L^q} 
\leq  C_q  \e^{ 3 + \frac{1 - 2N}{q}} .
\end{array}
$$
Hence $ \| \Delta \phi \|_{L^q} \leq C_q  \e^{ 3 + \frac{1 - 2N}{q}}  $ and   then 
 $\ds \Big\| \frac{ \p (\rho ^2)}{\p x_1} \cdot \Delta \phi \Big\|_{L^q} \leq C_q  \e^{ 6 + \frac{1 - 2N}{q}}  $.
 From \eqref{deriveq}  and the above estimates we infer that 
$\ds \Big\| \Delta \left(\frac{ \p \phi }{\p x_1} \right) \Big\|_{L^q} \leq C_q  \e^{ 4 + \frac{1 - 2N}{q}}  $.
As before, this implies $\ds  \Big\| \frac{ \p ^3 \phi}{\p x_1 \p x_i \p x_j } \Big\|_{L^q} \leq C_q  \e^{ 4 + \frac{1 - 2N}{q}}  $
for any $ i, j \in \{ 1, \dots, N \}$. 
By scaling we find
$$
\Big\| \frac{ \p ^3 \vp _{\e}}{\p z_1 ^3} \Big\|_{L^q} 
 + \e \Big\| \nabla_{z_{\perp}}  \frac{ \p ^2 \vp _{\e}}{\p z_1  ^2} \Big\|_{L^q}  
+  \e ^2 \Big\| \nabla_{z_{\perp}} ^2  \frac{ \p  \vp _{\e}}{\p z_1  } \Big\|_{L^q}  \leq C_q. 
$$
Then \eqref{goodestimate2} follows from the last estimate, \eqref{phasestimate1} and \eqref{phasestimate3}.
\carre

\subsection{Proof of Proposition \ref{convergence}} 

Let $(U_n, \e _n)_{n \geq 1}$ be a sequence as in Proposition \ref{convergence}.
We denote $ c_n = \sqrt{\cs ^2 - \e _n ^2}$.
By Corollary \ref{sanszero} we have $ \| \, |U_n| - r_0 \| _{L^{\infty}(\R^3 )} \to 0 $ as $ n \to \infty$, 
hence $|U_n| \geq \frac{ r_0}{2}$ in $ \R^3$ for all sufficiently large $n$, say $ n \geq n_0$. 
For $ n \geq n_0$  we have a lifting    as in Theorem \ref{res1} or in \eqref{ansatz}, that is 
$$
 U_n (x) = \rho _n (x) \ex^{i\phi _n(x)} 
=r_0 \left( 1 + \e _n ^2A_n(z) \right) \ex^{i \e _n \vp _n(z) }
= r_0 \sqrt{1+\e_n ^2 \BA_{n  }(z) }\ \ex^{i\e _n \vp_{n } (z)}, 
$$
$
\mbox{where }   z_1 = \e _n x_1 , \; z_\perp = \e_n ^2 x_\perp .
$
Let $\BW_n = \p_{z_1} \vp_n / \cs $. 
Our aim is  to show that $(\BW_n )_{n \geq n_0}$ is  a minimizing sequence for $\bS_*$ in the sense of Theorem \ref{gs}.
To that purpose we  expand the functional $E_{c_n} (U_n)$ in terms of 
the (KP-I) action of $\BW_n = \p_{z_1} \vp_n / \cs $. 
Recall that by \eqref{develo} we have
\begin{align*}
E_{c_n} (u_n) = & \ \e_n r_0^2 \int_{\R^3} \frac{1}{\e_n^2} 
\Big( \p_{z_1} \vp_n - c_n A_n \Big)^2 
+ (\p_{z_1} \vp_n)^2 ( 2 A_n + \e_n^2 A_n^2 ) 
+ |\nabla_{z_\perp} \vp_n |^2 ( 1 + \e_n^2 A_n )^2 
\nonumber \\
& \hspace{2cm} + (\p_{z_1} A_n)^2 
+ \e_n^2 |\nabla_{z_\perp} A_n|^2 + A_n^2 
+ \cs^2 \Big( \frac{\G}{3} - 1 \Big) A_n^3 
+ \frac{\cs^2}{\e_n^6} V_4( \e_n^2 A_n) \nonumber \\
& \hspace{2cm} - c_n A_n^2 \p_{z_1} \vp_n \ dz .
\end{align*}
By Proposition \ref{Born}, $(A_n)_{n \geq n_0 }$ is bounded in $W^{1,p}(\R^N)$ for all $ p \in (1, \infty)$, 
hence it is bounded in $L^{\ii }(\R^3)$. 
Since $ F( r_0 ^2 ( 1 + \e ^2 A_{\e } )) = F( r_0 ^2) - \cs ^2 \e ^2 A_{\e } + \BO ( \e ^4 A_{\e } ^2 ) 
= - c^2( \e) \e ^2 A_{\e } - \e^4 A_{\e } + \BO ( \e ^4 A_{\e }  ) $, 
from the second equation in \eqref{MadTW}, Lemma \ref{BornEnergy} and Proposition \ref{Born} we get
\beq
\label{approx1}
\| \p _{z_1} \vp _{n } - c _n A_{n } \|_{L^2} = \BO(\e_n ^2). 
\eeq
In particular, we have $ \ds \int_{\R^3} \frac{1}{\e_n^2} \Big( \p_{z_1} \vp_n - c_n A_n \Big)^2 \, dz = \BO(\e_n ^2) $ as $ n \to \ii$.

By Proposition \ref{Born}, $ \p_{z_1} \vp _n \in W^{2, p} (\R^N)$ for $ p \in (1, \ii)$. Integrating by parts  we have
$$
\int_{\R^N} ( \p _{z_1} A _n ) ^2 - \frac{( \p ^2_{z_1} \vp _n ) ^2}{ c_n^2 }\,  dz 
= -  \int_{\R^N} \left(   A _n   - \frac{ \p _{z_1} \vp _n}{ c_n }\right)   \left( \p _{z_1} ^2 A _n   + \frac{ \p ^3_{z_1} \vp _n}{ c _n }\right)\,  dz 
$$
From the above identity, the Cauchy-Schwarz inequality, \eqref{approx1} and Proposition \ref{Born} we get
$$
\Big\vert \int_{\R^N} ( \p _{z_1} A _n ) ^2 -  \frac{ ( \p_{z_1} ^2\vp _n )^2}{\cs ^2}  \, dz   \Big\vert
\leq \left( \frac{1}{ c_n ^2 } - \frac{1}{\cs ^2} \right) \int_{\R^N} (\p_{z_1} ^2 \vp _n )^2 \, dz + 
\Big\| A_n    - \frac{ \p _{z_1} \vp _n}{ c_n } \Big\|_{L^2} 
\Big\|  \p _{z_1} ^2 A _n   + \frac{ \p ^3_{z_1} \vp _n}{ c _n } \Big\|_{L^2}
= \BO( \e _n ^2). 
$$
Similarly, using \eqref{approx1}, H\"older's inequality and Proposition \ref{Born} we find
$$ 
\begin{array}{l}
\ds \Big| \int_{\R^3} A_n^2 -  \frac{(\p_{z_1} \vp_n)^2}{\cs^2} \, dz \Big| 
+ \Big| \int_{\R^3} A_n^3 - \frac{(\p_{z_1} \vp_n)^3}{\cs^3} \, dz \Big| 
\\
\\
\ds + \Big| \int_{\R^3} A_n^2 \p_{z_1} \vp_n - \frac{(\p_{z_1} \vp_n)^3}{\cs^2} \, dz \Big| 
+ \Big| \int_{\R^3} A_n ( \p_{z_1} \vp _n )^2 - \frac{(\p_{z_1} \vp_n)^3}{\cs}  \, dz \Big| 
= \BO(\e_n ^2) . 
\end{array}
$$
Since $(A_n)_{n \geq n_0}$  is bounded in $L^{\infty}(\R^3)$,  using Lemma  \ref{BornEnergy} we find
$$
\int_{\R^3 } |\nabla_{z_\perp} \vp_n |^2 ( 1 + \e_n^2 A_n )^2  \, dz 
= \int_{\R^3 } |\nabla_{z_\perp} \vp_n |^2 \, dz + \BO( \e _n ^2) 
= \cs ^2 \int_{\R^3 } |\nabla_{z_\perp} \p_{z_1}^{-1} \BW _n |^2 \, dz + \BO( \e _n ^2) .
$$
Recall that  $ V_4( \al ) = \BO(\al ^4) $ as $ \al \to 0$, hence Proposition \ref{Born} implies that
$$
 \int_{\R^3 } \e _n ^2 A_n ^2 (\p_{z_1} \vp _n )^2 + \e_n^2 |\nabla_{z_\perp} A_n|^2 +  \frac{\cs^2}{\e_n^6} V_4( \e_n^2 A_n)  \, dz = \BO (\e _n ^2).
$$
Inserting  the above estimates into \eqref{develo} 
we obtain
\beq
\label{ecun}
 \frac{E_{c(\e_n)} (U_n)}{\cs^2 r_0^2 \e_n} = 
\int_{\R^3} \Big| \nabla_{z_\perp} \p_{z_1}^{-1} \BW_n \Big|^2 
+ \frac{1}{\cs^2}\, (\p_{z_1} \BW_n )^2 
+ \frac{\G}{3}\, \BW_n^3 + \frac{1}{\cs^2}\, \BW_n^2 \, dz + \BO(\e_n ^2) 
= \bS ( \BW_n ) + \BO(\e_n ^2)  . 
\eeq
From the  above estimate and the upper bound on $E_{c_n} (U_n) = T_{c_n}$ given by Proposition \ref{asympto} $(ii)$ 
 we infer that
$$ 
\bS ( \BW_n ) =   \frac{E_{c(\e_n)} (U_n)}{\cs^2 r_0^2 \e_n} + \BO(\e_n ^2) 
= \frac{ T_{c_n}}{\cs^2 r_0^2 \e_n} + \BO(\e_n ^2)  
  \leq \bS_{\rm min} + \BO(\e_n ^2)  = \bS_{*} + \BO(\e_n ^2)   . 
$$   
Similarly we have
$$
\int_{\R^3}  |\nabla _{x_{\perp}} U_n |^2 \, dx 
= r_0 ^2 \e _n \int_{\R^3} ( 1 + \e _n ^2 A_n) ^2 |\nabla _{x_{\perp}} \vp _n |^2 + \e _n ^2 |\nabla _{x_{\perp}} A_n |^2 \, dz
=  r_0 ^2 \cs ^2 \e _n \int_{\R^3}  \Big| \nabla_{z_\perp} \p_{z_1}^{-1} \BW_n \Big|^2 \, dz + \BO( \e_n ^3).
$$
Since $U_n$ satisfies the Pohozaev identity       
$\ds E_{c_n}(U_n) = \int_{\R^3} | \nabla_{z_\perp} U_n |^2 \ dz $, 
comparing the above equation to the expression of $E_{c_n}(U_n) $ in \eqref{ecun} we find 
$$  
\int_{\R^3} \frac{1}{\cs^2}\, (\p_{z_1} \BW_n )^2    
+ \frac{\G}{3}\, \BW_n^3 + \frac{1}{\cs^2}\, \BW_n^2 \ dz = \BO(\e _n^2) . 
$$

In order to apply Theorem \ref{gs}, we have to check that there is $ m_1 >0 $ such that
for all $n$ sufficiently large there holds
$$ \int_{\R^3} \BW_n^2 + (\p_{z_1} \BW_n)^2 \ dz \geq m_1 . $$
By  Lemma \ref{minoinf}, there are $ k>0 $ depending only on $F$ and $ n_1 \geq n_0$ such that
$$ \forall n \geq n_1, \quad \n A_n \n_{L^\ii} \geq k . $$
Since $A_n$ tends to $0$ at infinity, after a translation we may assume  that 
$$ |A_n(0)| = \n A_n \n_{L^\ii} \geq k . $$
By Proposition \ref{Born} we know that 
for all $ p \in (1, \infty)$ there is $ C_p > 0 $ such that $ \| A_{n } \|_{W^{1, p } } \leq C_p$  
for any $ n \geq n_0$.
Then Morrey's inequality (see e.g. Theorem IX.12 p. 166 in \cite{brezis}) implies that for any $ \al \in (0,1)$ there is $ C_{\al } > 0$ such that  
for all $ n \geq n_0$  and all $x, y \in \R^3$ we have $|A_{n }(x) - A_{n }(y) | \leq C_{\al } |x -y |^{\al}.$
We infer that   $ |A_n| \geq k/2 $ in $B_r(0)$ for some $r>0$ 
independent of $n$, hence there is $m_1 > 0$ such that
$$
 \n A_n \n_{L^2} \geq \n A_n \n_{L^2(B_r(0))} \geq 2 m_1.
$$
From \eqref{approx1} it follows that $\| \BW _n - A_n \|_{L^2} \to 0 $ as $ n \to \infty$, hence 
$$ \n \BW_n \n_{L^2} \geq \n \BW_n \n_{L^2(B_r(0))} 
\geq m_1 \qquad \mbox{ for all } n \mbox{ sufficiently large.}  $$

Then Theorem \ref{gs} implies that there exist $\BW \in \mathscr{Y}(\R^3)$,
a subsequence of $(\BW_n)_{n \geq n_0}$ (still denoted $(\BW_n)_{n \geq n_0}$), 
 and  a sequence $(z^n)_{n\geq n_0} \subset \R^3$ such that 
$$
 \BW_n ( \cdot - z^n ) \to \BW \quad \quad \quad {\rm in}  \quad \mathscr{Y}(\R^3) . 
$$
Moreover, there is $ \si > 0 $ such that $ z \longmapsto \BW(z, \frac{1}{\si } z_{\perp})$  is a ground state 
(with speed $1/(2 \cs ^2)$) of  (KP-I).  We will prove that $ \si = 1$.

Let $ x^n = \left( \frac{z_1^n}{\e _n}, \frac{ z_{\perp}^n}{\e_n ^2} \right).$
We denote $ \tilde{\BW}_n = \BW_n( \cdot - z^n)$, $ \tilde{A}_n = A_n (\cdot - z ^n)$, $ \tilde{ \vp}_n = \vp_n (\cdot - z^n)$, $\tilde{U}_n = U_n (\cdot -x^n).$ 
It is obvious that $\tilde{U}_n$  satisfies (TW$_{c_n}$) and all the previous estimates hold with 
$ \tilde{A}_n$,   $ \tilde{ \vp}_n $ and  $\tilde{U}_n$ instead of $ A_n$, $\vp_n$ and $U_n$, respectively. 

Since $  \tilde{\BW}_n  = \frac{1}{\cs } \p_{z_1} \tilde{\vp}_n$ and $ \tilde{\BW}_n  \to \BW$ in $\mathscr{Y}(\R^3)$, we have
\beq
\label{conv1}
\p_{z_1} \tilde{\vp}_n \to \cs \BW, \qquad \quad 
\p_{z_1}^2  \tilde{\vp}_n \to \cs \p_{z_1}  \BW \qquad \mbox{ and } \qquad
\nabla_{z _{\perp}} \tilde{\vp}_n \to \cs \nabla_{z _{\perp}} \p_{z_1}^{-1} \BW 
\qquad \mbox{ in } L^2(\R^3). 
\eeq
Integrating by parts, then using the Cauchy-Schwarz inequality, Proposition \ref{Born} and \eqref{approx1} we find 
$$
\begin{array}{l}
\ds \int_{\R^3} \Big| \p _{z_1}^2 \tilde{\vp} _n - c_n \p_{z_1} \tilde{A}_n \Big| ^2 \, dz 
= - \int_{\R^3} ( \p _{z_1} \tilde{\vp} _n - c_n  \tilde{A}_n ) (  \p _{z_1}^3 \tilde{\vp} _n - c_n \p_{z_1}^2 \tilde{A}_n ) \, dz
\\
\\
\leq \| \p _{z_1} \tilde{\vp} _n - c_n  \tilde{A}_n  \|_{L^2} \|\p _{z_1}^3 \tilde{\vp} _n - c_n \p_{z_1}^2 \tilde{A}_n  \|_{L^2} 
= \BO(\e_n ^2), 
\end{array}
$$
hence $\| \p _{z_1}^2 \tilde{\vp} _n - c_n \p_{z_1} \tilde{A}_n  \|_{L^2} = \BO(\e _n) \to 0$. Since $ c_n \to \cs$, 
from  \eqref{approx1} and \eqref{conv1} we get 
\beq
\label{conv2} 
\tilde{A}_n \to \BW \qquad \mbox{ and } \qquad \p_{z_1} \tilde{A}_n \to \p_{z_1} \BW  
\qquad \mbox{ in  } L^2( \R^3) \quad \mbox{ as } n \to \infty.
\eeq
It is obvious that $ \tilde{A}_n$, $\tilde{\vp}_n $ and $ \e _n$ satisfy \eqref{desing}. 
Let $ \psi \in C_c^{\infty}(\R^3)$. 
We multiply  \eqref{desing} by $ \psi $, integrate by parts, then pass to the limit as $ n \to \infty$. 
We use Proposition \ref{Born}, \eqref{conv1} and \eqref{conv2} and after a straightforward computation we discover that $ \BW$ satisfies the equation (SW) 
in $\Do '( \R^3 )$. 
This implies that necessarily $ \si = 1$ and $ \BW $ is a ground state of speed $1/(2 \cs ^2)$ to (KP-I).
In particular, $ \BW $ satisfies  the Pohozaev identities \eqref{identites} and \eqref{Ident}. 

Since $ \tilde{\BW}_n \to \BW$ in  $ \mathscr{Y}(\R^3)  $, 
we have $ \bS ( \BW_n) = \bS ( \tilde{\BW}_n) \to \bS ( \BW) $ and \eqref{ecun} implies 
$$
\frac{ E_{c( \e _n)} ( U_n) }{\cs ^2 r_0 ^2 \e _n } = \bS ( \BW _n) + \BO( \e_n ^2) 
= \bS( \BW) + o(1) = \bS_{\rm min} +{o}(1), 
$$
that is \eqref{Ec} holds. 
Using the expression for the momentum in \eqref{momentlift}, then \eqref{conv1}, \eqref{conv2}, Proposition \ref{Born}
 and the Pohozaev identities  \eqref{identites} and \eqref{Ident} we get 
$$
- \frac{ \e_n}{r_0 ^2 \cs ^3 } Q(U_n)  = \frac{ \e_n}{r_0 ^2 \cs ^3 } \int_{\R^3} ( \rho _n ^2 - r_0 ^2) \frac { \p \phi_n}{\p x_1 } \, dx 
= \frac{1}{\cs ^3} \int_{\R^3} ( 2 A_n (z) + \e _n ^2 A_n ^2 (z) ) \frac{ \p  \vp _n}{ \p z_1} (z) \, dz 
\lra \frac{2}{\cs ^2} \int_{\R^3} \BW^2 (z) \, dz = \bS( \BW) .
$$
Hence $ - \cs Q(U_n) \sim r_0 ^2 \cs ^4 \bS_{\rm min} \e^{-1}  $ as $ n \to \ii$. 
Together with \eqref{Ec} this implies that $(U_n)_{n\geq n_0 } $   satisfies  \eqref{energy}. 

\medskip

By Proposition \ref{Born} we know that $ (\tilde{A}_n)_{n \geq n_0}$, $ (\p_{z_1} \tilde{A}_n)_{n \geq n_0}$, 
$ (\p_{z_1} \tilde{\vp }_n)_{n \geq n_0}$ and $ (\p_{z_1} ^2 \tilde{\vp }_n)_{n \geq n_0}$
are bounded in $L^p(\R^3) $ for $1 < p < \infty$. 
From \eqref{conv1}, \eqref{conv2} and standard interpolation in $L^p $ spaces we find  as $ n \to \ii$
\beq
\label{conv3}
\tilde{A}_n \to \BW,  \qquad
\p_{z_1} \tilde{A}_n \to \p_{z_1} \BW,  \qquad
\p_{z_1} \tilde{\vp }_n \to \cs \BW  \quad \mbox{ and } \quad 
\p_{z_1}^2  \tilde{\vp }_n \to \cs \p_{z_1} \BW 
\quad \mbox{ in } L^p 
\eeq
for any $ p \in (1, \ii)$.

Proceeding as in \cite{BGS1} (see Lemma 4.6 p. 262 and Proposition 6.1 p. 266 there)
one can prove that for any multiindex $\al \in \N^N$ with $ |\al | \leq 2$, the sequences  
$ (\p ^{\al } \tilde{A}_n)_{n \geq n_0}$, $ (\p ^{\al } \p_{z_1} \tilde{A}_n)_{n \geq n_0}$, 
$ ( \p ^{\al }  \p_{z_1} \tilde{\vp }_n)_{n \geq n_0}$ and $ (\p ^{\al }  \p_{z_1} ^2 \tilde{\vp }_n)_{n \geq n_0}$
are bounded in $L^p(\R^3) $ for $1 < p < \infty$. 
Then by interpolation we see that \eqref{conv3} holds in $W^{1,p}(\R^3) $ for all $ p \in (1, \infty)$. 
\carre

\subsection{Proof of Theorem \ref{res1} completed in the case $\bs{N=2}$}

Assume that $ N =2$. Let $ (U_n, c_n )$ be a sequence of travelling waves to (NLS) 
satisfying assumption (b) in Theorem \ref{res1} such that $ c_n \to \cs$ as $ n \to \ii$. Let  $ \e _n = \sqrt{\cs ^2 - c_n ^2}$.
By Theorem \ref{th2dposit} we have $ \ds \int_{\R^2} | \nabla U_n |^2 \, dx \to 0 $ as $ n \to \ii$ and then 
Lemma \ref{liftingfacile} implies that $\| \, | U_n | - r_0 \|_{L^{\ii}} \to 0 $; 
in particular, for $n$  sufficiently large we have a lifting 
$U_n (x) = \rho_n(x) \ex^{i \phi _n (x)} = r_0 \Big( 1 + \e_n^2 A_n (z) \Big) 
\ex ^{ i\e _n \vp_n (z) } $ 
as in \eqref{ansatzKP} and the conclusion of Proposition \ref{Born} holds for $A_n$ and $ \vp _n$. 
As in the proof of Proposition \ref{convergence} we obtain
\beq
\label{approx2}
\| \p _{z_1} \vp _{n } - c _n A_{n } \|_{L^2} = \BO(\e_n ^2) 
\qquad \mbox{ and } \qquad 
\| \p _{z_1} ^2 \vp _{n } - c _n \p_{z_1} A_{n } \|_{L^2} = \BO(\e_n ) \qquad \mbox{ as } n \to \ii. 
\eeq

Let $ k_n = \ds \int_{\R^2} |\nabla U_n(x)|^2 \, dx$. We denote 
$\BW_n = \cs^{-1} \p_{z_1} \vp_n $. By \eqref{approx2} we have  $\| \BW_n -  A_n \|_{L^2} = \BO(\e_n^2)$. 
As in the proof of Proposition \ref{convergence} we find
$ \ds \Big|  \int_{\R^2} (\p_{z_1} A_n )^2 - (\p_{z_1} \BW _n )^2 \, dz \Big|
= \Big\vert \int_{\R^2} ( \p _{z_1} A _n ) ^2 -  \frac{ ( \p_{z_1} ^2\vp _n )^2}{\cs ^2}  \, dz   \Big\vert = \BO( \e _n ^2).$
 Using \eqref{approx2} and Proposition  \ref{Born} we get
\begin{align}
\label{kn}
k_n = & \ \int_{\R^2} |\nabla U_n|^2 \ dx = 
\e_n r_0^2 \int_{\R^2} (\p_{z_1} \vp_n)^2 (1 + \e_n^2 A_n)^2 
+ \e_n^2 (\p_{z_1} A_n)^2 + \e_n^2 (\p_{z_2} \vp_n)^2 (1 + \e_n^2 A_n)^2 
+ \e_n^4 (\p_{z_2} A_n)^2 \ dz 
\nonumber
\\ 
= & \ \e_n r_0^2  \int_{\R^2} (\p _{z_1} \vp_n)^2 \ dz 
+ \e_n^3 r_0^2  \int_{\R^2} 
\Big( 2 A_n (\p_{z_1} \vp_n ) ^2 + (\p_{z_1} A_n )^2 
+ ( \p_{z_2} \vp _n )^2 \Big) \ dz 
+ \BO(\e_n^5) 
\nonumber
\\ 
 = & \ \e_n r_0^2 \cs^2 \int_{\R^2} \BW_n^2 \ dz 
+ \e_n^3 r_0^2 \cs^2 \int_{\R^2} \Big( 2 \BW_n^3 
+ \frac{1}{\cs^2} ( \p_{z_1}  \BW_n )^2 
+ ( \p_{z_2}\p_{z_1}^{-1} \BW_n )^2 \Big) \ dz 
+ \BO(\e_n^5) .
\end{align}
Inverting this expansion we find the following expression of $ \e _n$ in terms of $ k_n$:
\be
\label{pepsi}
\e_n = \frac{k_n}{r_0^2 \cs^2 \| \BW_n \|_{L^2}^2 } 
- \frac{k_n^3}{r_0^6 \cs^6 \| \BW_n \|_{L^2}^8} 
\int_{\R^2} \Big( 2 \BW_n^3 + \frac{1}{\cs^2} ( \p_{z_1}  \BW_n )^2 
+ ( \p_{z_2}\p_{z_1}^{-1} \BW_n )^2 \Big) \ dz + \BO(k_n^5) .
\ee

Recall that  the mapping $U_n (c_n \cdot )$ is a minimizer of the functional $I(\psi ) = Q( \psi)+ \ds \int_{\R^2} V( |\psi |^2) \, dx $ 
under the constraint $\ds \int_{\R^2} |\nabla \psi |^2 \, dx = k_n$. 
Using this information, Proposition \ref{asympto} $(i)$, the fact that 
$c_n^2 = \cs^2 - \e_n^2 $ and \eqref{pepsi} we get 
\begin{align}
\label{lapubelle}
c_n Q(U_n) + \int_{\R^2} V(|U_n|^2) \ dx = & \, 
c_n ^2 I(U_n (c_n \cdot ) ) = c_n^2 I_{\rm min} (k_n) 
\nonumber 
\\
\leq & \,  c_n^2 \left( - \frac{k_n}{\cs^2} - \frac{4 k_n^3}{27 r_0^4 \cs^{12} \bS^2_{\rm min}} + \BO(k_n^5) \right)
\nonumber 
\\
= & \, - k_n + \frac{k_n^3}{r_0^4 \cs^6 \| \BW_n  \|_{L^2}^4} 
- \frac{4 k_n^3}{27 r_0^4 \cs^{10} \bS_{\rm min}^2} + \BO(k_n^5) .
\end{align}
Moreover, using the Taylor expansion \eqref{V}, we find
$$ \int_{\R^2} V(|U_n|^2) \ dx = r_0^2 \cs^2 \e_n 
\int_{\R^2} \Big( A_n^2 + \e_n^2 \Big[ \frac{\G}{3} - 1 \Big] A_n^3 
+ \frac{V_4(\e_n^2 A_n)}{\e_n^4} \Big) \ dz $$
and by \eqref{momentlift} we  have
$$ Q(U_n) = - \e_n r_0^2 \int_{\R^2} \Big( 2 A_n + \e_n^2 A_n^2 \Big) 
\frac{\p \vp_n}{\p z_1} \ dz . $$
Taking into account \eqref{approx2} and the equality $c_n^2 = \cs^2 - \e_n^2$, then using expansion 
of $ \e _n$ in terms of $ k_n$ \eqref{pepsi}  we get
\begin{align}
\label{99}
c_n Q(U_n) + & \ \int_{\R^2} V(|U_n|^2) \ dx  
\nonumber 
\\
= & \ r_0^2 \cs^2 \left( \e_n  \int_{\R^2} 
\Big( - 2 A_n \BW_n + A_n^2 \Big) \ dz 
+ \e_n^3  \int_{\R^2} 
\Big( - A_n^2 \BW_n + \Big[ \frac{\G}{3} - 1 \Big] A_n^3 
+\frac{1}{\cs^2} A_n \BW_n \Big) \ dz 
+ \BO(\e_n^5) \right) 
\nonumber 
\\
= & \ r_0^2 \cs^2 \left( \e_n  \n \BW_n - A_n \n_{L^2}^2 
- \e_n  \int_{\R^2} \BW_n^2 \ dz 
+ \e_n^3  \int_{\R^2} 
\Big[ \frac{\G}{3} - 2 \Big] \BW_n^3 + \frac{\BW_n^2}{\cs^2} \ dz 
+ \BO(\e_n^5) \right) 
\nonumber 
\\
= & \ r_0^2 \cs^2 \left(  - \e_n  \int_{\R^2} \BW_n^2 \ dz 
+ \e_n^3  \int_{\R^2} 
\Big[ \frac{\G}{3} - 2 \Big] \BW_n^3 + \frac{\BW_n^2}{\cs^2} \ dz 
+ \BO(\e_n^5) \right)
\nonumber
\\
= & \ - k_n 
+ \frac{k_n^3}{ r_0^4 \cs^4 \| \BW_n \|_{L^2}^6} \bS(\BW_n) 
+ \BO(k_n^5) .
\end{align}
Inserting \eqref{99} into \eqref{lapubelle} we discover 
$$ \frac{k_n^3}{ r_0^4 \cs^4 \| \BW_n \|_{L^2}^6} \bS(\BW_n) 
+ \BO(k_n^5) \leq\frac{k_n^3}{r_0^4 \cs^6 \| \BW_n \|_{L^2}^4} 
- \frac{4 k_n^3}{27 r_0^4 \cs^{10} \bS_{\rm min}^2} + \BO(k_n^5) , $$
that is
$$ \bS(\BW_n) \leq \frac{1}{\cs^2} \| \BW_n \|_{L^2}^2 
- \frac{4}{27 \cs^{6} \bS_{\rm min}^2} \| \BW_n \|_{L^2}^6 
+ \BO(k_n^2) $$
or equivalently 
\be
\label{topbelle}
 \mathscr{E} (\BW_n) = \bS(\BW_n) - \frac{1}{\cs^2} \int_{\R^2} \BW_n^2 \ d z 
 \leq - \frac{1}{2 \bS_{\rm min}^2} 
\Big( \frac{2}{3} \Big)^3 
\cdot \Big( \frac{1}{\cs^{2}} \n \BW_n \n _{L^2}^2 \Big)^3 + \BO(k_n^2) .
\ee
As in the proof of Proposition \ref{convergence}, it follows from Lemma \ref{minoinf} and Proposition \ref{Born} that there are 
 some positive constants $ m_1, \, m_2 $ such that 
$$m_1 \leq \| \BW_n \|_{L^2}^2 \leq  m_2  \qquad \mbox{ for all sufficiently large } n. $$
Denote $  \la_n =  \frac{\| \BW_n \|_{L^2}^2}{\cs^2} .$
Passing to a subsequence if necessary we may assume that  $ \la _n \to \la $, where $ \la \in (0,+\ii)$.
Let 
$$ {\BW}_n ^{\#}  (z) = \frac{\mu^2}{\la_n^2} 
\BW_n \Big( \frac{\mu}{\la_n} z_1,\frac{\mu^2}{\la_n^2} z_2 \Big),  $$
where $\mu$ is as in Theorem \ref{gs2d}.  Then ${ \BW}_n ^{\# }$  satisfies
$$ \int_{\R^2} \frac{1}{\cs^2} \, ({\BW}_n^{\#}) ^2 \ dz = 
\frac{\mu}{\la_n} \int_{\R^2} \frac{1}{\cs^2} \, {\BW}_n^2 \ dz 
= \mu \quad \quad \quad {\rm and} \quad \quad \quad 
\mathscr{E} ({\BW}_n^{\# }) = \frac{\mu^3}{\la_n^3} \mathscr{E} (\BW_n) . $$
Plugging this into \eqref{topbelle} and recalling that 
$ \mu = \frac32 \bS_{\rm min}$, we infer that 
$$ \mathscr{E} ({\BW}_n^{\#}) = 
\frac{\mu^3}{\la_n^3} \mathscr{E} (\BW_n) 
\leq - \frac{1}{2 \bS_{\rm min}^2} 
\Big( \frac{2\mu }{3} \Big)^3 + \BO(k_n^2) 
= - \frac{1}{2} \bS_{\rm min} + \BO(k_n^2) . $$
Therefore  $({\BW}_n^{\#})_{n \geq n_0} $ is a minimizing sequence for 
\eqref{minimiz}. 

\medskip

By Theorem \ref{gs2d} we infer that there exist   a subsequence of 
$ ({\BW}_n^{\#} )_ {n \geq n_0}$,  still denoted 
$ ({\BW}_n^{\#})_{n \geq n_0} $, a sequence $ (z^n )_{n \geq n_0} = ( z_1^n, z_2^n)_{n \geq n_0} \subset \R^2 $ and 
 a ground state $\BW$  (with speed $1/(2\cs^2)$) of (KP-I) such that 
 $ {\BW}_n^{\#} ( \cdot - z^n) \lra \BW$ 
 strongly in $\mathscr{Y}(\R^2) $
 as $ n \to \ii$. 
 
Let $ x^n = \left( \frac{ \mu}{\e _n \la _n } z_1^n, \frac{ \mu ^2 }{\e _n ^2 \la _n ^2 } z_2^n \right)$
and $ \tilde{U}_n = U( \cdot - x^n)$, 
$ \tilde{A}_n (z) = A_n \left( z_1 - \frac{ \mu}{\la _n } z_1 ^n, z_2 - \frac{ \mu ^2}{\la _n ^2} z_2 ^n \right)$, 
$ \tilde{ \vp }_n (z) = \vp_n \left( z_1 - \frac{ \mu}{\la _n } z_1 ^n, z_2 - \frac{ \mu ^2}{\la _n ^2} z_2 ^n \right)$, 
$ \tilde{ \BW }_n (z) = \BW_n \left( z_1 - \frac{ \mu}{\la _n } z_1 ^n, z_2 - \frac{ \mu ^2}{\la _n ^2} z_2 ^n \right)$. 
We denote $ \tilde{ \BW}(z)= \frac{ \la ^2}{\mu ^2} \BW ( \frac{ \la }{\mu} z_1, \frac{ \la ^2}{\mu ^2} z_2)$. 
It is obvious  that $ \tilde{U}_n (x) = r_0 \left( 1 + \e _n ^2\tilde{A}_n (z) \right) \ex ^{ i \e _n \tilde{\vp }_n (z )}$ 
is a solution to (TW$_{c_n}$) with the same properties as $ U_n$ and the functions 
$ \tilde{A}_n$, $ \tilde{\vp}_n$, $ \tilde{\BW }_n$ satisfy the same estimates as $A_n$,  $ \vp _n$ and $ \BW _n$, respectively. 
Moreover, we have $ \tilde{\BW }_n = \frac{1}{\cs } \p_{z_1} \tilde{\vp} _n$ 
and $  \tilde{\BW }_n \lra \tilde{ \BW}$ strongly in $\mathscr{Y}(\R^2)  $ as $ n \to \ii$.

It is clear that $ \tilde{A}_n$, $\tilde{\vp}_n $ and $ \e _n$ satisfy \eqref{desing}. 
For any fixed  $ \psi \in C_c^{\infty}(\R^3)$  
we mutiply  \eqref{desing} by $ \psi $, integrate by parts, then pass to the limit as $ n \to \infty$. 
Proceeding   as in the proof of Proposition \ref{convergence} we find that $ \tilde{\BW }$ satisfies  equation (SW) in $ \Do'( \R^2)$. 
We know that $ \BW$ also solves (SW) and  comparing the  equations for $ \BW $ and $ \tilde{\BW}$ we infer that 
$ \left( \frac{ \la ^3}{\mu ^3} - \frac{ \la ^5}{\mu ^5} \right) \p_{z_1} \BW = 0 $ in $ \R^2$ . 
Since $ \p_{z_1} \BW \not= 0$, $ \la > 0 $ and $ \mu > 0$, we have necessarily $ \la = \mu$, that is $ \tilde{\BW} = \BW$. 

In particular, we have $ \bS ( \BW _n) = \bS ( \tilde{\BW}_n) \lra \bS ( \BW) = \bS_{\rm min} $ as $ n \to \ii$. 
Since $ \ds \int_ {\R^2} |\nabla U_n |^2 \, dx = k_n$, using \eqref{99} and \eqref{kn} we get 
$$
E(U_n) + c_n Q( U_n) =  \frac{k_n^3}{ r_0^4 \cs^4 \| \BW_n \|_{L^2}^6} \bS(\BW_n) + \BO(k_n^5)
\sim \e_n ^3 r_0 ^2 \cs ^2 \bS_{\rm min} 
\qquad \mbox{ as } n \to \ii .
$$
Hence \eqref{Ec} holds. 
As in the proof of Proposition \ref{convergence} we have 
$$
\begin{array}{l}
\ds Q( U_n) = - \int_{\R^2} (\rho_n ^2 - r_0 ^2 ) \frac{ \p \phi}{\p x_1}  
= - r_0 ^2 \e _n  \int_{\R^2}   ( 2 A_n (z) + \e _n ^2 A_n ^2 (z) ) \frac{ \p  \vp _n}{ \p z_1} (z) \, dz 
\\
\\
\ds \sim -2  r_0 ^2 \cs  \e_n  \int_{\R^2} \BW^2 (z) \, dz = - 3 r_0 ^2 \cs ^3 \bS( \BW) \e _n. 
\end{array}
$$
The above computation and \eqref{Ec} imply \eqref{energy}. 

Finally, the convergence in \eqref{conv3} as well as the similar property in $W^{1,p} (\R^2)$ are proven  exactly 
as in the three dimensional case.
\carre

\section{The higher dimensional case}

\subsection{Proof of Proposition \ref{dim6}}

We argue by contradiction. Suppose that the assumptions of Proposition  \ref{dim6} hold and  there is  a sequence 
$(U_n)_{n \geq 1}  \subset  \BE $ of nonconstant solutions to  (TW$_{c_n}$) 
 such that  $E_{c_n}(U_n) \to 0$ as $n \to +\ii$. 
By Proposition \ref{lifting} $(ii)$ we have  $|U_n| \to r_0 > 0 $ uniformly in 
$\R^N$.  Hence for $n$ sufficiently large we have the lifting 
$ U_n(x) = \rho_n(x) \ex^{i\phi_n(x)} . $
We write 
$$ \BB_{n  } =  \frac{|U_n| }{r_0 } - 1 , 
\qquad \mbox{ so that } \qquad
\rho _n = r_0( 1 + \BB _n) \qquad \mbox{ and } \qquad \BB _n \to 0 
\quad \mbox{ as } n \to \ii.
$$
Recall that $U_n$ satisfies the Pohozaev identities \eqref{Pohozaev}. The identity $P_{c_n} (U_n) = 0 $ can be written as 
$$
\int_{\R^N} \Big| \frac{ \p U_n}{\p x_1} \Big|^2 + \frac{N-3}{N-1} |\nabla _{x_{\perp}} U_n |^2 \, dx 
+ c_n Q(U_n) + \int_{\R^N} V(|U_n| ^2) \, dx = 0 .
$$
Using the formula \eqref{momentlift} for $Q(U_n)$ and the Taylor expansion \eqref{V} for $V( r_0 ^2 ( 1 + \BB_n)^2 ) $ we get 
\begin{align*}
r_0 ^2 
\int_{\R^N} & \Big| \frac{ \p \BB _n}{\p x_1} \Big|^2 + ( 1 + \BB _n) ^2 \Big| \frac{ \p  \phi _n}{\p x_1} \Big|^2
+ \frac{N-3}{N-1} |\nabla _{x_{\perp}} \BB _n |^2 
+ \frac{N-3}{N-1} ( 1 + \BB_n )^2 |\nabla _{x_{\perp}} \phi _n |^2 
\\
 & - c_n ( 2 \BB_n + \BB _n ^2 )  \frac{ \p  \phi _n}{\p x_1} 
+ \cs ^2 \left( \BB_n ^2 + \Big( \frac{ \G }{3} -1 \Big) \BB_n ^3 + V_4 (\BB _n) \right) \, dx = 0, 
\end{align*}
where $V_4 ( \al ) = \BO( \al ^4) $ as $ \al \to 0$. 
After rearranging terms, the above equality yields 
\begin{align*}
& \int_{\R^N} 
( \p_{x_1} \phi_n - c_n \BB _n )^2 
+ (\p_{x_1} \BB _n )^2 
+ \frac{N-3}{N-1} |\nabla_{x_\perp} \phi_n|^2 ( 1 + \BB _n )^2 
+ \frac{N-3}{N-1} |\nabla_{x_\perp} \BB _n|^2 + \e_n^2 \BB _n ^2 \ dx 
\nonumber \\ & = - \int_{\R^6} 
(\p_{x_1} \phi_n)^2 ( 2 \BB _n  + \BB_n ^2 ) 
+ \cs^2 \Big( \frac{\G}{3} - 1 \Big) \BB_n^3 + \cs^2 V_4( \BB_n ) 
- c_n \BB_n ^2 \p_{x_1} \phi_n \ dx \nonumber \\ & 
 = - \Big[ \frac{\G}{3} \cs^2 - \e_n^2 \Big] \int_{\R^N} \BB_n ^3 \ dz 
- \cs^2 \int_{\R^N} V_4( \BB_n ) \ dx 
- \int_{\R^N} (\p_{x_1} \phi_n)^2 \BB_n ^2 \ dx \nonumber 
\\ & \quad \quad 
+ \int_{\R^N} \BB_n  \Big( ( \p_{x_1} \phi_n - c_n \BB_n )^2 
-3 c_n \BB_n  (\p_{x_1} \phi_n - c_n \BB_n )  \Big) \ dx 
\end{align*}
and this can be written as 
\beq
\label{Dev3}
\begin{array}{l}
\ds \int_{\R^N} 
( \p_{x_1} \phi_n - c_n \BB _n )^2 
+ (\p_{x_1} \BB _n )^2 
+ \frac{N-3}{N-1} |\nabla_{x_\perp} \phi_n|^2 ( 1 + \BB _n )^2 
+ \frac{N-3}{N-1} |\nabla_{x_\perp} \BB _n|^2 + \e_n^2 ( 1 - \BB_n ) \BB _n ^2 \ dx 
\\ \\
= \ds  -  \frac{\G}{3} \cs^2  \int_{\R^N} \BB_n ^3 \ dz 
- \cs^2 \int_{\R^N} V_4( \BB_n ) \ dx 
- \int_{\R^N} (\p_{x_1} \phi_n)^2 \BB_n ^2 \ dx 
\\ \\\quad \quad \ds 
+ \int_{\R^N} \BB_n  \Big( ( \p_{x_1} \phi_n - c_n \BB_n )^2 
-3 c_n \BB_n  (\p_{x_1} \phi_n - c_n \BB_n )  \Big) \ dx .
\end{array}
\eeq
For $n$ sufficiently large we have $ \frac 12 \BB_n \leq ( 1 - \BB_n ) \BB _n ^2 \leq \frac 32 \BB _n ^2$ 
and then all the terms in the left-hand side of \eqref{Dev3} are nonnegative.
We will find an  upper bound  for the right-hand side of \eqref{Dev3}. First we  notice 
that the third integral there is nonnegative. Since  $\BB _n \to 0$  in $L^\ii$ and 
$V_4(\a) = \BO(\a^4) $ as $\a \to 0$, we have
\beq
\label{good1}
\Big|   \cs^2 \int_{\R^N} V_4( \BB _n) \ dx  \Big|
\leq C \| \BB_n  \|_{L^4} ^4
\leq C \| \BB_n  \|_{L^\ii} \| \BB_n  \|_{L^3}^3 .
\eeq
Using the fact that $\| \BB_n  \|_{L^\ii} \leq 1/4$ for 
$n$ large enough and the inequality $2ab \leq a^2 + b^2$, we get
\beq
\label{good2}
  \int_{\R^N} \BB_n  \Big( ( \p_{x_1} \phi_n - c_n \BB_n )^2 
- 3 c_n \BB_n (\p_{x_1} \phi_n - c_n \BB_n ) \Big) \ dx \leq 
\frac12 \int_{\R^ N} ( \p_{x_1} \phi_n - c_n \BB_n )^2 \ dx 
+ 9 \cs^2 \int_{\R^ N} \BB_n ^4 \ dx . 
\eeq

It is easy to see that $ \BB _n \in H^1 ( \R^N)$ (see the Introduction of \cite{CM1}).
We recall  the critical Sobolev embedding: for any $ h \in H^1 ( \R^N) $ (with $N\geq 3$) there holds 
\beq
\label{Sobol}
 \| h \|_{L^{\frac{2N}{N-2}}} \leq C 
\| \p_{x_1} h \|_{L^2}^{\frac{1}{N}} 
\| \nabla_{x_\perp} h \|_{L^2}^{\frac{N-1}{N}} . 
\eeq

Assume first that $ N \geq 6$. Then $ 2^* = \frac{ 2N}{N-2} \leq 3$. 
Using the Sobolev embedding \eqref{Sobol} and the fact that $\| \BB _n \| _{L^{\ii}} $ is bounded we get 
\beq
\label{good3}
\| \BB _n \|_{L^3} ^3 \leq \| \BB _n \|_{L^{\ii }}^{ 3 - 2^*} \| \BB _n \|_{L^{2^*}}^{ 2^*} 
\leq C \| \p_{x _1} \BB _n \| _{ L^2 } ^{\frac{2^*}{N}} \| \nabla_{x_{\perp}} \BB _n \| _{L^2}^{ \frac{ 2^*(N-1)}{N}} .
\eeq
Using  the inequalities  $ \| \BB_n \|_{L^4}^4 \leq 
\| \BB_n \|_{L^\ii} \| \BB_n \|_{L^3}^3 $ and 
$1+\BB_n  \geq 1/2$ for $n$ large, we deduce from \eqref{Dev3} that
\beq
\label{sobolomenthe}
\int_{\R^N} 
( \p_{x_1} \phi_n - c_n \BB _n )^2 
+ (\p_{x_1} \BB _n )^2 
+ |\nabla_{x_\perp} \phi_n|^2 
+ |\nabla_{x_\perp} \BB _n |^2 + \e_n^2 \BB _n ^2 \ dx 
\leq C \n \BB _n \n _{L^3}^3 . 
\eeq
From \eqref{sobolomenthe} and \eqref{good3} we obtain
\be
\label{sobolomenthos}
\n \nabla_{x_\perp} \phi_n \n _{L^2}^2 +
\n \p_{x_1} \BB _n \n_{L^2}^2 
+ \n  \nabla_{x_\perp} \BB_n  \n _{L^2}^2 
\leq C \n  \BB_n  \n_{L^3}^3 \leq C 
\n  \p_{x_1} \BB_n \n _{L^2}^{\frac{2}{N-2}} 
\n  \nabla_{x_\perp} \BB_n  \n _{L^2}^{\frac{2N-2}{N-2}}.
\ee

Assume now that ($ N=4$ or $N=5$) and $ \G \neq 0$. From \eqref{Dev3}, \eqref{good1} and 
\eqref{good2} we get 
\beq
\label{sobolomenthe1}
\int_{\R^N} 
( \p_{x_1} \phi_n - c_n \BB _n )^2 
+ (\p_{x_1} \BB _n )^2 
+ |\nabla_{x_\perp} \phi_n|^2 
+ |\nabla_{x_\perp} \BB _n |^2 + \e_n^2 \BB _n ^2 \ dx 
\leq C \n \BB _n \n _{L^4}^4 . 
\eeq
We have $ 2^* = 4 $ if $ N =4$ and $ 2^* = \frac{10}{3} < 4$ if $ N=5$. 
By the Sobolev embedding we have  
\beq
\label{good4}
\| \BB _n \|_{L^4 }^4 \leq \| \BB_n  \|_{L^\ii}  ^{ 4 - 2^* } \| \BB_n  \|_{L^{2^*}}^{2^*} 
\leq C   \| \BB_n  \|_{L^{2^*}}^{2^*}  
\leq C   \| \p_{x _1} \BB _n \| _{ L^2 } ^{\frac{2^*}{N}} \| \nabla_{x_{\perp}} \BB _n \| _{L^2}^{ \frac{ 2^*(N-1)}{N}} .
\eeq
The two inequalities above give
\be
\label{sobolomenthos1}
\n \nabla_{x_\perp} \phi_n \n _{L^2}^2 +
\n \p_{x_1} \BB _n \n_{L^2}^2 
+ \n  \nabla_{x_\perp} \BB_n  \n _{L^2}^2 
\leq C \n  \BB_n  \n_{L^4}^4 \leq C 
\n  \p_{x_1} \BB_n  \n _{L^2}^{\frac{2}{N-2}} 
\n  \nabla_{x_\perp} \BB_n  \n _{L^2}^{\frac{2N-2}{N-2}}.
\ee
From either \eqref{sobolomenthos}  or \eqref{sobolomenthos1} we obtain 
$$
\n \p_{x_1} \BB _n \n_{L^2}^2  \leq C 
\n  \p_{x_1} \BB_n  \n _{L^2}^{\frac{2}{N-2}} 
\n  \nabla_{x_\perp} \BB_n  \n _{L^2}^{\frac{2N-2}{N-2}}, 
$$
which gives $ \n \p_{x_1} \BB _n \n_{L^2}^ {\frac{2N -6}{N-2}} \leq C 
\n  \nabla_{x_\perp} \BB_n  \n _{L^2}^{\frac{2N-2}{N-2}}$, or equivalently 
\beq
\label{good5}
\n \p_{x_1} \BB _n \n_{L^2} \leq C \n  \nabla_{x_\perp} \BB_n  \n _{L^2}^{\frac{N-1}{N-3}}.
\eeq
Now we plug \eqref{good5} into \eqref{sobolomenthe}  or \eqref{sobolomenthe1} to discover
$$
\n  \nabla_{x_{\perp}} \BB_n  \n _{L^2}^ 2 
\leq  C
\n  \p_{x_1} \BB_n  \n _{L^2}^{\frac{2}{N-2}} 
\n  \nabla_{x_{\perp}} \BB_n  \n _{L^2}^{\frac{2N-2}{N-2}}
\leq C \n  \nabla_{x_{\perp}} \BB_n  \n _{L^2}^{\frac{2(N-1)}{N-3}}.
$$
Since $ \frac{2(N-1)}{N-3} > 2$ we infer that there is a constant $ m > 0 $ such that 
$\n  \nabla_{x_{\perp}} \BB_n  \n _{L^2} \geq m$ for all sufficiently large $n$. 
On the other hand $U_n$ satisfies the Pohozaev identity  $ P_{c_n }(U_n) = 0$, hence for large $n$ we have 
$$
E_{c_n} (U_n) = \frac{2}{N-1} \int_{\R^N} |\nabla_{x_{\perp} } U_n |^2 \, dx 
\geq \frac{2}{N-1} r_0 ^2 \int_{\R^N} |\nabla_{x_{\perp} } \BB _n |^2 \, dx 
\geq \frac{1}{N-1} r_0 ^2 m^2. 
$$
This contradicts the assumption that $E_{c_n} (U_n) \to 0 $ as $ n \to \ii$. 
The proof of Proposition \ref{dim6} is complete.
\carre

\begin{rem} \rm We do not know  whether  $T_{c}$ 
tends to zero  or not as $c \to \cs $  if  $N=4$ or $N=5$  and $\G \neq 0$.
\end{rem}

\subsection{Proof of Proposition \ref{vanishing}}

Let $ N \geq 4 $ and let
$ (U_n,c_n)_{n \geq 1} $ be a sequence of   nonconstant, finite energy solutions  solution of 
(TW$_{c_n}$) such that  $ E_{c_n}(U_n) \to 0$. By Proposition 
\ref{lifting} $(ii)$  we have $|U_n| \to r_0 > 0 $ uniformly in $\R^N$, 
hence for $n$ sufficiently large  we  may write
$$ U_n(x) = \rho_n(x) \ex^{i\phi_n(x)} 
= r_0 \Big( 1 + \a_n A_n (z) \Big) 
\exp \Big( i \beta_n \vp_n (z) \Big) \quad \quad \quad \mbox{ where }
z_1 = \la_n x_1, \quad z_\perp = \s_n x_\perp , $$
 and $ \a_n = \frac{1}{r_0} \| \rho_n - r_0  \|_{ L^{\ii}}  \to 0 $. 
Using the Pohozaev identity $ P_{c_n} (U_n) = 0$ and \eqref{blancheneige}  we have
$$
\frac{2}{N-1} \int_{\R^N} |\nabla_{x _{\perp} } U_n (x) |^2 \, dx = E(U_n) + c_n Q(U_n) 
= \frac{2}{N} \int_{\R^N} |\nabla \rho _n|^2 \ dx.
$$
Since $U_n \in \BE $ and $ U_n$ 
 is not constant, we have $ \ds \int_{\R^N} |\nabla_{x _{\perp} } U_n (x) |^2 \, dx > 0$ and the above identity implies 
 that $ \rho _n $ is not constant. The equality 
 $ E( U_n) + c_n Q( U_n) = \ds \frac{2}{N} \int_{\R^N} |\nabla \rho _n|^2 \ dx$
can be written as 
$$
\left( 1 - \frac 2N \right) \int_{\R^N} |\nabla \rho _n|^2 \, dx
+ \int_{\R^N} \rho _n ^2 |\nabla \phi _n|^2 \, dx
+ c_n Q( U_n) + \int _{\R^N} V(\rho _n ^2) \, dx = 0.
$$
Since $ \rho _n \to r_0 $ uniformly in $ \R^N$ as $ n \to \ii$, for $ n$ large we have $ V( \rho _n ^2) \geq 0 $ and 
from the last identity we infer that $ \ds 0 > c_n Q( U_n) = \int_{\R^N} (r_0 ^2 - \rho _n ^2) \frac{ \p \phi}{\p x_1} \, dx $, 
which  implies $ \| \p _{x_1 } \phi _n \|_{L^2}  > 0$. 
We must have $ \| \nabla _{x_{\perp}} \phi _n \|_{L^2} > 0 $ (otherwise $\phi $ would depend only on $ x_1$, 
contradicting the fact that $ \ds \int_{\R^N} |\nabla \phi _n|^2 \, dx $ is finite). 

 The choice of $ \al _n$ implies $ \| A_n \| _{L^{\ii}} =1 $. Since  $A_n$, $\p_{z_1} \phi _n$ and 
$\nabla_{z_\perp} \phi_n $ are nonzero, by scaling it is easy to see that 
\be
\label{normal}
\n A_n \n_{L^2} = \n \p_{z_1} \vp_n \n_{ L^2} = 
\n \nabla_{z_\perp} \vp_n \n_{L^2} = 1 
\ee
if and only if
$$ \la_n \s_n^{N-1} = \frac{\n \, |U_n| - r_0 \n_{L^\ii}^2}{\n \, |U_n| - r_0 \n_{L^2}^2} , 
\quad \quad 
\la_n \beta_n = \n \p_{x_1} \phi_n \n_{L^2} 
\frac{\n \, |U_n| - r_0 \n_{L^\ii}}{\n \, |U_n| - r_0 \n_{L^2}} , 
\quad \quad \beta_n \s_n = \n \nabla_{x_\perp} \phi_n \n_{L^2} 
\frac{\n \, |U_n| - r_0 \n_{L^\ii}}{\n \, |U_n| - r_0 \n_{L^2}} . $$
Since $N\geq 3$, the above equalities allow to compute 
$\la_n$, $\beta _n$  and $ \s_n$. 
Hence the scaling parameters $(\a_n,\beta_n, \la_n,\s_n)$ are uniquely determined  if \eqref{normal} holds and $ \| A_n \| _{L^{\ii}} =1 $.

\medskip

The Pohozaev identity $ P_{c_n }(U_n) = 0 $ gives 
\begin{align}
\label{Dev2}
& \int_{\R^N} 
\la_n^2 \beta_n^2 (\p_{z_1} \vp_n)^2 \Big( 1 + \a_n A_n \Big)^2 
+ \a_n^2 \la_n^2 (\p_{z_1} A_n)^2 
\nonumber \\ & 
+ \frac{N-3}{N-1} \beta_n^2 \s_n^2 |\nabla_{z_\perp} \vp_n|^2 
\Big( 1 + \a_n A_n \Big)^2 
+ \frac{N-3}{N-1} \a_n^2 \s_n^2 |\nabla_{z_\perp} A_n|^2 
+ \frac{1}{r_0 ^2} V\Big(r_0^2(1 + \a_n A_n)^2 \Big)\ dz
\nonumber \\ & \hspace{1cm} 
= 2 c_n \int_{\R^N} 2 \la_n \a_n \beta_n A_n \p_{z_1} \vp_n 
+ \la_n \a_n^2 \beta_n A_n^2 \p_{z_1} \vp_n \ dz .
\end{align}
By \eqref{normal}, the right-hand side of \eqref{Dev2} is 
$\BO(\la_n \a_n \beta_n)$. Since $\a_n \to 0$ and 
$\| A_n \|_{L^\ii}=1$ for $n$ large enough we have  $1+\a_n A_n \geq 1/2$, 
and by \eqref{V} we get $V(r_0^2(1 + \a_n A_n)^2 )  \geq \frac 12 r_0 ^2 \cs ^2 \al _n ^2 A_n ^2$. 
If  $N\geq 3$ 
all the terms in the left-hand side of \eqref{Dev2} are non-negative
and we infer that
$$
 \int_{\R^N} \la_n^2 \beta_n^2 (\p_{z_1} \vp_n)^2 
+ \a_n^2 A_n^2 \ dz = \BO(\la_n \a_n \beta_n). 
$$
From the normalization \eqref{normal} it follows that
$$ \la_n^2 \beta_n^2 = \BO(\la_n \a_n \beta_n), 
\quad \quad \quad {\rm and} \quad \quad \quad \a_n^2 = \BO(\la_n \a_n \beta_n), $$
which yields 
\be
\label{tutu1}
C_1 \leq \frac{\la_n \beta_n}{\a_n} \leq C_2 \qquad \mbox{ for some } C_1, \, C_2 > 0.
\ee
Let $\theta _n = \frac{ \la _n \beta _n}{\a _n}$.
We  use the Taylor expansion \eqref{V} for the potential $V$,  multiply 
\eqref{Dev2} by  $\frac{1}{\a_n^2}$ and write the resulting equality in the form
\begin{align*}
& \int_{\R^N} 
\Big( \theta_n \p_{z_1} \vp_n - c_n A_n \Big)^2 
+ \la_n^2 (\p_{z_1} A_n)^2 
+ \frac{N-3}{N-1} \frac{\theta_n^2 \s_n^2}{\la_n^2} |\nabla_{z_\perp} \vp_n|^2 
\Big( 1 + \a_n A_n \Big)^2 
+ \frac{N-3}{N-1} \s_n^2 |\nabla_{z_\perp} A_n|^2 
\nonumber \\ & \hspace{2cm} 
+ (\cs^2-c_n^2) A_n^2 \ dz 
\nonumber \\ & = - \int_{\R^N} 
\theta_n^2 \a_n (\p_{z_1} \vp_n)^2 \Big( 2 A_n + \a_n A_n^2 \Big) 
+  \cs^2 \a_n \Big( \frac{\G}{3} - 1 \Big) A_n^3 
+ \cs^2 \frac{V_4( \a_n A_n)}{\a_n^2} 
- 2 c_n \theta_n \a_n A_n^2 \p_{z_1} \vp_n \ dz .
\end{align*}
By \eqref{normal} and \eqref{tutu1} the right-hand side of the above equality is $\BO(\a_n)$. 
If $N \geq 3$  all the terms in the left-hand side are nonnegative.
In particular, we get 
$ \ds (\cs^2 - c_n^2 ) \int_{\R^N}  A_n ^2 \, dz 
= \cs^2 - c_n^2 = \BO(\a_n) , $
so that $c_n \to \cs$. Assuming that $N \geq 4$, we also infer that 
$$ \int_{\R^N} \la_n^2 (\p_{z_1} A_n)^2 
+ \frac{\s_n^2}{\la_n^2} |\nabla_{z_\perp} \vp_n|^2 \ dz = \BO(\a_n ) . $$
Together with \eqref{normal} and \eqref{tutu1}, this implies 
\be
\label{tutu2}
\frac{\s_n^2}{\la_n^2} = \BO(\a_n) 
\quad \quad \quad {\rm and} \quad \quad \quad 
\int_{\R^N} (\p_{z_1} A_n)^2 \ dz = \BO \Big( \frac{\a_n}{\la_n^{2}} \Big) .
\ee
The Pohozaev identity $P_{c_n}(U_n) = 0$ and \eqref{normal} imply that for each $n$ such that 
$ 1 + \a_n A_n \geq \frac 12 $ we have
\begin{align}
\label{final1}
E_{c_n}(U_n) = & \ \frac{2}{N-1} 
\int_{\R^N} |\nabla_\perp U_n|^2 \ dx \nonumber\\ \nonumber \\
= & \ \frac{2 r_0 ^2 }{(N-1)\la_n \s_n^{N-1}} \int_{\R^N} 
\beta_n^2 \s_n^2 |\nabla_{z_\perp} \vp_n|^2 
\Big(1+\a_n A_n \Big)^2 + \a_n^2 \s_n^2 |\nabla_{z_\perp} A_n|^2 \ dz \nonumber \\ \nonumber \\
\geq & \ \frac{r_0 ^2 \a_n^2 \theta _n^2 }{2(N-1) \la_n^3 \s_n^{N-3}} \int_{\R^N} 
|\nabla_{z_\perp} \vp_n|^2 \ dz \geq \frac{\a_n^2}{C \la_n^3 \s_n^{N-3}}.
\end{align}
However, in view of \eqref{tutu2}  we have
\beq
\label{final2}
 \frac{\a_n^2}{\la_n^3 \s_n^{N-3}} 
= \frac{\a_n^2}{\la_n^N ( \s_n / \la_n )^{N-3}} 
\geq \Big( \frac{\a_n}{\la_n^2} \Big)^{N/2} \frac{\a_n^2}{C \a_n^{N/2} \a_n^{(N-3)/2}} 
= \Big( \frac{\a_n}{\la_n^2} \Big)^{N/2} \frac{1}{C \a_n^{(2N-7)/2}} . 
\eeq
Notice that $ \a_n^{(2N-7)/2} \to 0 $ as $\a_n \to 0$ because  $N \geq 4$.
The fact that  $E_{c_n} (U_n) \lra 0 $, \eqref{final1} and \eqref{final2}  imply that 
$ \frac{\a_n}{\la_n^2} \to 0 $ 
{\rm as} $\ n \to +\ii . $
Then  using  \eqref{tutu2} we find
$$ \int_{\R^N} (\p_{z_1} A_n)^2 \ dz = \BO \Big( \frac{\a_n}{\la_n^{2}} \Big) 
\to 0 $$
and the proof is complete. 
\carre

\bigskip

\noindent 
{\bf Acknowledgement:} We greatfully acknowledge the support of the French ANR (Agence Nationale de la Recherche) 
under Grant ANR JC { ArDyPitEq}.


\end{document}